\documentclass[a4paper,10pt]{article}

\usepackage[hidelinks]{hyperref}
\usepackage[dvipsnames]{xcolor}
\usepackage{tikz}
\usepackage{bm,bbm}
\usepackage{natbib}
\usepackage{amsmath}
\usepackage{amssymb}
\usepackage{amsfonts}
\usepackage{amsthm}
\usepackage{mathrsfs}
\usepackage[title]{appendix}
\usepackage{textcomp}
\usepackage{manyfoot}
\usepackage{booktabs}
\usepackage{listings}
\usepackage{mathtools}
\usepackage{color}
\usepackage{dsfont}
\usepackage{animate}
\usepackage{etoolbox}
\usepackage{comment}
\usepackage{pgf}
\usetikzlibrary{calc}
\usetikzlibrary{patterns}
\usetikzlibrary{arrows}
\usetikzlibrary{decorations.pathreplacing}
\usetikzlibrary{pgfplots.groupplots}
\usepackage[utf8]{inputenc}
\usepackage{pgfplots}
\makeatletter
\pgfplots@warn@for@filter@discardsfalse
\makeatother
\usepackage{subcaption}

\usepackage{tcolorbox}
\usepackage{multirow}
\usepackage[ruled,vlined,linesnumbered,algosection,resetcount]{algorithm2e}
\usepackage{blkarray}
\usepackage{arydshln}

\usepackage{graphicx}
\usepackage{placeins}
\pgfplotsset{plot coordinates/math parser=false}
\usepgfplotslibrary{statistics}
\usepgfplotslibrary{fillbetween}
\pgfplotsset{select coords between index/.style 2 args={
    x filter/.code={
        \ifnum\coordindex<#1\fi
        \ifnum\coordindex>#2\fi
    }
}}

\usepackage{readarray}
\usepackage[margin=2cm]{geometry}

\theoremstyle{thmstyleone}
\newtheorem{theorem}{Theorem}
\newtheorem{proposition}[theorem]{Proposition}

\theoremstyle{thmstyletwo}

\theoremstyle{thmstylethree}

\theoremstyle{plain}
\newtheorem{lemma}[theorem]{Lemma}

\theoremstyle{definition}

\theoremstyle{remark}
\def\R{\mathbb R}
\def\P{\mathbb P}

\def\E{\mathbb E}
\def\T{\mathbb T}

\def\00{\mathbf 0}

\def\dd{\,\mathrm{d}}

\def\SS{\mathcal S}
\def\bet{\begin{theorem}}
\def\ent{\end{theorem}}
\def\bec{\begin{corollary}}
\def\enc{\end{corollary}}
\def\bep{\begin{proof}}
\def\enp{\end{proof}}
\def\f{\frac}

\def\a{\alpha}
\def\g{\gamma}
\def\la{\lambda}

\def\su{\subseteq}
\def\ms{\mathsf}
\def\co{\colon}

\def\mc{ \mathcal}
\def\ff{\infty}
\def\PP{\mc P}

\def\one{\mathds1}

\def\d{\mathrm d}

\renewcommand\leq{\leqslant}
\def\r{\rho}
\renewcommand\geq{\geqslant}
\renewcommand\le{\leqslant}
\renewcommand\ge{\geqslant}

\def\bs{\boldsymbol}

\def\e{\varepsilon}
\def\bel{\begin{lemma}}
\def\enl{\end{lemma}}
\def\im{\item}
\def\been{\begin{enumerate}}
\def\enen{\end{enumerate}}

\def\sm{\setminus}

\def\k_d{\kappa_d}

\def\tff{\uparrow\infty}

\def\co{\colon}
\def\ti{\times}
\def\k{\kappa}

\def\bepr{\begin{proposition}}
\def\enpr{\end{proposition}}
\def\vp{\varphi}

\def\os{o}

\def\ni{\noindent}

\def\gt{G^{\ms{th}, \eta}}
\def\gth{G^{\ms{th}, \eta}}

\def\b{\beta}
\def\Var{\ms{Var}}

\def\sna{S_n^{\ge}}
\def\snb{S_n^{\le}}
\def\snc{S_n^{(1)}}
\def\snd{S_n^{(2)}}

\def\Cov{\ms{Cov}}
\def\De{\Delta}
\def\TT{\mc T}
\def\CC{\mc C}
\def\Nl{\mathbf N_{\ms{loc}}}
\def\dego{D_{\ms{out}}}
\def\degi{D_{\ms{in}}}
\def\w{\wedge}

\providecommand{\keywords}[1]
{
  \small	
  \textbf{\textit{Keywords: }} #1
}

\providecommand{\classification}[1]
{
  \small	
  \textbf{\textit{MSC Classification: }} #1
}

\title{On the topology of higher-order\\age-dependent random connection models}

\makeatletter
\renewcommand\@date{{
  \vspace{-\baselineskip}
  \large\centering
  \begin{tabular}{@{}c@{}}
    \hspace{0.7cm} Christian Hirsch\textsuperscript{1,2} \hspace{0.7cm} \\
    \texttt{hirsch@math.au.dk}
  \end{tabular}
  \qquad and \qquad
  \begin{tabular}{@{}c@{}}
    Peter Juhasz\textsuperscript{1} \\
    \texttt{peter.juhasz@math.au.dk}
  \end{tabular}

  \bigskip

  \textsuperscript{1}\small{Department of Mathematics, Aarhus University, Aarhus, Denmark}\par
  \textsuperscript{2}DIGIT Center, Aarhus University, Aarhus, Denmark

  \bigskip
}}
\makeatother

\begin{document}

\maketitle

\abstract{
  In this paper, we investigate the potential of the age-dependent random connection model (ADRCM) with the aim of representing higher-order networks.
  A key contribution of our work are probabilistic limit results in large domains.
  More precisely, we first prove that the higher-order degree distributions have a power-law tail.
  Second, we establish central limit theorems for the edge counts and Betti numbers of the ADRCM in the regime where the degree distribution is light tailed.
  Moreover, in the heavy-tailed regime, we prove that asymptotically, the recentered and suitably rescaled edge counts converge to a stable distribution.
  We also propose a modification of the ADRCM in the form of a thinning procedure that enables independent adjustment of the power-law exponents for vertex and edge degrees.
  To apply the derived theorems to finite networks, we conduct a simulation study illustrating that the power-law degree distribution exponents approach their theoretical limits for large networks.
  It also indicates that in the heavy-tailed regime, the limit distribution of the recentered and suitably rescaled Betti numbers is stable.
  We demonstrate the practical application of the theoretical results to real-world datasets by analyzing scientific collaboration networks based on data from arXiv.
}

\vspace{0.3cm}

\classification{60D05, 60G55, 60F05}

\pgfkeys{/pgf/number format/.cd,1000 sep={\,}}

\FloatBarrier
\section{Introduction}
\label{sec:int}

\begin{figure} [h] \centering \includegraphics[width=0.38\textwidth]{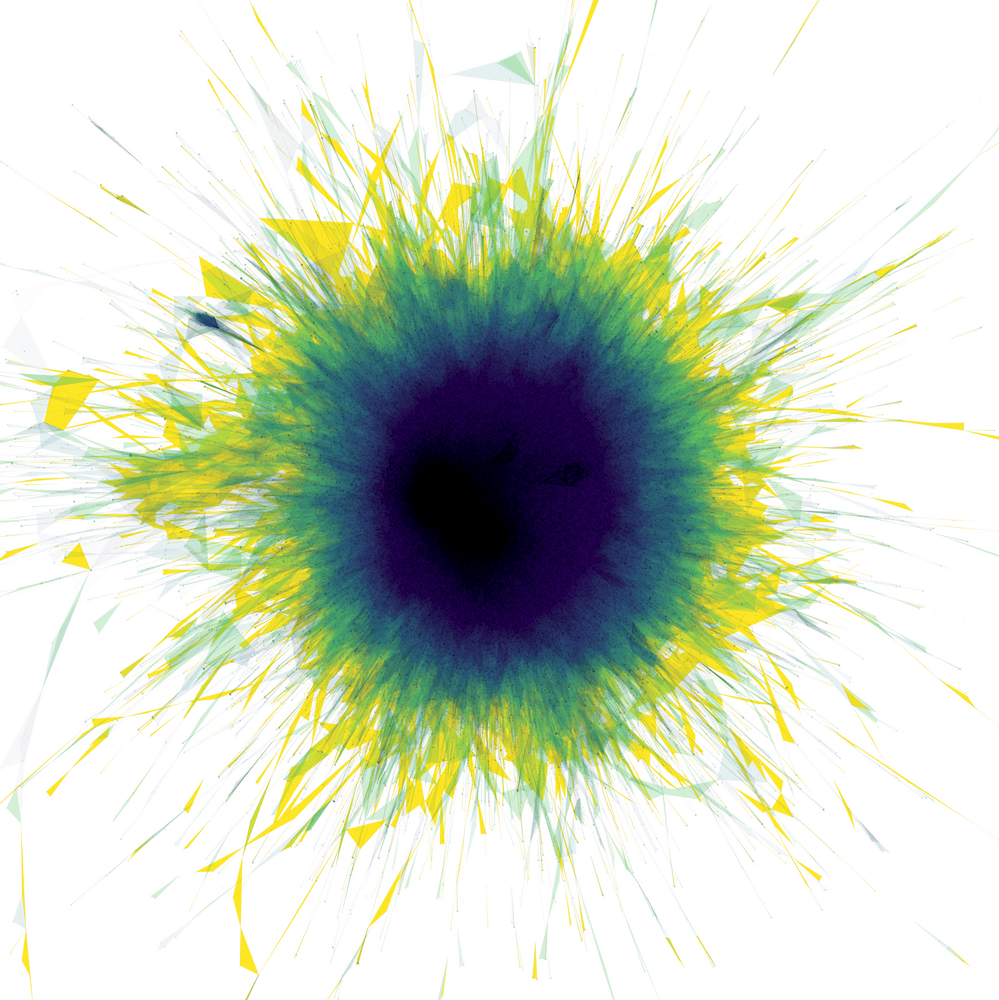} \caption{The largest component of a higher-order network of scientists.} \label{data_demo} \end{figure}

% complex networks
In recent decades, the field of complex networks has emerged as a powerful framework for analyzing systems whose properties cannot be understood by studying their parts in isolation.
The human brain, collaboration among researchers, the interaction of chemical elements, technological infrastructures, or the evolution of species are some examples for complex systems in which studying the relationships between the parts is inevitable~\citep{battiston2020networks,holland1976local}.
For instance, a collaboration network of scientists based on data from arxiv is illustrated in Figure~\ref{data_demo}, where vertices represent authors of publications and each document is represented by a simplex.

% importance of network modeling
Apart from a descriptive approach, it is often desirable to develop a stochastic model for generating synthetic networks.
The key advantage of creating such a stochastic model representation is that it enables a more refined analysis and a tool to deeper understand properties of a complex system.
Through this approach, it becomes feasible to reveal effects that might remain hidden in an actual dataset, particularly if its size is not large enough.
For an excellent discussion on complex network models as null models, we also refer the reader to~\cite{litvak}.

% simple networks
The traditional way of modeling complex systems relies on binary networks where the parts of the system are represented by vertices, and their relationships are represented by connections between them.
At the turn of the century, the field of complex networks experienced a rapid growth due to the insight that networks occurring in a wide variety of disciplines share a common set of key characteristics.

% preferential attachment
In their seminal work, \cite{barabasi} discovered that many key empirical features of complex networks are explained by a surprisingly simple \emph{preferential attachment model}.
Loosely speaking, it provides mathematical precision to the idea that many real networks emerge through a ``rich-get-richer'' mechanism.
In addition to the broad impact of network science in the application domains, complex networks also became the subject of intense research activity in mathematics, where rigorous mathematical proofs were provided of many of the effects that were previously empirically identified in network science~\citep{hf}.
In particular, the analysis of large preferential attachment models has become a highly fruitful research topic leading to important results such as the limiting distributions of various network characteristics~\citep{der1,der2}.

% spatial embedding
One of the shortcomings of the standard preferential attachment models is that they lead to tree-like structures, thus failing to reproduce the clustering property observed in real-world networks.
To address this issue, among others, spatial variants of the preferential attachment models have been proposed~\citep{jacMor1,jacMor2}.
Here, the network nodes are embedded in Euclidean space so that the preferential-attachment rule can take into account the spatial positions.
While the embedding produces the desired clustering effects, it makes the mathematical analysis more complicated.
Subsequently, it was realized by~\cite{glm2} that the decisive scale-free and clustering properties of spatial preferential attachment mechanism could also be realized by a simplified construction rule.
In the \emph{age-dependent random connection model (ADRCM)}, the connection probability to an existing vertex now depends on the age rather than the in-degree of that vertex.
In particular, knowing the age of a vertex does not require any information on the network structure in the neighborhood of that vertex.
This gives a far larger degree of spatial independence, which substantially simplifies many of the mathematical derivations.
Later,~\cite{komjathy,komjathy2} described a more general framework for incorporating weights into the connection mechanism.

% higher-order networks
As traditional network analysis was designed to study pairwise relationships between entities, simple network models are not capable of modeling higher-order interactions in which more than two entities are involved.
The study of higher-order network models has recently gained special attention due to its ability to capture these multibody relationships a simple network model cannot handle.
Among others, the need for higher-order relationships arise in scientific collaboration networks where the joint publication of three authors is not identical to three distinct papers with two authors each.
Beyond collaboration networks, the study of group relationships could already explain several phenomena like the synchronization of neurons  or the working mechanism of supply chain routes, see~\citep{xu2016representing}.

% modeling higher-order networks: simplicial-complexes
We model higher-order networks using simplicial complexes, where the relationships are represented with simplices of various dimensions.
The key benefit of modeling higher-order networks with simplicial complexes is that we can describe the networks using tools from topological data analysis (TDA).
This form of analysis was carried out in a number of studies~\citep{bianconi,carstens,patania,scolamiero}.

% contributions: model, theory, applications
While these studies investigate different datasets and rely on different TDA tools to analyze them, none of these works proposes a mathematical model to represent such higher-order networks.
Previously, \cite{fims} considered a stochastic model for higher-order complex networks.
However, as this model relies on a form preferential attachement mechanism, even the derivation of the asymptotic degree distribution is highly involved.
In contrast, since the ADRCM relies on a far simpler connection mechanism, in the present paper we are able to derive results that are substantially more refined than the degree distribution.
The main contributions of the present work are as follows:
\been
\im We begin by rigorously proving that the higher-order degree distributions follow a power law.
\im As a basis for hypothesis tests, we provide central limit theorems (CLTs) and stable limit theorems for the edge count and Betti numbers.
\im Recognizing the limitations of the ADRCM, we propose a model extension of the ADRCM capable of matching both any given admissible vertex and edge degree exponents.
\im Since these results are proved in the limit for large networks, we support the validity of these results for finite-size networks through conducting a simulation study.
\im Showing the convergence of the related quantities for finite-size networks, we proceed by developing statistical tests for finite networks based on the number of triangles and the Betti numbers for different parameter regimes.
\im Finally, we illustrate the use of these hypothesis tests for analyzing real-world collaboration networks.
\enen
We now expand on the above points in more detail and refer to Section~\ref{sec:mod} for details.

% ADRCM with simplices
As discussed earlier, the ADRCM stands out as an appealing model due to its ability to replicate key features -- power-law distributed vertex degrees and a high clustering coefficient -- observed in real-world complex networks, while also being mathematically tractable.
In light of this, our approach utilizes the ADRCM as a foundation and endows it with a higher-order structure by forming the \emph{clique complex}.
That is, the simplices in this complex are the cliques of the underlying graph.
A set of $k + 1$ vertices forms a $k$-simplex if and only if it is a $k$-clique, i.e., if and only if there is an edge between every pair of the $k + 1$ vertices.

% higher-order degree distributions
While for binary networks, the degree distribution is arguably the most fundamental characteristics, for higher-order networks, it is essential to understand also higher-order adjacencies.
Hence, to extend the concept of degree distributions to higher-order networks, we draw upon the concept of \emph{generalized degrees} introduced by~\cite{flavor}.
For $m' \ge m$, one considers the distribution of the number of $m'$-simplices containing a typical $m$-simplex as a face.
For instance, the standard vertex degree corresponds to the scenario where $(m, m') = (0, 1)$.
One of the fundamental findings by~\cite{glm2} is that in the ADRCM, the vertex-degree distribution satisfies a power law.
As another example, for $(m, m') = (1, 2)$ the \emph{edge degree} counts the number of triangles adjacent to a given edge.
In Theorem~\ref{thm:gen_deg}, we show that the generalized degrees also adhere to a power-law distribution.
Furthermore, we relate the exponents of the higher-order degree distribution to the exponent governing the vertex-degree distribution.
We pay special attention to the edge-degree distribution, since the formation of triangles in complex spatial networks is also of high interest for determining the clustering coefficient, as explored by~\cite{litvak,pvh}.

% edge count
In our second main result, we find that the distribution of the recentered and rescaled edge count in the ADRCM converges to a normal distribution for light-tailed degree distribution and to a stable distribution for heavy-tailed degree distributions (Theorems~\ref{thm:ass} and~\ref{thm:alpha}).
Based on our simulation study, we conjecture that these asymptotic results extend to higher-dimensional simplices.

% Betti numbers
Next, turning our focus to the features relevant in TDA, we continue with the analysis of the distribution of the Betti numbers of the clique complexes generated by the ADRCM.
\cite{siu2023many} derive asymptotic expressions for the growth rate of the expected Betti numbers in non-spatial preferential attachment models.
In contrast, the focus of our study is on the fluctuations around the expectation, enabling the application of hypothesis tests.
In Theorem~\ref{thm:bet}, we prove a CLT for the Betti numbers if the degree distribution is sufficiently light-tailed.
We also conjecture that for other values of the model parameters, the distribution of Betti numbers follows a stable distribution.
Again, this hypothesis gains credibility through the results of our simulation study, supporting the above conjecture.

% new model
By analyzing the empirical distributions within the arXiv data set, we find that the relation between the exponents governing vertex and edge degrees from Theorem~\ref{thm:gen_deg} to be too rigid to be applicable in real-world scenarios.
To address this limitation, we present a model extension that provides a larger flexibility to the original ADRCM by introducing a new parameter.
The main challenge of establishing this result is to ensure that we can independently adjust the edge-degree exponent, while keeping the vertex-degree exponent intact.
More precisely, we proceed as follows:
First, we increase both the vertex and the edge-degree exponents by adjusting the original parameters of the ADRCM, so that the edge-degree exponent reaches the desired value.
Then, we apply a dependent thinning operation involving the random removal of a fraction of certain edges that do not affect the edge-degree exponent, but which decrease the vertex-degree exponent.
These steps lead to the desired greater flexibility between vertex and edge degrees formalized in Theorem~\ref{thm:tgen_deg}.

% simulation
Our theoretical results presented above hold in the limit for very large networks.
For applications to real data, we accompany our theoretical results by a simulation study. 
\begin{itemize}
\item{
First, we explore the finite-size effects on higher-order degree distributions by examining the rate of convergence of the degree distribution exponents to their theoretical limits.
We see that the fluctuations of the exponents around their theoretical values decrease with increasing network size.
The simulations also reveal that apart from their fluctuations, the exponents also have a bias due to the finite size.
An interesting aspect of the simulation study is that, through Palm theory, we are able to simulate typical simplices in infinite networks that are free of finite-size effects.
}
\item{
Simulating three sets of networks with different model parameters, we validate our theoretical results regarding the edge-count distributions.
Furthermore, we also estimate the parameters of the distributions that are not explicitly derived in the theorems.
Finally, we discover the finite-size effects that are the most prominent in certain boundary cases.
}
\item{
As for the exploration of the edge-count distribution, we conduct a similar analysis for the Betti numbers.
This analysis supports our conjectures concerning the stable distribution of Betti numbers.
}
\end{itemize}

% data set analysis: analysis, fitting of parameters, hypothesis tests
Next, we demonstrate the application of the theorems on four real-world collaboration networks of scientists based on arXiv data.
After a general exploratory analysis, we analyze the vertex and edge-degree distribution exponents.

% triangle count, hypothesis tests
Based on the fitted vertex-degree exponents of the datasets, we fix the model parameters to use the ADRCM for further analysis of collaboration networks.
These fitted parameters guarantee that the vertex-degree exponent and the edge count are modeled correctly.
Thus, instead of the edge count, we conduct hypothesis tests based on the triangle count, where the null hypothesis is that the dataset is well described by the ADRCM.
Similar tests are also conducted for the Betti numbers.

The results of the hypothesis tests reveal that the topological structure of scientific collaboration networks is highly complex.
In particular, an elementary two-parameter model, such as the ADRCM, is not enough to capture all aspects of higher-order networks.

% Betti numbers, hypothesis tests

% organization
The rest of the manuscript is organized as follows.
Section~\ref{sec:mod} presents our main theoretical results regarding the higher-order networks generated by extending the ADRCM model to a clique complex.
Sections~\ref{sec:gen_deg},~\ref{sec:bet},~\ref{sec:alpha},~\ref{sec:tgen_deg} contains the proofs of the theorems stated in Section~\ref{sec:mod}.
Section~\ref{sec:sim} details the simulation study to demonstrate the validity of the asymptotic results discussed in Section~\ref{sec:mod} for finite networks.
Section~\ref{sec:dat} illustrates the application of the ADRCM model to higher-order networks of scientific collaborations.
Lastly, Section~\ref{sec:conc} includes a summary and ideas for directions of further research.

\section{Model and main results}
\label{sec:mod}

The higher-order network model discussed in this paper is an extension of the ADRCM, which we now recall for the convenience of the reader.
In this network model, vertices arrive according to a Poisson process and are placed uniformly at random in Euclidean space.
Two vertices are connected with a probability given by the profile function, which is a function of the distance and the ages of the vertices.

As justified below, we restrict our attention to the special case of latent Euclidean dimension is $d = 1$ and profile function $\vp(r) = \one_{[0, 1]}(r)$.
Let $\PP =\{ P_i \} = \{ (X_i, U_i) \}_{i \ge 1}$ be a unit-intensity Poisson point process on $\R \ti [0, 1]$, let $\beta > 0$ and $0 < \g < 1$.
Then, for $(x, u), (y, v)\in \PP$ with $u \le v$, there is an edge from $(y, v)$ to $(x, u)$, in symbols $(y, v) \to (x, u)$, if and only if
\[ |x - y| \le \frac{\b}{2} u^{-\g} v^{\g - 1}, \]
where $\b > 0$ is a parameter governing the edge density.
We henceforth denote this network by $G := G(\PP)$.

We stress that the framework developed by~\cite{glm2} allows to treat arbitrary dimensions and far more general connection functions.
However, the results by~\cite{glm2,pvh} indicate that many of the key network properties, such as the scaling of the vertex degree or the clustering coefficient, depend neither on the dimension nor the connection function.
We expect that a similar observation holds for higher-order characteristics and therefore decided to work with the simplest form of the ADRCM, greatly reducing the level of technicality in the presentation of the proofs.

While $G$ determines the binary vertex-interactions, in many applications, higher-order interactions play a crucial role.
The key idea for taking this observation into account is to extend $G$ to a simplicial complex.
The most popular approach for achieving this goal relies on the \emph{clique complex}~\citep{ctda}.
Here, a set of $k + 1$ vertices forms a $k$-simplex if and only if it is a $k$-clique, i.e., if and only if there is an edge between every pair of the $k + 1$ vertices.
To ease readability, we will henceforth also write $G = G(\PP)$ not only for the binary ADRCM network but also  for the clique complex generated by it.

While~\cite{glm2} analyze a number of key characteristics of the ADRCM considered as a binary network, we focus on the simplicial structure.
Specifically, we deal with the higher-order degrees and the Betti numbers, which we introduce in Sections~\ref{ssec:deg} and~\ref{ssec:test}, respectively.

%
% DEGREE DISTRIBUTIONS
%
\subsection{Higher-order degree distribution}
\label{ssec:deg}

% general degree
Arguably the most fundamental characteristics of complex networks is the degree distribution.
While the standard degree distribution provides an important summary of a complex network, it ignores higher-order structures.
Therefore,~\cite{courtney} argue to consider generalized degrees that are able to convey information on the adjacency structure of simplices of varying dimensions.

% vertex degree distribution
To define the typical vertex degree, the idea is to add to $\PP$ a distinguished typical vertex of the form $\os = (0, U)$ where $U$ is uniform in $[0, 1]$ and independent of $\PP$~\citep{glm2}.
We let $G_* = G(\PP \cup \{\os\})$ denote the ADRCM constructed on the extended vertex set.
Then, the typical vertex degree is that of $\os$ in $G_*$.
We define the tail of the vertex degree distribution $d_{0, 1} (k)$ to be
$$ d_{0, 1} (k) = \P \big( \deg_1 (\os) \ge k \big), $$
i.e., the probability that the vertex degree at the typical vertex exceeds $k \ge 0$.
In this context,~\cite{glm2} proved that the ADRCM is scale free in the sense that the degree distribution satisfies a power law:
$$ \lim_{k\tff} \log(d_{0, 1} (k)) / \log(k) = - \frac{1}{\g}. $$
While the higher-order vertex degrees provide a more refined picture than the standard vertex degrees, it is also important to go beyond vertices by considering higher-dimensional simplices as well.

To study generalized degrees, we define the higher-order degree of an $m$-simplex $\Delta \subseteq G$ as
$$ \deg_{m'} (\Delta) := | \{ \sigma \in G : \sigma \supseteq \Delta, \, |\sigma| = m' + 1 \} |, $$
represented as the number of $m'$-simplices containing $\Delta$.
For instance, $(m, m') = (0, 1)$ recovers the standard vertex degree and the higher-order vertex degree $\deg_{ m'} (v)$ of the vertex $v$ denotes for the number of $m'$-simplices that are incident to $v$.

% Palm representation
To study the generalized degree distributions, we consider typical simplices via the concept of Palm distribution.
Here, we describe the specific setting needed in the present paper, and refer the reader to~\cite{poisBook} for a more general introduction to Palm theory.
For $m \ge 0$, we can consider the $m$-simplices $\De_m = \{P_0, \dots, P_m\}$ in $G$ as a marked point process by centering $\De_m$ at its oldest vertex $c(\De_m)$.
Let $\TT_m(\PP)$ denote the family of $m$-simplices in the clique complex on $G$.
Then, the expectation of a function $f$ of the typical $m$-simplex $\De_m^*$ is given by
\begin{align} \label{eq:palm}
  \E [ f \left( \De_m^*, \, \PP \right) ] = \f1{\la_m} \E \Big[ \sum_{\De \in \TT_m(\PP)} \one \left\{ c \left( \De \right) \in \left[ 0, 1 \right] \right\} f \left( \De - c \left( \De \right), \, \PP - c \left( \De \right) \right) \Big],
\end{align}
where $\la_m > 0$ is the simplex density and where $f \co \CC_m \times \mathbf{N}_{\mathrm{loc}} \to [0, \ff)$ is any measurable function from the space $\CC_m$ of distinct $(m + 1)$-tuples of points in $\R\ti[0, 1]$ and the space of locally finite point processes to $[0, \ff)$, and which is symmetric in the first $m + 1$ arguments.

% generalized vertex degree distribution
In the present paper, we extend \cite[Proposition 4.1]{glm2} result by considering the generalized vertex degree distribution
$$ d_{m, m'} (k) = \P \big( \deg_{m'} (\De_m) \ge k \big), $$
represented as the distribution of the number of $m'$-simplices incident to $\os$.

\renewcommand{\theenumi}{\alph{enumi}}

%
% THM DEGREE DISTRBUTION
%
\bet[Power law for the typical vertex \& edge degree]
\label{thm:gen_deg}
Let $\g \in (0, 1)$ and $m' \ge m \ge 0$.
Then,
$$ \lim_{k\tff} \log(d_{m, m'} (k)) / \log(k)=  m - \frac{m + 1}{\g}. $$
\ent

% hypothesis tests
\subsection{Central and stable limit theorems}
\label{ssec:test}
As outlined in Section~\ref{sec:int}, to decide whether a given model is a good fit for a dataset, it is important to be able to carry out statistical hypothesis tests.
In this work, we discuss possible hypothesis tests that become asymptotically exact for growing networks.
While higher-order degree distributions are an important tool for describing higher-order networks, they only provide a highly restricted view of the topological structure.
The idea behind TDA is to rely on invariants from algebraic topology for extracting more refined shape-related information.
In this context, one of the most celebrated characteristics is the Betti numbers, which, loosely speaking, can be interpreted as the number of topological holes in a dataset.
For a more detailed explanation on Betti numbers, we refer the reader to~\cite{shirai}.
One attractive way to deveop a hypothesis test is to show that the considered test statistic becomes asymptotically normal.
This is the content of the following theorem.
Here, we write $\b_{n, q}$ for the $q$th Betti number of the clique complex $G\big(\PP \cap [0, n]\big)$.

%
% THM BETTI NUMBERS
%
\bet[CLT for the Betti numbers]
\label{thm:bet}
Let $q \ge 0$ and $0 < \g < 1/4$.
Then, $\Var(\b_{n, q})^{-1/2}(\b_{n, q} - \E[\b_{n, q}])$  converges in distribution to a standard normal distribution.
\ent

A disadvantage of Theorem~\ref{thm:bet} is that our proof imposes a substantial constraint on the parameter $\g$.
In particular, Theorem~\ref{thm:bet} considers a regime where the variance of the degree distribution is finite, while for many real-world datasets it is infinite.
Note that for large values $\g$, the ADRCM gives rise to extremely long edges, which makes it difficult to control spatial correlations over long distances, which is the main challenge in the proof.
While we expect that by a more careful argumentation in the proof of Theorem~\ref{thm:bet}, the range of $\g$ could be extended, we conjecture that the asymptotic normality breaks down for values of $\g > 1/2$.
To provide evidence for this conjecture, we now illustrate that a similar effect occurs for a more elementary test statistic, namely, the edge count
$$ S_n := \left| \left\{ (y, v) \to (x, u) : \, (y, v), \, (x, u) \in \PP, \, x \in [0, n] \right\} \right|.$$

%
% THM SIMPLICES
%
The key observation is that depending on whether $\g$ is smaller or larger than $1/2$, the variance of the $S_n$ at a typical vertex is either finite or infinite.
Hence, we should only expect a CLT in the finite variance regime.
We show that this is indeed the case.

%
% THM SIMPLEX COUNT CLT
%
\bet[CLT for the edge count]
\label{thm:ass}
Let $\g < 1/2$.
Then, $\Var(S_n)^{-1/2}(S_n - \E[S_n])$  converges in distribution to a standard normal distribution.
\ent

For $\g > 1/2$ since the degree distribution is heavy-tailed, the right tails in the edge count are more pronounced than those of a normal distribution.
For many combinatorially defined network models like the configuration model, the degrees are taken iid from a given distribution.
Hence, here the limiting vertex-degree distribution follows from the classical stable central limit theorem~\citep[Theorem 4.5.2]{whitt}.
We also refer the reader to~\cite{litvak} for a discussion in this direction.
However, we are not aware of any existing corresponding results for spatial networks, which often feature strong spatial correlations between the individual vertex degrees.
Hence, the main challenge in the proof of Theorem~\ref{thm:alpha} is to understand and overcome these correlations in order to extend the results from the combinatorial networks to spatial network models.

%
% THM EDGE COUNT SLT
%
\bet[Stable limit law for the edge count]
\label{thm:alpha}
Let $\g \in (1/2, 1)$.
Then, $n^{-\g}(S_n - \E[S_n])$ converges in distribution to $\SS$, where $\SS$ is a $1/\g$-stable distribution.
\ent

%
% MODEL EXTENSION, TADRCM
%
\subsection{Model extensions}
\label{ssec:ext}

Theorem~\ref{thm:gen_deg} expresses the power-law exponent of the vertex degree distribution and the edge degree distribution in terms of $\g$.
However, as we will illustrate in Section~\ref{sec:dat}, when analyzing datasets of scientific collaborations, the relation between the vertex and edge exponents suggested in Theorem~\ref{thm:gen_deg} may often be violated in real datasets.
More precisely, for a given vertex degree exponent, we found the edge degree in the data to be substantially more heavy-tailed than suggested in Theorem~\ref{thm:gen_deg}.
In other words, real datasets exhibit a larger proportion of edges incident to a large number of triangles than what can be realized by the ADRCM.
Alternatively, we could choose $\g$ so as to match the power-law exponent of the edge degree in the data.
In this case, however, the vertex degrees of the fitted model would exhibit too heavy tails.

To address this shortcoming, we propose a model extension \emph{thinned age-dependent random connection model (TADRCM)}, where we remove some edges so that the power-law exponent of the edge degrees is not affected.
The key observation is that for edges with high edge degrees, typically both endpoints are very old.
However, only a very small proportion of vertices connect to more than one very old vertex.
More precisely, we say that an edge $(z, w) \to (x, u)$ is \emph{protected} if $w \le 2u$ or if there exists a vertex $(y, v)$ with $v \le 2u \le 4v$ with $(z, w) \to (y, v)$.
An edge is \emph{exposed} if it is not protected.
Then, we define the TADRCM $\gth$ of $G$, by removing independently exposed edges.
The key idea is to use a retention probability of $u^{\eta}$, where $\eta > 0$ is a new model parameter.
Our next main result is the following analogue of Theorem~\ref{thm:gen_deg} for the thinned model, where $d^{\ms{th}, \eta}_{m, m'}$ is defined as $d_{m, m'}$, except that we use the TADRCM instead of the ADRCM.

%
% THM THINNED MODEL DEGREE DISTRIBUTION
%
\bet[Power law for the thinned typical vertex and edge degree]
\label{thm:tgen_deg}
Let $\g \in (1/2, 1)$ and $\eta > 0$ be such that $2/\g - 1 > 1/ (\g - \eta)$. Then,
$$ \lim_{k\tff} \log(d^{\ms{th}, \eta}_{0, m'} (k)) / \log(k)= - 1 / (\g - \eta) \quad \text{ and } \quad \lim_{k\tff} \log(d^{\ms{th}, \eta}_{1, m'} (k)) / \log(k)= 1 - 2/\g. $$
\ent

We stress that alternative approaches also exist to enhance the flexibility of the ADRCM.
For instance,~\cite{pvh} introduce a different extension, focusing on clustering properties.
However, in the scope of our work, we found the thinning-based model more convenient for two reasons.
First, through Theorem~\ref{thm:tgen_deg}, the parameters $\g$ and $\eta$ are related very transparently to the vertex and the edge degrees, which simplifies substantially fitting the model to datasets.
In contrast,~\cite{pvh} discuss a model where the connection between the model parameters and degree exponents is less obvious, and it is not immediate, which combination of higher-order degrees can be realized in the model.
Second, when carrying out the proofs, it is convenient that in the ADRCM the outdegree is Poisson-distributed independently of the vertex age.
Although we expect that our proofs could be adapted to the extension from~\cite{pvh} some of the steps would require more work.

\section{Proof of Theorem~\ref{thm:gen_deg} -- power-law exponents for higher-order simplex}
\label{sec:gen_deg}

In this section, we establish Theorem~\ref{thm:gen_deg}, i.e., we compute the power-law exponents for the higher-order simplex degrees in the ADRCM.
To reach this goal, we consider separately the lower and upper bounds in Sections~\ref{ss:gdl} and~\ref{ss:gdu}, respectively.

To prepare the proof, we start with an integral representation for the distribution of the typical $m$-simplex $\Delta_m^*$.
While~\eqref{eq:palm} provides a conceptually clean definition of the expectation of a function of a typical $m$-simplex, it is not ideal for carrying out actual computations.
For this reason, we derive an alternative representation in Proposition~\ref{pr:palm} below.

In this representation, we write $o := o_0 := (0, u)$ with $u \in [0, 1]$ for the typical vertex at the origin and
$$ \bs o_m := (o_1, \dots, o_m) := \big((y_1, v_1), \dots, (y_m, v_m)\big) \in \T^m := (\R\ti [0, 1])^m $$
for the remaining vertices.
Then, we let $g_m(u, \bs o_m)$ be the indicator of the event that $(o_0, \bs o_m)$ forms an $m$-simplex in the ADRCM with $u\le v_1\le \cdots \le v_m$.
Henceforth, we let $I_r (x) := [-r/2 + x, r/2 + x]$ denote the interval of side length $ r > 0$ centered at $x \in \R$.
We let $\Nl$ denote the family of all locally finite subsets of $\T$.

%
%PR PALM
%
\bepr[Distribution of the typical $m$-simplex]
\label{pr:palm}
Let $m \ge 1$.
Then,
$$\E [ f(\De_m^*, \mc P) ] = \frac {\int_0^1\int_{\T^m} \E[f(\{ o, \bs o_m\}, \mc P \cup \{ o, \bs o_m\}) ]  g_m (u, \bs o_m) \d \bs o_m\d u}{\int_0^1\int_{\T^m}  g_m (u, \bs o_m) \d \bs o_m \d u},
$$
for any measurable $f\co \CC_m \times \Nl \to [0,\ff)$, which is translation-covariant in the sense that $f((x + y, u), \vp + y) = f((x, u), \vp)$ for every $ (x, u) \in \T, y \in \R$ and $\vp \in \Nl$.
\enpr

%
%LEM GM
%
To ensure that the Palm version is well-defined, we need to show that the denominator is finite.
We formulate this property as a separate auxiliary result.
First, define the function
$$ \mu_m(u) :=\int_{\T^m} g_m(u, \bs o_m)  \d \bs o_m. $$
For instance,  $\mu_0 \equiv 1$ and also for $m = 1$ the expression simplifies.
To that end, we write
$$ M(p) := \{p' \in \T \co p' \to p\} $$
for the set of all space-time points connecting to $p \in \T$.
Then,
\begin{align} \label{eq:mu} \mu(u) := \mu_1(u) = |M(o)| = \int_u^1 |I_{\b u^{-\g}v^{\g - 1}} (0)| \d v = \f{\b}\g (u^{-\g} - 1) \end{align}
is the expected in-degree of the typical vertex. That is, $\mu_1(u) = \E[ D_{\ms{in}}(u)]$, where
$$ D_{\ms{in}}(u) := \big| \PP \cap M(o) \big| $$
is the in-degree of the typical vertex $o$.

For general $m \ge1$, we can derive the small-$u$ asymptotics.

%
%LEM GM
%
\bel[Asymptotics for $\mu_m(u)$] \label{lem:gm}
Let $m \ge 1$, $\g \in (0, 1)$ and $\eta > 0$.
Then, $\mu_m(u) \in O(u^{-\g - \eta})$.
\enl

%
%PRF GM
%
\bep
 First, 
$$\int_{\T} g_m(u, \bs o_m)\d  o_m  \le g_{m - 1}(u, \bs o_{m - 1})\int_0^1 \big|I_{\b v_{m - 1}^{-\g} v_m^{\g - 1}}(y_{m - 1})\big|\d v_m = \f{\b}\g g_{m - 1}(u, \bs o_{m - 1})v_{m - 1}^{-\g}.$$
Next, 
$$\int_{\T} g_{m - 1}(u, \bs o_{m - 1}) v_{m - 1}^{-\g - \eta} \d o_{m - 1} \le g_{m - 2}(u, \bs o_{m - 2}) \int_{v_{m - 2}}^1 \b v_{m - 2}^{-\g} v_{m - 1}^{-1 - \eta}\d v_{m - 1} \le \f{\b}\eta g_{m - 2}(u, \bs o_{m - 2})   v_{m - 2}^{-\g -\eta}.$$
Hence, iterating this bound yields that
$\mu_m(u) \le {\b^m}u^{-\g - \eta}/({\g \eta^{m - 1}}),$
as asserted.
\enp

\bep[Proof of Proposition~\ref{pr:palm}]
Let $g_m'( P_0, \dots, P_m )$ be the indicator of the event that $\{P_0, \dots, P_m\}$ forms an $m$-simplex in the ADRCM with $U_0 \le\cdots \le U_m$.
Let $A \su \R$ be a Borel set with $|A| = 1$.
Then
\[\la_m \E [ f( \De_m^*, \PP) ] = \E \Big[ \sum_{\substack{P_0, \dots, P_m \in \PP\text{ distinct} \\ U_0 \le \cdots \le U_m}}\one\{X_0 \in A\} f \big( \{ P_0, \dots, P_m \}, \PP \big) \, g_m'( P_0, \dots, P_m ) \Big], \]
Then, writing $\bs p = (p_1, \dots, p_m)$, by the Mecke formula~\citep[Theorem 4.4]{poisBook},
\begin{align*}
\la_m\E [ f( \De_m^*, \PP) ] = \int \limits_{A \ti [0, 1]} \int \limits_{\T^m}  \E\big[f ( \{ p_0, \bs p \}, \PP ) \big] &g_m'(\{ p_0, \bs p \}) \d \bs p \, \d p_0.
\end{align*}
As $|A| = 1$, a substitution $\mathbf{p}_m = \mathbf{o}_m + p_0$ and an application of Fubini's theorem give that
\[ \begin{split}
\la_m\E [ f( \De_m^*, \PP) ] = \int \limits_0^1 \int \limits_{\T^m} &\E\big[f ( \{ o, \bs o_m \}, \PP )\big]
g_m ( u, \bs o_m ) \, \d  \bs o_m  \, \d u.
\end {split} \]
Hence, evaluating this equality for $f = 1$ concludes the proof.
\enp

%
%SS GDL
%
\subsection {Proof of lower bound} \label{ss:gdl}

Next, we prove the lower bound by relying on the Palm representation derived in Proposition~\ref{pr:palm}.
More precisely, we produce specific configurations of $u, \bs o_m$ that occur with sufficiently high probability and such that $\P (\deg_{m'}(u, \bs o_m) \ge k)$ is bounded away from $0$.

%
%PRF LOW
%
\bep[Proof of Theorem~\ref{thm:gen_deg}, lower bound]
To ease notation, we put $\b' := \b/2$. 
First, let $p:= \P\big(\PP([0, \b'] \ti [3/4, 1]) \ge m'\big)$ denote the probability that a $(\b' \ti 0.25)$-box contains at least $m'$ Poisson points.
Furthermore, let $M:= \lceil 2/p\rceil$. Then, consider the set $B_k \su [0, 1]^{m + 1}\ti \R^m$ given by
$$ B_k := B_k '\ti [0, \b' k]^m := \big(\prod_{j \le m + 1} [( j/(Mmk))^{1/\g},((j + 1)/(Mmk))^{1/\g} ]\big) \ti [0, \b' k]^m. $$
Since $|B_k| \in \Omega( k^{-(m + 1)/\g + m})$, we only need to verify the following two items for every $(u, \bs o_m) \in B_k$.
\been
\im The $(o, \bs o_m)$ points form an m-simplex in the clique complex of the ADRCM.
\im It holds that $\P\big(\deg_{ m'}(u, \bs o_m) \ge k\big) \ge 1/2$
\enen

%1
For part (a), note that every $(u,\bs o_m) \in B_k$ indeed defines an $m$-simplex since
$$ \max_{i\le m}|y_i| \le \b' k \le \b' ((Mk)^{-1/\g})^{-\g}\;\text{ and }\max_{i, j\le m}|y_i - y_j| \le \b'k \le \b' ((Mk)^{-1/\g})^{-\g}. $$

%2
For part (b), we note that the events $E_{i, k}: =\big\{\PP\big([i\b', (i + 1)\b'] \ti [3/4, 1]\big) \ge m'\big\}$ are independent for $i \le kM$.
Moreover, let $N_k := \sum_{i \le kM}\one\{E_{i, k}\}$ be the number of the events that occur.
Then, $N_k$ is a binomial random variable with $kM$ trials and success probability $p$.
Since $kMp \ge 2k$, the binomial concentration result implies that $\P(N_k \ge k) \ge 1/2$ holds for sufficiently large $k$.

Hence, it suffices to show that almost surely, $N_k \le\deg_{m'}(u, \bs o_m)$.
To achieve this goal, we first note that for fixed $i \le kM$ any two points in $[i\b', (i + 1)\b'] \ti [0, 1]$ are connected by an edge.
Moreover, we claim that  any $(Z, W) \in [0, \b' kM] \ti [3/4, 1]$ connects to $o$ and to every $o_i $, $i \le m$.
Now,
$$ |Z - 0| \le \b' kM \le \b'( (kM)^{-1/\g})^{-\g}\;\text{ and }\;\max_{i \le m}|Z - y_i| \le \b' kM \le \b'( (kM)^{-1/\g})^{-\g}. $$
This concludes the proof  since the Poisson concentration inequality~\citep[Lemma 1.2]{penrose} implies that $\P(\PP(C_k) \ge k)\to1$  as $k \tff$.
\enp

%
%SS GDU
%
\subsection{Proof of upper bound} \label{ss:gdu}
In this subsection, we prove the upper bound for the simplex degree in Theorem~\ref{thm:gen_deg}.
First, to provide the reader with a gentle introduction, we present the case of the in-degree, which was considered previously by~\cite[Proposition 4.1]{glm2}.
In fact,~\cite[Lemma 4]{glm} is slightly more refined than Theorem~\ref{thm:gen_deg} in the sense that it provides not only the asymptotics for the tail probabilities but for the entire probability mass function.
Nevertheless, we include the short argument here because it makes the presentation self-contained and provides a intuition for the more complicated higher-order case.

The key observation is that conditioned on the arrival time $u$ of the typical vertex $o = (0, u)$, the in-degree is Poisson distributed.
Indeed, by the restriction theorem, the in-neighbors form a Poisson point process for fixed $u$~\citep[Theorem 5.2]{poisBook}.

\bep[Proof of upper bound for indegree]
First, note that if $\mu(u) \le k/2$ -- where $\mu(u)$ is the expected in-degree of the typical vertex introduced in~\eqref{eq:mu} --, then by the Poisson concentration inequality, the probability $\P ( D_{\ms{in}}(u) \ge k)$ vanishes exponentially fast in $k$.
Hence, we may assume that $u \le \mu^{-1}(k/2)$.
Noting that~\eqref{eq:mu} gives that $\mu^{-1}(k/2) \in O(k^{-1/\g})$ concludes the proof.
\enp

To tackle the general case, we proceed in two steps.
First, we reduce to the case where $m' = m + 1$, and then deal with this case. For the reduction step, the key idea is that the out-degree of a given vertex is Poisson distributed with a finite parameter~\citep{glm2}.
Hence, the number of simplices containing a given point as its youngest vertex has rapidly decaying tail probabilities.
In particular, there are only a few simplices containing a given vertex as its youngest vertex as this number is bounded from above by the outdegree of the vertex at hand.

We want to show that for the higher-order degree of the typical vertex, $o = (0, U)$,
\begin{align} \label{eq:credm}
\limsup_{k \tff}\f1{\log(k)}\log \P (\deg_{ m'}(\De_m) \ge k) \le m - \frac{m + 1}{\g}.
\end{align}
Hence, using Proposition~\ref{pr:palm}, we see that~\eqref{eq:credm} is equivalent to
\begin{align} \label{eq:cred0}
\limsup_{k \tff}\f1{\log(k)}\log\int_0^1 \vp_{k, m, m'}(u)\d u \le m - \frac{m + 1}{\g}.
\end{align}
where
$$ \vp_{k, m, m'}(u) := \int_{\T^m} \P \big( \deg_{m'} (u, \bs o_m) \ge k \big) g_m(u, \bs o_m) \d \bs o_m. $$

%
%PRF RED
%
\bep[Proof of reduction to $m' = m + 1$]
Let $M(\bs o_m) := \bigcap_{j \le m}M(o_j)$ denote the common in-neighbors of $o_1, \dots, o_m$.
Then, the goal of this step is to reduce the problem to deriving the asserted power-law bound for the expression $\P\big(\PP(M(\bs o_m)) \ge k\big)$.
First, Lemma~\ref{lem:gm} gives that $\vp_{k, m, m'}(u) \in O( u^{-\g})$.
Hence, we may assume that $u \ge k^{-2K}$, where $K$ is chosen such that $(1 - \g)K = (m + 1)/\g - m$.

Now, we note that any $m'$-simplex containing the typical $m$-simplex consists of the $m + 1$ vertices of the typical simplex and $m' - m$ additional Poisson points.
In particular, the number of $(m'-m)$-simplices containing the typical vertex $o$ as its youngest vertex is at most $D_{\ms{out}}(o)^{m' - m}$.
Moreover,
\begin{align} \label{eq:out}
\P\big(D_{\ms{out}}(o)^{m' - m}\ge k\big) = \P\big(D_{\ms{out}}(o)\ge k^{\frac{1}{m ' - m}}\big),
\end{align}
which decays stretched exponentially by~\cite[Proposition 4.1]{glm2} and Poisson concentration~\citep[Lemma 1.2]{penrose}.

Hence, it suffices to consider the number $N_{m, m'}$ of $m'$ simplices incident to the typical $m$ simplex with the property that the youngest vertex is one of the $m' - m$ Poisson points $P_i \in \PP$.
Again, the number of $(m' - m)$-simplices having $P_i$ as its youngest vertex is bounded above by $\dego(P_i)^{m' - m}$.
Hence, we have for any $\e > 0$ that
$$\P(N_{m, m'} \ge k) \le \P\Big(\sum_{P_i \in M(\bs o_m)}\hspace{-.3cm}\dego(P_i)^{m' - m} \ge k\Big) \le \P\big(\PP(M(\bs o_m)) \ge k^{1 - \e}\big) + \P\big(\max_{P_i \in M(o)}\hspace{-.1cm} \dego(P_i)^{m' - m} \ge k^{\e}\big).$$
In particular, by the Mecke formula,
\begin{align*}
	\P\big(\max_{P_i \in M(o)} \dego(P_i)^{m' - m} \ge k^{\e}\big) &\le \int_{\R} \int_u^1 \P\big(\dego(x)^{m' - m} \ge k^{\e}\big) \one\{(x, v) \in M(o) \}\d x\d v \\
&= \P\big(\dego(o)\ge k^{\e/(m'-m)} \big)\mu(u).
\end{align*}
Now, the Poisson concentration inequality shows that the probability on the right-hand side decays to $0$ exponentially fast in $k$, whereas the assumption $u \ge k^{-2K}$ gives a polynomial upper bound on $\mu(u)$.
In particular, this step reduces the proof to bounding the expression $\P\big(\PP(M(\bs o_m)) \ge k^{1 - \e}\big)$.
\enp

%
% PRF M' = M +1
%
It remains to consider $m' = m + 1$. During the proof, it is important to control the conditional mean
\begin{align} \label{eq:mpp}
\mu(p, q):= \big| M(p, q) \big| := \big| M(p) \cap M(q) \big|,
\end{align}
i.e., the area of the set of space-time points connecting to both $p = (x, u)$ and $q = (y, v)$ where we henceforth assume that $v \ge u$.

%
% LEM MNU
%
\bel[Bound on $\mu(p, p')$]
\label{lem:mnu}
Put $s(r, u) := (\b u^{-\g}/|r|)^{1/(1 - \g)}$ and $s_\w(r, u) := s(r, u) \w 1 $. Then,
$$\mu(p, q)\le \f{\b}\g v^{-\g} s_\w (x - y, u)^\g \one\{u|x - y|\le \b\}.$$
\enl
\bep
The key observation is that $|x - y| \le |z - x| + |z - y| \le \b u^{-\g}w^{\g - 1}$
for every  $(z, w) \in M(p,q)$.
Hence, $u \le w \le s_\w(x - y, u)$.
In particular,
$$\mu(p, q)\le \one\{u \le s_\w(x - y, u) \} \int_v^{s_\w(x - y, u) }|I_{\b v^{-\g}w^{\g - 1}}| \d w .$$
Hence, a computation of the integral concludes the proof.
\enp

Next, we need to bound suitable integrals on $s(r, u)$.
%
% LEM
%
\bel[Integrals of $s(r, u)$]
\label{lem:sur}
Let $\g < 1$ and $0 < \eta < 1 - \g < \r$.
Then,
\been
\im $\int_{0}^\ff s_\w(r, u)^\r  \d r \in O(u^{-\g}).$
\im $\int_0^{{\b /u}} s_\w(r, u)^{1 - \g - \eta}  \d r \in O(u^{-\g - \eta})$.
\im $\int_0^{\b/r}s_\w(r, u)^\g \d u \in O(r^{-1 -\g})$ if $\g < 1/2$.
\enen
\enl
\bep
{\bf Part (a).}
We distinguish two cases.
First, note that $ |I_{\b u^{-\g}}| \in O(u^{-\g})$
Hence, we may  assume that $y \in \R \sm I_{\b u^{-\g}}$ so that $s_\w(r, u) = s(r, u)$.
Then, as asserted,
$$ \int_{\R \sm I_{\b u^{-\g}}}s(r, u)^\r \d r \le 2\int_{\b u^{\g}}^\ff (u^{\g} r/\b)^{-\r/(1 - \g)} \d r \in O(u^{-\g}). $$

{\bf Part (b).}
We compute that
$$ \int_{0}^{\b/u}s(r, u)^{1 - \g - \eta}\d r \le 2\int_{0}^{\b u^{-1}} (u^{\g} r/\b)^{-(1 - \g - \eta)/(1 - \g)} \d r \le \f{2(1 - \g)}\eta u^{-\g  - \eta}. $$

{\bf Part (c).}
We compute that
$$ \int_0^{\b/r}s(r, u)^{\g}\d u = (r/\b)^{-\g/(1 - \g)}\int_0^{\b/r} u^{-\g^2/(1 - \g)} \d u. $$
The latter integral is of the order $O(r^{\zeta})$ where
$ \zeta = -\g/({1 - \g}) + {\g^2}/({1 - \g}) - 1 = -1 -\g, $
as asserted.
\enp

Finally, we complete the proof of the upper bound in the case $m' = m + 1$.

\bep[Proof of upper bound; $m' = m + 1$]
We need to bound the tail probabilities of the Poisson random variable  $D_m'(u, \bs o_m) := \PP\big(M(o, \bs o_m)\big)$, which has parameter $\mu_m'(u, \bs o_m) := \E[D_m'(u, \bs o_m)]$.
Note that
\[ \P(D_m'(u, \bs o_m) \ge k) \le \P\big(D_m'(u, \bs o_m) \ge k,\, \mu_m'(u, \bs o_m) \le k/2\big) + \one\big\{\mu_m'(u, \bs o_m) \ge k/2\big\}, \]
where by Poisson concentration, the first probability on the right vanishes exponentially.
Moreover, $\mu_m'(u, \bs o_m) \le \min_{n \le m}\mu(o_{n - 1}, o_n)$, where we set $\mu(o_{-1}, o_0):= \mu(u)$.
Hence, it remains to bound
$$\int_0^1 \int_{\T^m}\prod_{n \le m} \one\{\mu(o_{n - 1}, o_n) \ge k\} \d \bs o_m \d u.$$
We start with the innermost integral.
Here, by Lemma~\ref{lem:mnu}, we deduce that if $\mu(o_{m - 1}, o_m)\ge k$, then
$$v_m \le (\b/\g)^{1/\g} k^{-1/\g}s_\w(y_{m - 1} - y_m,v_{m - 1},).$$
Therefore,
\begin{align}
	\label{eq:pot}
	\int_{\T}\one\{\mu(o_{m - 1}, o_m) \ge k\} \d o_m\le (\b/\g)^{1/\g} k^{-1/\g} \int_{-\ff}^\ff s_\w(y_{m - 1} - y_m,v_{m - 1})  \d y_m.
\end{align}
Hence, applying part a of Lemma~\ref{lem:sur} shows that for some $c > 0$,
$$\int_0^1 \int_{\T^m}\prod_{n \le m} \one\{\mu(o_{n - 1}, o_n) \ge k\} \d \bs o_m \d u \le ck^{-1/\g}\int_0^1 \int_{\T^{m - 1}}\prod_{n \le m - 1} \one\{\mu(o_{n - 1}, o_n) \ge k\} v_{m - 1}^{-\g}\d \bs o_{m - 1} \d u.$$

We now continue to compute the integral over $o_{m - 1}$, which is the next innermost integral.
More generally, we claim that for every $n \ge 1$ and  sufficiently small $\eta > 0$, we have that
$$\int_{\T}\one\{\mu(o_{n - 1}, o_n) \ge k\} v_n^{-\g - \eta} \d o_n\in O\big( k^{-(1 - \g - \eta)/\g} v_{n - 1}^{-\g - \eta}\big). $$

First, similarly as in~\eqref{eq:pot}, for some $c > 0$,
$$\int_{\T}\one\{\mu(o_{n - 1}, o_n) \ge k\} v_n^{-\g - \eta} \d o_n\le c k^{-(1 - \g -\eta)/\g}\int_{-I_{\b v_{n - 1}^{-1}}} s_\w(y_{n - 1} - y_n,v_{n - 1})^{1 - \g -\eta}  \d y_n.$$
Therefore, by part (b) of Lemma~\ref{lem:sur},
$$\int_{\T}\one\{\mu(o_{n - 1}, o_n) \ge k\} v_n^{-\g - \eta} \d o_n\in O\big( k^{-(1 - \g - \eta)/\g} v_{n - 1}^{-\g - \eta}\big).$$
as asserted.
\enp

\section{Proof of Theorem~\ref{thm:bet} -- CLT for Betti numbers}
\label{sec:bet}

%
%CLT
%
In this section, we prove Theorem~\ref{thm:bet}.
The idea is to proceed similarly as~\cite[Theorem 5.2]{shirai} and apply the general Poisson CLT from~\citep[Theorem 3.1]{yukCLT}.
While the general strategy is similar to that chosen by~\cite[Theorem 5.2]{shirai}, the long-range dependencies in the ADRCM require more refined argumentation.
Therefore, we provide additional details here.
For a locally finite set $\vp \su \R \ti [0, 1]$ we let $\b(\vp) = \b_q(\vp)$ denote the $q$th Betti number computed of the ADRCM on $\vp$.
To state the conditions of~\citep[Theorem 1]{yukCLT} precisely, we introduce the add-one cost operator
$$ \delta (\vp, u) := \b(\vp \cup \{(0, u)\}) - \b(\vp). $$
Now, to apply~\citep[Theorem 3.1]{yukCLT}, we need to verify the following two conditions:
\been
\im It holds that $\sup_{n \ge 1}\E[ \delta (\PP_n, U)^4] < \ff$ ({\bf moment condition}).
\im It holds that $\delta (\PP \cap W_n, U)$ converges almost surely to a finite limit as $n \tff$ ({\bf weak stabilization}), where $W_n = [-n/2, n/2] \times [0, 1]$.
\enen

We now verify separately the weak stabilization and the moment condition.
In both cases, we follow the general strategy outlined by~\cite[Theorem 5.2]{shirai}.
To make the presentation self-contained, we provide the most important steps of the proof.
First, we consider the moment condition.
%
%PRF MOM
%
\bep[Proof of Theorem~\ref{thm:bet}, moment condition]
As in the proof by~\cite[Theorem 5.2]{shirai}, we note that $\De(\vp, U)$ is bounded above by the number of $q$- and $(q + 1)$-simplices containing $o$.
Thus,
\[ \begin{aligned}
\E[ \delta (\vp, U)^4] &= \int_0^\ff \P[ \delta (\vp, U) \ge s^{1/4}] \d s \le \int_0^\ff \P[ \mathrm{deg}_q(o) + \mathrm{deg}_{q+1}(o) \ge s^{1/4}] \d s \\
&\le \int_0^\ff d_{0, q}(s^{1/4}/2) + d_{0, q + 1}(s^{1/4}/2) \d s,
\end{aligned} \]
where the last inequality holds as if $\mathrm{deg}_q(o) + \mathrm{deg}_{q+1}(o) \ge s^{1/4}$, then at least one of $\mathrm{deg}_q(o)$ or $\mathrm{deg}_{q+1}(o)$ is larger than $s^{1/4} / 2$.
Now, by Theorem~\ref{thm:gen_deg}, both $d_{0, q}$ and $d_{0, q + 1}$ have tail index $1/\g > 4$, thereby showing the finiteness of the above integral.
\enp
Second, we consider the weak stabilization.
%
%PRF STAB
%
\bep[Proof of Theorem~\ref{thm:bet}, weak stabilization]
For $n \ge 1$, we write
$$\b_n := \dim(Z_n) - \dim(B_n) := \dim(Z(\PP_n)) - \dim(B(\PP_n))$$
for the Betti number of the ADRCM constructed on $\PP_n$, noting that this characteristic is the dimension difference of the corresponding cycle space $Z_n$ and boundary space $B_n$, respectively.
Similarly, we set
$$\b_n' := \dim(Z_n') - \dim(B_n') := \dim(Z(\PP_n\cup \{o\})) - \dim(B(\PP_n\cup \{o\})),$$
where we have now added the typical vertex, $o = (0, U)$.
Hence, it suffices to show the weak stabilization with respect to $\dim(Z_n)$ and $\dim(B_n)$ separately.
We now discuss the case of $\dim(Z_n)$, noting that the arguments for $\dim(B_n)$ are very similar.
To check weak stabilization, we show that the sequence $\dim(Z_n') - \dim(Z_n)$ is increasing and bounded.

%
%BOUND
%
First, to show that $\dim(Z_n') - \dim (Z_n)$ is bounded, we note that $\dim(Z_n') - \dim(Z_n) \le \mathrm{deg}_{q,n}(o)$, where, $\mathrm{deg}_{q,n}(o)$ denotes the number of $q$-simplices in $W_n$ containing the typical vertex $o$.
This is because the $q$-simplices constructed from $\PP_n\cup\{o\}$ can be decomposed into the set of $q$-simplices containing the typical vertex $o$ and into the family of all simplices formed in $\PP_n$.
We refer to the arguments by~\cite[Lemma 2.9]{shirai} for the rigorous result.
Now, almost surely, there exists $n_0\ge1$ such that for $n \ge n_0$, the neighbors of $o$ do not change any further.
In particular, 
 $\dim(Z_n') - \dim(Z_n)\le |K_n^0|$.

%
%INC
%
Second, we show that $\dim(Z_n') - \dim (Z_n)$ is nondecreasing.
To that end, we take $n_2 \ge n_1$ and consider the canonical map
$$ Z'_{n_1, q} \to Z'_{n_2, q} / Z_{n_2, q}, $$
where the index $q$ refers to the dimension of the cycle space.
Then, any cycle contained in the kernel of this map consists of simplices formed by vertices in $\PP_n$.
In other words, the kernel equals $Z_{n_1, q}$, which shows that the induced map
$$ Z'_{n_1, q} / Z_{n_1, q} \to Z'_{n_2, q} / Z_{n_2, q} $$
is injective.
In particular,
$\dim(Z_{n_1}') - \dim (Z_{n_1}) \le \dim(Z_{n_2}') - \dim (Z_{n_2}),$
as asserted.
\enp

%
%SEC ALPH
%
\section{Proof of Theorems~\ref{thm:ass} and~\ref{thm:alpha} -- asymptotics of edge counts}
\label{sec:alpha}

In this section, we prove Theorems~\ref{thm:ass} and~\ref{thm:alpha}.
In both results, the idea is to write 
$$ S_n = \sum_{i \le n}T_i :=\sum_{i \le n} \sum_{P_j \in [i - 1, i] \ti [0, 1]} \degi(P_j), $$
i.e., to express the edge count $S_n$ as the sum of the indegrees of all vertices contained in $[i - 1, i] \ti [0, 1]$.
For the proofs of Theorems~\ref{thm:ass} and~\ref{thm:alpha}, it will be important to compute variances of suitable sums of indegrees.
For ease of reference, we therefore state such bounds as a general auxiliary result.
To make this precise, we henceforth let
$$ S(B) = \sum_{P_j \in B} \degi(P_i) $$
denote the indegree sum for all vertices contained in the space-time region $B \su \T$.

%
%LEM VAR
%
\bel[Variance of accumulated indegrees]
\label{lem:var}
Let $\g \ne 1/2$, $A \su \R$ and $u_* > 0$.
Then, there exists a constant $c_{\ms{VI}} > 0$ such that
$$ \Var(S(A \ti [u_*, 1])) \le c_{\ms{VI}}\Big(|A|(1 + u_*^{1 - 2\g}) + \iint_{A\ti A}\int_{u^*}^{1 \wedge (\b/|x - y|)}s_\wedge(u, |x - y|)^\g  \d (x , y) \d u\Big) $$
\enl
\bep
First, we note that
$$ \degi\big((x, u), \PP \cup \{(x, u), (y, v)\}\big) = \degi\big((x, u), \PP \cup \{(x, u)\}\big) +  \one \{(y, v) \in M(x, u)\}\big). $$
Hence, by the Mecke formula~\citep[Theorem 4.4]{poisBook} with  $B:= A \ti [u, 1])$,
\begin{align*}
	\Var(S(B)) =& \int_B \E[\degi(x, u)^2] \d (x, u) + \int_B\int_B\Cov\big(\degi(x, u), \degi(y, v)\big) \d(x, u) \d(y, v)\\
	&+ 2\int_B\E[\degi(x, u)]\big|\{(y, v)\in B\co (x, u) \in M(y, v)\}\big| \d (x, u),
\end{align*}
where $\degi(x, u)$ denotes the in-degree of the vertex $(x, u)$ on the ADRCM constructed on $\PP \cup \{(x, u)\}$.
Now, note that $|\{(y, v)\in B\co (x, u) \in M(y, v)\}| \in O(1)$ and that $\E[\degi(x, u)] \le \E[\degi(x, u)^2]$.
Hence, it suffices to bound the sum
$$ \int_B \E[\degi(x, u)^2] \d (x, u) + \int_B\int_B\Cov\big(\degi(x, u), \degi(y, v)\big) \d(x, u) \d(y, v), $$
and we deal with the two summands separately.

%2 MOM
We start by bounding $\E[\degi(x, u)^2]$.
Since $\degi(x,u)$ is a Poisson random variable with mean $\mu(u) \in O(u^{-\g})$, the Poisson concentration inequality shows that $\E[\degi(x,u)^2] \in O(u^{-2\g})$.
Now, we note that $\int_{u_*}^1 u^{-2\g} = (1 - 2\g)^{-1}(1 - u_*^{1 -2\g})$, which is of order $O(u^{1 - 2\g})$ for $\g > 1/2$ and of order $O(1)$ for $\g < 1/2$.

%COV
To bound the covariance recall from~\eqref{eq:mpp} that a point $(z, w)$ connects to both $(x, u)$ and $(y, v)$ if and only if $(z, w) \in M\big((x, u),  (y, v)\big)$.
Hence, by the independence property of the Poisson process,
$\Cov\big(\degi(x, u), \degi(y, v)\big) = \mu\big((x, u), (y, v)\big)$.
In particular, applying Lemma~\ref{lem:mnu} concludes the proof.
\enp

%
%PRF ASS
%
First, we prove the CLT for the simplex count in the regime $\g < 1/2$, where we will rely on a general CLT for associated random variables~\citep[Theorem 4.4.3]{whitt}.
Here, we recall that the real-valued random variables $T_1, \dots, T_k$ are \emph{associated} if
$$\Cov\big(f_1(T_1, \dots, T_k), f_2(T_1, \dots, T_k)\big) \ge 0.$$
holds for any coordinatewise increasing functions $f_1, f_2\co \R^k \to [0,\ff)$.

\bep[Proof of Theorem~\ref{thm:ass}]
Since the in-degrees are an increasing function in the underlying Poisson point process, we conclude from the Harris-FKG theorem~\citep[Theorem 20.4]{poisBook} that the random variables $\{T_i\}$ are associated.
Hence, to apply~\citep[Theorem 4.4.3]{whitt}, it remains to prove that $\Var(T_1) < \ff$ and $\sum_{k\ge 2}\Cov(T_1, T_{k}) < \ff$.
The finiteness of $\Var(T_1)$ follows from Lemma~\ref{lem:var} so that it remains to consider the covariance sum.

%
%COV SUM
%
 We prove that $\Cov(T_1, T_k) \in O(k^{-1 - \g})$, recalling that $\g < 1/2$.
 Proceeding similarly as in Lemma~\ref{lem:var}, and setting $a := |x - y|$, we need to show that $\int_0^{\b/a} s_\wedge(a,u)^\g   \d u \in O(k^{-\g - 1})$.
 Hence, applying part (c) of Lemma~\ref{lem:sur} concludes the proof.
\enp

%
%RESULT
%
Next, we prove Theorem~\ref{thm:alpha}, i.e., the stable limit theorem for the edge count.
Before proving Theorem~\ref{thm:alpha}, we stress that while there are several general limit results in the literature for deriving the distributional convergence to $\a$-stable limits~\citep{basrak,dst,wolff}, these do not apply in our setting.
More precisely, it is difficult to verify~\citep[Condition 3.3]{basrak} since the ADRCM is mixing but not $\phi$-mixing.
Second,~\cite[Theorem 7.8]{dst} give a general convergence result of Poisson functionals to $\a$-stable random variables with $\a \in (0, 1)$.
However, this corresponds to the case where $\g > 1$, which is not possible due to the model constraints.
While~\cite[Remark 7.9]{dst} state that in principle, the method should generalize to $\a \in (1, 2)$, the ensuing computations may lead to difficulties that are difficult to tackle.
Third,~\cite{wolff} derive a general limit result for U-statistics based on iid input.
However, in our setting, we work in a growing domain, so the distributions change after every step.

%
%IID
%
Before starting the proof of Theorem~\ref{thm:alpha}, it will be convenient to review the classical stable limit theorem for iid sequences from~\citep[Theorem 4.5.2]{whitt}.
To ease presentation, we restrict to the present setting of nonnegative random variables.
More precisely, let $\{X_i\}_i$ be iid nonnegative random variables such that $\P(X_i > x) \sim Ax^{-\a}$ for some $\a \in (1, 2)$ and $A >  0$.
Then, $n^{-1/\a}(\sum_{i \le n}X_i - n\E[X_1])$ converges in distribution to an $\a$-stable random variable $\SS$.

%
%TRUNC IID
%
A key step proving Theorem~\ref{thm:alpha} is a truncation argument, which we first discuss in the iid case.
\bel[Truncation in the iid case]
\label{lem:trunc_iid}
Let  $\{X_i\}_i$ be iid random variables with $\P(X_i > x) \sim Ax^{-\a}$ for some $\a \in (1, 2)$ and $A > 0$.
Then, for every $a< 1/\a$, 
$$n^{-1/\a}\Big(\sum_{i \le n}X_i\one\{X_i \le n^a\} - n\E[X_1\one\{X_1 \le n^a]\Big)\xrightarrow{L^2}0.$$
\enl
\bep
Since the $X_i$ are iid, the claim follows by showing that $\Var(X_i\one\{X_i \le n^a\}) \in o(n^{2/\a - 1})$.
Now,
$$\E[X_i^2\one\{X_i^2 \le n^{2a}\}] = \int_0^{n^{2a}}\P(X_1^2 \in[ r, n^{2a}]) \d r.
$$
Since $\P(X_1^2 \ge r) \asymp A r^{-\a/2}$ we note that $\int_1^{n^{2a}}\P(X_1^2 > r)\d r \in O(n^{2a(1 - \a/2)})$. Hence, observing that $2a(1 - \a/2) < 2/\a - 1$ concludes the proof.
\enp

%
%DEG
%
Now, we return to the case of the edge count in the ADRCM.
The idea of proof is to decompose $S_n$ as $\sna + \snb$, where $\sna$ and $\snb$ contain the contributions of the young and the old vertices, respectively.
More precisely, for $u_n := n^{-0.9}$ put
\[ \sna := \sum_{\substack{P_i \in [0, n]\ti [u_n, 1]}}\degi(P_i),
\quad \text{ and } \quad
\snb := \sum_{\substack{P_i \in [0, n]\ti [0, u_n]}}\degi(P_i). \]

%
%PRF
%
First, we control the deviations of $\sna$ via the Chebyshev inequality.

%
%PR SNA
%
\bepr[$\sna$ is negligible]
\label{pr:sna}
It holds that $n^{-\g}(\sna - \E[\sna])$ converges to 0 in probability.
\enpr
\bep
To prove the claim, we apply Lemma~\ref{lem:var} with $A = [0, n]$ and $u_* = u_n$.
In particular, the first summand in Lemma~\ref{lem:var} is then of order $O(nu_n^{1-2\g})$.
Now, since $-0.9(1 - 2\g) + 1 < 2\g$, we get $nu_n^{1-2\g} \in o(n^{-2\g})$.
Hence, it suffices to bound the second summand in Lemma~\ref{lem:var}.
Here, we can apply part a of Lemma~\ref{lem:sur} which shows that $\int_{\T} s_\w(a, u)^\g \d(a, u) \in O(1)$,
thereby concluding the proof.
\enp

Second, we approximate $\snb$ by a sum of iid Pareto random variables so that we can apply the stable CLT~\citep[Theorem 4.5.2]{whitt}.
%
%PR SNB
%
\bepr[$\snb$ converges to a stable distribution]
\label{pr:snb}
It holds that $n^{-\g}(\snb - \E[\snb])$ converges in distribution to a stable random variable.
\enpr

%
%PRF ALPHA
%
To maintain a clear structure, we conclude the proof of Theorem~\ref{thm:alpha} before establishing Proposition~\ref{pr:snb}.
\bep[Proof of Theorem~\ref{thm:alpha}]
By Proposition~\ref{pr:sna}, $n^{-\g}(\sna - \E[\sna])$ tends to 0 in probability, and $n^{-\g}(\snb - \E[\snb])$ tends in distribution to a stable  random variable.
Hence, also
$$n^{-\g}(S_n - \E[S_n])= n^{-\g}(\sna - \E[\sna])  + n^{-\g}( \snb - \E[\snb])$$
 tends in distribution to a stable random variable.
\enp

It remains to prove Proposition~\ref{pr:snb}.
That is, the renormalized sum of the large in-degrees converges to a stable distribution.
To make this precise, we introduce two further approximations, namely $\snc$ and $\snd$ that we define now.
In these approximations, we replace the in-degree by its expectation, and replace the Poisson number of points in $[0, n]$ by a fixed number, respectively.
More precisely, we set
$$ \snc := \sum_{\substack{X_i \in [0, n]\\ U_i \le u_n}}\mu(U_i)\;\text{ and }\; \snd := \sum_{\substack{i \le n\\ U_i \le u_n}}\mu(U_i). $$
The key step in the proof of Proposition~\ref{pr:snb} is to show that each of these expressions is close in $L^1$-norm.

%
%LEM BCDE
%
\bel[$\snb, \snc$ and $\snd$]
\label{lem:bcde}
It holds that $\E[|\snb - \snc|] + \E[|\snc - \snd|] \in o(n^\g)$.
\enl

Before proving Lemma~\ref{lem:bcde}, we explain how to conclude the proof of Proposition~\ref{pr:snb}.

%
%PRF PR SNB
%
\bep[Proof of Proposition~\ref{pr:snb}]
By Lemma~\ref{lem:bcde}, it suffices to show that $n^{-\g}(\snd - \E[\snd])$ converges in distribution to a stable random variable.
We note that by construction, the summands $\mu(U_i)$, $i\le n$ are iid.
Moreover, by Lemma~\ref{lem:mnu}, $\mu(u) \sim (\b/\g) u^{-\g}$.
Hence, an application of Lemma~\ref{lem:trunc_iid} concludes the proof.
\enp

It remains to prove Lemma~\ref{lem:bcde}.

\bep[Proof of Lemma~\ref{lem:bcde}]
We prove the two parts separately.
\medskip

\ni{$\bs{\E[|\snb - \snc|].}$}
By the Mecke formula, it suffices to show that
$$n\int_{[0, u_n]}\E\big[|\degi(o) - \mu(u)|\big] \d u \in o(n^\g).$$
To achieve this goal, we use that the centered moment of a Poisson random variable with parameter $\la$ is given by $2 \la^{\lfloor \la\rfloor + 1}e^{-\la}/\lfloor \la \rfloor!$.
Specializing to $\la = \mu(u)$ and applying the Stirling formula shows that  $\E[|\degi(o) - \mu(u)|] \in O(u^{-0.6\g})$.
Therefore,
$$n\int_{[0, u_n]}\E[\big|\degi(o) - \mu(u)\big|] \d u \in O(n u_n^{1 - 0.6\g}).$$
Now, since $1 - 0.9 (1 - 0.6)\g < \g$, we deduce that $n u_n^{1 - 0.6\g} \in o(n^\g)$, thereby concluding the proof.
\medskip

\ni{$\bs{\E[|\snc - \snd|].}$}
Let $N$ be a Poisson random variable with parameter $N$.
Then,
$$\E[|\snc - \snd|] \le \E\big[|N - n|\big] \int_0^{u_n}\mu(u)\d u.$$
First, the CLT for iid random variables gives that $\E\big[|N - n|\big] \in O(\sqrt n)$.
Furthermore,
$$\int_0^{u_n}\mu(u)\d u \in O(u_n^{1- \g}).$$
Therefore, $\E\big[\big|\snc - \snd\big|\big] \in O(\sqrt n u_n^{1 - \g})$.
Hence, noting that $1/2 -0.9 + 0.9 \g < \g$ concludes the proof.
\enp

\pagebreak
\section{Proof of Theorem~\ref{thm:tgen_deg}}
\label{sec:tgen_deg}

We deal with parts (a) and (b) of Theorem~\ref{thm:tgen_deg} separately.

%
%PRF A
%
\subsection{Proof of part (a)}
We start with part (a).
In the assertion, we need to establish an upper and a lower bound for the probability that the typical degree in the thinned graph $\gt$ is large.
First, we discuss the lower bound since in the proof, we can ignore the distinction between exposed and protected edges.

%
%PRF A LOW
%
\bep[Proof of Theorem~\ref{thm:tgen_deg}(a), lower bound]
Let $G'$ be obtained by independent edge thinning, where all edges of the ADRCM $G$ are eligible to be removed.
Moreover, the retention probability of an edge $(Y, V) \to o = (0, U)$ is set as $U^\eta$.
Then, $G' \su \gt$ so that
$\P\big(\deg_{\gt, \ms{in}}(o) \ge k \big) \ge \P\big(\deg_{G', \ms{in}}(o) \ge k \big)$.

Now, the thinning theorem for Poisson point processes implies that, conditioned on $U$, the retained in-neighbors form a Poisson point process.
Hence, conditioned on $U$, the in-degree is a Poisson random variable with a mean $U^\eta \mu(U)$.
 Moreover,
 $$\P\big(\deg_{G', \ms{in}}(o) \ge k \big) \ge \P(U^\eta \mu(U) \ge 2k) - \P\big(\deg_{G', \ms{in}}(o) \le k , U^\eta\mu(U) \ge 2k \big). $$
 By Poisson concentration, the second probability on the right decays exponentially in $k$, whereas~\eqref{eq:mu} yields the asserted $\P(U^\eta\mu(U) \ge 2k) \in O\big(k^{-1/(\g - \eta)}\big)$.
\enp

%
%LOW
%
Next, we prove the upper bound for the tail probabilities of the vertex degrees.
That is, an upper bound for the probability that the typical degree is very large.
Since in the model $\gt$ only the exposed edges are thinned out, this is more difficult than the lower bound.
Loosely speaking, we need to ensure that the number of protected edges is negligible so that it does not matter whether or not they are considered in the thinning.
To achieve this goal, in Lemma~\ref{lem:pr}, we bound the power-law exponent of the number of protected edges leading to the typical node $o$.

%
%LEM NONEXP
%
\bel[Power-law for the vertex degree of protected edges]
\label{lem:pr}
It holds that
$$\lim_{k\tff} \log(d^{\ms{pr}, \eta}_{k})/\log(k)=  1 - 2/\g. $$
\enl

%
%PRF A UP
%
Before proving Lemma~\ref{lem:pr}, we conclude the proof of the upper bound in part a of Theorem~\ref{thm:tgen_deg}.
\bep[Proof Theorem~\ref{thm:tgen_deg}(a), upper bound]
As in the proof of the lower bound, we let $G'$ be the graph obtained by independent edge thinning, where we allow all edges to be thinned.
Moreover, we let $I_{\ms{pr}}$ denote the number of protected edges incident to $o$.
Then,
$$\P\big(\deg_{\gt}(o) \ge k \big) \le \P\big(\deg_{G'}(o) \ge k/2 \big) + \P\big(I_{\ms{pr}} \ge k/2 \big).$$
By Lemma~\ref{lem:pr}, the second probability on the right hand side is of order at most $k^{1 - 2/\g + o(1)}$.
Hence, to conclude the proof, we need to show that the first probability is of order at most $k^{-1/(\g - \eta) + o(1)}$.
To that end, we proceed as in the proof of the lower bound.
More precisely,
$$\P\big(\deg_{G'}(o) \ge k \big) \le  \P\big(\deg_{G'}(o) \ge k , U^\eta\mu(U) \le k/2 \big) + \P(U^\eta\mu(U) \ge k/2).$$
Again, by Poisson concentration, the first probability on the right-hand side decays exponentially in $k$, whereas the second one is of order $k^{-1/(\eta - \g) + o(1)}$, as asserted.
\enp

Now, we prove Lemma~\ref{lem:pr}.
The idea is to carefully distinguish between different cases of how an edge can be protected, and then to bound each of the resulting probabilities separately.
%
%PRF
%
\bep[Proof of Lemma~\ref{lem:pr}]
Our goal is bound the probability that the number of protected edges leading to $o$ is at least $k \ge 1$.
By definition,  it suffices to bound $\P(\PP^{(1)} \ge k)$, and $\P(\PP^{(2)} \ge k)$, where
\begin{enumerate}
\im $\PP^{(1)}:= \big\{(Z, W) \in M(o)\co U \le W \le 2U  \big\}.$
\im $\PP^{(2)}:= \big\{(Z, W) \in M(o)\co (Z, W) \to (Y, V)  \text{ for some $(Y, V) \in \PP$ with $V \le 2U \le 4V$}\big\}.$
\end{enumerate}
We now deal with the two cases separately and heavily rely on the result from~\cite[Proposition 4.1]{glm2}, that conditioned on $U = u$, the in-degree of $o$ is Poisson-distributed with mean in $\mu(u)$.
\medskip

%
%C1
%
\ni{$\bs{\P(|\PP^{(1)}| \ge k).}$} We note that conditioned on $U = u$, the quantity $|\PP^{(1)}|$ is a Poisson random variable with mean $\int_u^{2u}|I_{\b u^{-\g}v^{\g - 1}}|\d v = {\b}(2^\g - 1)/\g$.
Hence, $\P\big(|\PP^{(1)}| \ge k\big)$ decays exponentially fast in $k$.
\medskip

%
%C3
%
\ni{$\bs{\P(|\PP^{(2)}| \ge k).}$}
Note that if $(Z, W) \to (0, U)$ and $(Z, W) \to (Y, V)$, then
$$|Y| \le \b|Z|/2 + \b|Z - Y|/2\le \b U^{-\g} W^{\g -1}.$$
For any $\e > 0$, the probabilities
$\P\big(\PP\big([-\b U^{-1},\b U^{-1}] \ti [U, 2U]\big) \ge k^\e\big)$
decay at stretched exponential speed.
Indeed, conditioned on $U = u$, the random variable $\PP\big([-\b u^{-1}, \b u^{-1}] \ti [u, 2u]\big)$ is Poisson distributed with mean $2\b$ so that the asserted decay is a consequence of the Poisson concentration inequality.

Therefore, recalling~\eqref{eq:mpp}, it suffices to bound
$$\int_0^1\int_{u/2}^1\int_{-\ff}^\ff \P\big(\PP\big(M\big((0, u), (y, v)\big)\big) \ge k^{1- \e} \big) \d y\d v\d u.$$
Again, applying the Poisson concentration inequality reduces this task to bounding
\begin{align} \label{eq:p2i}
	\int_0^1\int_{u/2}^1\int_{-\ff}^\ff \one\{\mu\big((0, u), (y, v)\big) \ge k^{1- \e} \big) \d y\d v\d u.
\end{align}
Since $\mu\big((0, u), (y, v)\big) \le \mu( u) $, we conclude that if $\mu\big((0, u), (y, v)\big) \ge k^{1- \e}$, then $u \le c k^{-(1 - \e)/\g}$ for some $c > 0$.
Moreover, we deduce from Lemma~\ref{lem:mnu} that $\mu\big((0, u), (y, v)\big)\le (\b/\g)v^{-\g}s_\w(y, u)^\g$.
Therefore,~\eqref{eq:p2i} is bounded above by
\begin{align*}
(\b/\g)^{1/\g} k^{-(1 - \e)/\g} \int_0^{ck^{-(1 - \e)/\g}}\int_{-\ff}^\ff s_\w(y, u) \d y\d u.
\end{align*}
Hence, an application of part (a) of Lemma~\ref{lem:sur} concludes the proof.
\enp

%
%PRF B
%
\subsection{Proof of part (b)}
Next, we proof part (b) of Theorem~\ref{thm:tgen_deg}.
That is, the thinning operation does not affect the power-law exponent of the edge degree.
Loosely speaking, the idea is that even after removing all exposed edges, the protected edges are sufficient to sustain a positive proportion of all the triangle leading to the high edge degree in the ADRCM.

As in the proof of part (a), we show upper and lower bounds for the tail probabilities separately.
We start with the proof of the upper bound.
Intuitively, it is not surprising that removing edges reduces the edge degrees.
Nevertheless, to make the presentation self-contained, we give a rigorous proof.

%
%PRF UPPER
%
\bep[Proof of part b) of Theorem~\ref{thm:tgen_deg}, upper bound]
The key idea is to use the Palm representation of the typical edge degree.
More precisely,
\begin{align*}
d_{1, k} &= \P\big(\deg_{2}(\De_1) \ge k\big) = \f1{\la_2}\E\Big[\sum_{\substack{(X,U), (Y,V)\in \PP\\ (Y,V) \to (X,U)}} \one\{X \in [0, 1]\}\one\{\deg_2\big((X, U), (Y, V)\big) \ge k\}\Big],
\end{align*}
where $\la_2 > 0$ denotes the edge intensity of $G$.
Similarly, by writing $\to^{\eta}$ to indicate a directed edge in the graph $\gt$, we get that
\begin{align*}
d_{1, k}' &=\f1{\la_2^{(\eta)}}\E\Big[\sum_{\substack{(X,U), (Y,V)\in \PP\\ (Y,V) \to^{\eta} (X,U)}} \one\{X \in [0, 1]\}\one\{\deg_G\big((X, U), (Y, V)\big) \ge k\}\Big],
\end{align*}
where $\la_2^{(\eta)}$ is the edge intensity of the thinned graph.
Now, noting that $\gt$ is a subgraph of $G$ implies that $d_{1, k}' \la_2^{(\eta)} \le d_{1, k}\la_2$.
In particular,
$ \limsup_{k\tff} k^{-1} \log(d_{1, k}') \le 1 - 2/\g, $
as asserted.
\enp

%
%LOWER
%
The lower bound is more delicate since we need to show that triangles formed by the protected edges are sufficient to sustain the original edge degree even after the thinning.
By monotonicity, it suffices to establish the asserted lower bound for the graph $G'' := G^{(\ff)}$, i.e., the graph where only the protected edges are retained.
%
%PRF LOWER
%
\bep[Proof of part b) of Theorem~\ref{thm:tgen_deg}, lower bound]
As a preliminary observation, we note that a directed edge of the form $(Y, V) \to (X, U)$ with $V\le2 U$ is never exposed.
Hence,
as in the non-thinned case in Theorem~\ref{thm:gen_deg}, we need to derive a lower bound for the expression
$$ \int_0^1\int_u^{2u}\int_{-\ff}^\ff T(u, v, y) \d y \d v \d u, $$
where $ T(u, v, y):=  \P\big(\deg_2\big((y, v), (o, u)\big) \ge k\big). $
To achieve this goal, we derive a lower bond for $T(u, v, y)$ when $(u, v, y)$ is in the domain
$$ B_k := [(\b/(64k))^{1/\g}, (\b/(32k))^{1/\g}] \ti [(\b/(32k))^{1/\g}, (\b/(16k))^{1/\g}] \ti [0, 8k]. $$
First, note that $(y, v) \to (0, u)$ for every $(u, v, y) \in B_k$ as $(\b/2) u^{-\g}v^{\g - 1} \ge 16k \ge y.$
Since $|B_k| \in O(k^{1 -2/\g })$, it therefore suffices to show that $T(u, v, y)$ is uniformly bounded away from 0 for $(u, v, y) \in B_k$.

To achieve this goal, we first note that any point $(z, w) \in C_k := [0, 8k] \ti [3/4, 1]$ connects to both $(o, u)$ and $(y, v)$.
Indeed,
$$ |z - 0| \le 8k \le (\b/2)((\b/(32k))^{1/\g})^{-\g}\;\text{ and }\;|z - y| \le 8k \le (\b/2)((\b/(16k))^{1/\g})^{-\g}. $$
Noting that $v \le 2u$ implies that both edges are protected and therefore also exist in $\gth$.
Now, we conclude since the Poisson concentration inequality implies that $\P(\PP(C_k) \ge k)\to1$ as $k \tff$.
\enp

\FloatBarrier
\section{Simulation study}
\label{sec:sim}

\pgfmathdeclarefunction{normalpdf}{3}{\pgfmathparse{1/(#3*sqrt(2*pi))*exp(-((#1-#2)^2)/(2*#3^2))}}

% goal: verify theorems in finite simulated data
This section serves as a bridge between the theory and its applications to real-world data.
Specifically, we study to what extent the methods and limit theorems derived for the ADRCM apply to finite networks.
Our Monte Carlo approach involves simulating multiple networks with identical model parameters.
Subsequently, we calculate various network properties and subject them to statistical analysis, often entailing parameter estimation for theoretical probability distributions.
Relying on Palm calculus, we also explore the simulation of typical simplices in infinite networks to examine fluctuations of different quantities around the limit, devoid of finite size effects.

%\FloatBarrier
%
%SS METH
%
\subsection{Simulation methods}
\label{ss:meth}

%\FloatBarrier
%\subsubsection{Simulating finite networks}

% method to simulate finite sized networks: birth times, rescaling, positions, connections, expansion
To simulate a finite network, we follow a step-by-step process as outlined below.
\begin{enumerate}
\item{
We begin by fixing the network size, setting the volume $V$ of the sampling window equal to the expected vertex number in the network.
The vertex number $N$ is drawn from a Poisson distribution with parameter $V$.
This step determines the actual vertex number in the network.
}
\item{
Next, we generate the birth times of the $N$ vertices.
Conditioned on the vertex number, the birth times are uniformly distributed.
Thus, the birth times are generated by drawing $N$ iid uniformly distributed random variables from the interval $[0, 1]$.
For each vertex, its position is also generated independently and uniformly across the entire sampling window.
This process corresponds to sampling the spatial Poisson point process conditioned on the point count.
}
\item{
Connections between vertices are created based on the following condition.
For every pair of vertices $(x, u)$ and $(y, v)$, where $u \le v$, a connection is formed if the distance between the vertices satisfies $|x - y| \leq \tfrac12 \b u^{-\g} v^{\g - 1}$.
This criterion governs the establishment of connections in the network.
}
\item{
Finally, the generated binary network is expanded to a clique complex.
This simplicial complex allows for topological analysis and examination of higher-order network properties.
}
\end{enumerate}

\begin{figure} [h!] \centering
  \includegraphics[width=0.5\textwidth]{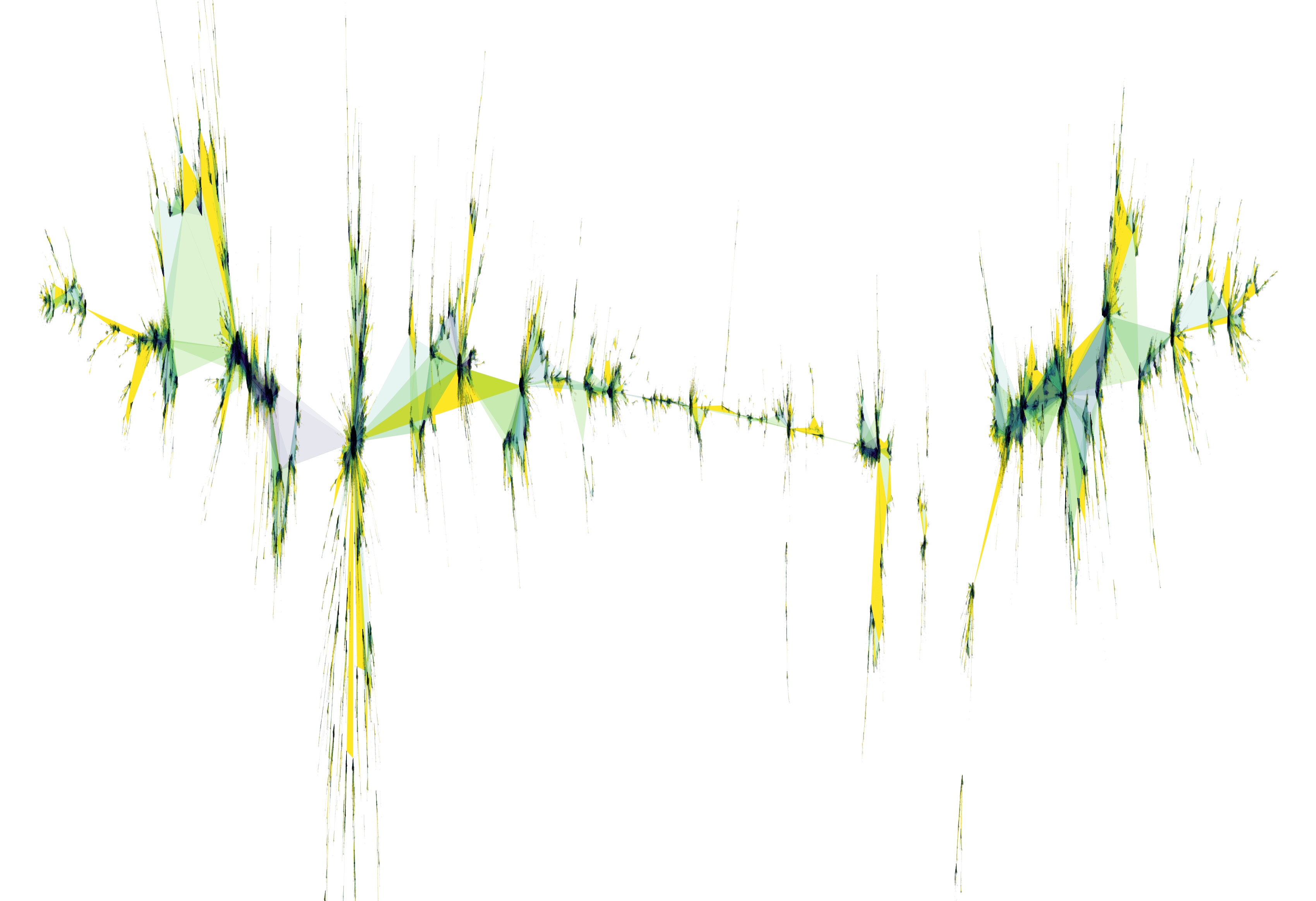}
  \caption{The largest component of a network sample generated by the ADRCM}
  \label{model_sample}
\end{figure}

Figure~\ref{model_sample} shows the largest component of a generated network of size $1\,000\,000$ with $\g = 0.7$.

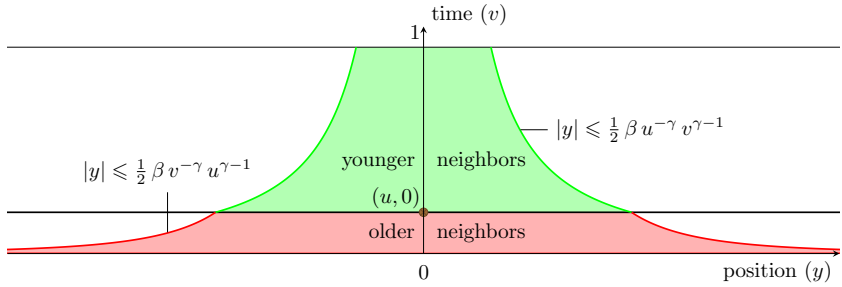
\begin{figure} [h!] \centering \resizebox{0.65\textwidth}{!}{
  \begin{tikzpicture}
    \begin{axis}[
        scale=2,
        axis lines=middle, % axes are at the origin
        axis on top=true, % draw the axis lines on top of the plot
        unit vector ratio*=1 1 1, % unit ratio of the axes are the same
        xmin=-2, xmax=2,
        ymin=0,  ymax=1.1,
        tick style={draw=none},
        xtick={0.001},
        xticklabels={$0$},
        ytick={0.2,1},
        yticklabels={$(u{,} \, 0)$, $1$},
        xlabel=position ($y$),
        x label style={yshift=-0.6cm},
        ylabel=time ($v$),
        y tick label style={
          xshift=+1.5mm,
          yshift=+2.5mm
        },
        y label style={yshift=0.5cm},
      ]
      \addplot [black] {1};
      \addplot [black, thick] {0.2};
      \addplot coordinates {(0, 0.2)};

      \addplot [red,   thick, samples=100, domain=0.0:0.2, name path=olderleft   ] (+0.2 * 1. *   x^(-0.3) * 0.2^(0.3 - 1), x); % 0.5 * beta * v^(-g) * u^(g - 1)
      \addplot [red,   thick, samples=100, domain=0.0:0.2, name path=olderright  ] (-0.2 * 1. *   x^(-0.3) * 0.2^(0.3 - 1), x); % 0.5 * beta * v^(-g) * u^(g - 1)
      \addplot [green, thick, samples=100, domain=0.2:1.0, name path=youngerleft ] (+0.2 * 1. * 0.2^(-0.3) *   x^(0.3 - 1), x); % 0.5 * beta * u^(-g) * v^(g - 1)
      \addplot [green, thick, samples=100, domain=0.2:1.0, name path=youngerright] (-0.2 * 1. * 0.2^(-0.3) *   x^(0.3 - 1), x); % 0.5 * beta * u^(-g) * v^(g - 1)

      \addplot [red!30]   fill between [of=olderleft   and olderright,   reverse=true,];
      \addplot [green!30] fill between [of=youngerleft and youngerright, reverse=true,];

      \draw[-, shorten >=0.1mm, shorten <=0.1mm] (axis cs:-1.2311, 0.30) node[right, align=center, anchor=south] {$|y| \leq \frac12 \, \b \, v^{-\g} \, u^{\g - 1}$} -- (axis cs:-1.2311, 0.1);
      \draw[-, shorten >=0.1mm, shorten <=0.1mm] (axis cs:+0.6   , 0.60) node[right, align=left  , anchor=west ] {$|y| \leq \frac12 \, \b \, u^{-\g} \, v^{\g - 1}$} -- (axis cs:+0.4635, 0.6);
      \draw[-, shorten >=0.1mm, shorten <=0.1mm] (axis cs:-0.3   , 0.10) node[right, align=left  , anchor=west ] {older   \hspace{0.125cm} neighbors};
      \draw[-, shorten >=0.1mm, shorten <=0.1mm] (axis cs:-0.43  , 0.45) node[right, align=left  , anchor=west ] {younger \hspace{0.125cm} neighbors};

    \end{axis}
  \end{tikzpicture}}
  \caption{
    Simulation of the Palm distribution.
    A typical vertex is placed to the origin with fixed birth time $u$.
    The typical vertex connects to older vertices in the red shaded area, whereas younger vertices connect to the typical vertex in the green shaded area of the graph.
  }
  \label{infinite_network_plot}
\end{figure}

%\FloatBarrier
%\subsubsection{Simulating typical simplices in infinite networks}
To avoid the influence of finite size effects and simulate typical simplices in infinite-size networks, we use Palm calculus.
The main idea is to focus only on the immediate neighborhood of a typical vertex placed at the origin, thereby eliminating the presence of finite-size effects for the central vertex.
In this neighborhood, other vertices can form connections with the central vertex and with each other as well.
Any vertex that cannot form a connection with the vertex at the origin is not considered.
The simulation of a single network is visualized in Figure~\ref{infinite_network_plot}.

\begin{enumerate}
\item{
{\bf Typical vertex.}
We begin by randomly placing a vertex $(u,0)$ at the origin of the sampling window with a uniformly distributed birth time $u$.
}
\item{
{\bf Simulaton of older vertices.}
We create older vertices to which the vertex $(u,0)$ connects by simulating a homogeneous Poisson point process in the red shaded area.
The number of older vertices in the red area born up to time $v_0 \leq u$ is Poisson distributed with parameter
\[
  \qquad
		\int \limits_0^{v_0} |I_{ \b v^{-\g} u^{\g - 1} }| \dd v
= \frac{\b}{1-\g} \, u^{\g-1} v_0^{1-\g},
\]
To generate the birth times $\{ v_i \}$ of the points, we simulate a homogeneous Poisson point process $\{ w_i \}$ in the domain $[0, \, u^{1 - \g}]$ with intensity $ {\b u^{\g-1}}/({1-\g})$.
The cardinality of $\{ w_i \}$ will have the same distribution as the point count in the red area.
We then transform $\{ w_i \}$ to receive the set of birth times: $\{ v_i \} = \{ w_i^{1 / (1-\g)} \}$.
The transformation ensures that the birth times $\{ v_i \}$ have the required density.
The positions $y_i$ of the vertices are chosen uniformly in the respective domain $[ - \frac12 \b v_i^{-\g} u^{\g - 1}, \ \frac12 \b v_i^{-\g} u^{\g - 1} ]$.
}
\item{
		{\bf Simulaton of younger vertices.}
Simulation of the younger neighbors of the typical vertex is similar.
The number of younger vertices in the green area born up to time $v_0 \geq u$ is again Poisson distributed with parameter
\[
  \qquad
  \int \limits_u^{v_0} |I_{ \b v^{-\g} u^{\g - 1} }| \dd v
= \frac{\b}{\g} \left( u^{-\g} v_0^{\g} - 1 \right),
\]
To generate the birth times $\{ v_i \}$ using a homogeneous Poisson point process $\{ w_i \}$ in the domain $[u^{\g}, \, 1]$ with intensity ${\b u^{-\g}}/{\g}$, which means that the number of elements in $\{ w_i \}$ will have the same distribution as the number of younger vertices connecting to the typical vertex.
Then, we transform $\{ w_i \}$ as before to get the birth times $\{ v_i \} = \{ w_i^{1 / \g} \}$.
The positions $y_i$ are chosen uniformly in $[ - \frac12 \b u^{-\g} v_i^{\g - 1}, \ \frac12 \b u^{-\g} v_i^{\g - 1} ]$.
}
\item{{\bf Clique complex.}
As before, the generated simple graph is expanded to a clique complex.
Note that those simplices in the clique complex are not subject to finite-size effects, which include the central vertex $(0,u)$ at the origin.
}
\end{enumerate}

\FloatBarrier
\subsection{Higher-order degree distributions of the ADRCM}

% topic of the section
First, we illustrate that the higher-order degree distributions converge to their theoretical limit for increasing network size.

% power-law distributions
To estimate the parameters of power-law distributions, we face two problems.
First, the domain in which the power-law distribution holds is not identical to the entire domain of the data.
As discussed in Section~\ref{sec:mod}, the degree of a typical vertex is a Poisson random variable whose parameter is itself a heavy-tailed random variable.
Thus, the power-law distribution will only be visible for empirical values that are larger than a minimal value $x_{\min}$, from where the influence of the Poisson distributions is negligible.
On the other hand, $x_{\min}$ cannot be too large since in this case the estimation of the power-law exponent becomes too inaccurate due to the low number of values above $x_{\min}$.
Considering these two effects, we carried out a pilot study and found $x_{\min} = 30$ to be a good compromise.
With the domain at hand, we can estimate the exponent $a$ of the power-law distributions via maximum likelihood~\citep{clauset2009power}.
In our setting, this means that
\[ \hat a = 1 + n \, \Big[ \sum_{i \le n} \log \Big( \frac {x_i} {x_{\min} - \frac 12} \Big) \Big]^{-1}, \]
where the index $i$ goes over the data points $x_i \ge x_{\min}$.

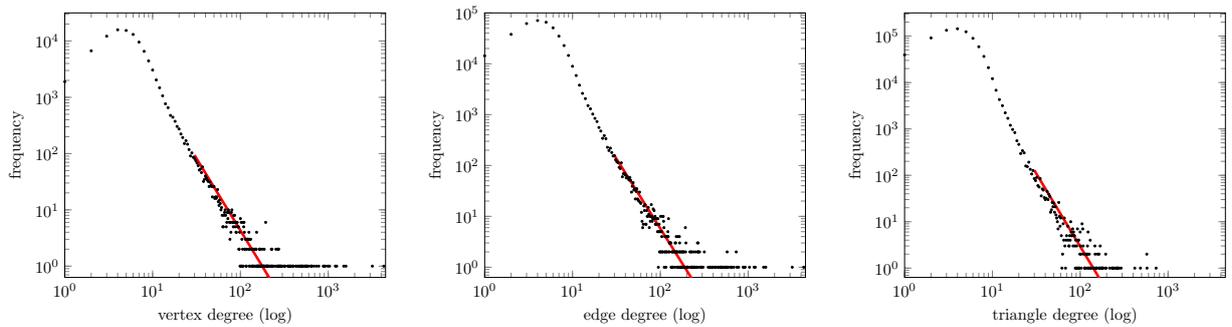
\begin{figure} [h!] \centering
  \hspace*{\fill}
  \begin{subfigure}[b]{0.30\textwidth} \resizebox{\textwidth}{!}{
    \begin{tikzpicture}
      \begin{axis} [
          xmin=1, xmax=4500, ymin=10^(-0.2), ymax=10^(4.5),
          smooth, xmode=log, ymode=log,
          xlabel=vertex degree (log),
          ylabel=frequency,
          ylabel style={
            yshift=-0.2cm
          },
          legend style={at={(0,0)}, anchor=south west, at={(axis description cs:0.025,0.025)}}, legend cell align={left},
        ]
        \addplot[black, only marks, mark options={scale=0.3}] table [x index=0, y index=1, header=false, col sep=comma] {data/model_sample/total_degree_distribution/value_counts.csv};
        \addplot[red, ultra thick, no marks, domain={30:4300}, forget plot] {10^(5.75)*x^(-2.5575)};  % theoretical: -2.4286
      \end{axis}
    \end{tikzpicture}}
  \end{subfigure}
  \hfill
  \begin{subfigure}[b]{0.30\textwidth} \resizebox{\textwidth}{!}{
    \begin{tikzpicture}
      \begin{axis} [
          xmin=1, xmax=4500, ymin=10^(-0.2), ymax=10^(5.0),
          smooth, xmode=log, ymode=log,
          ylabel=frequency,
          xlabel=edge degree (log),
          ylabel style={
            yshift=-0.2cm
          },
          legend style={at={(0,0)}, anchor=south west, at={(axis description cs:0.025,0.025)}}, legend cell align={left},
        ]
        \addplot[black, only marks, mark options={scale=0.3}] table [x index=0, y index=1, header=false, col sep=comma] {data/model_sample/ho_degree_distribution_1/value_counts.csv};
        \addplot[red, ultra thick, no marks, domain={30:4500}, forget plot] {10^(6.25)*x^(-2.7444)};  % theoretical: -2.8571
      \end{axis}
    \end{tikzpicture}}
  \end{subfigure}
  \hfill
  \begin{subfigure}[b]{0.30\textwidth} \resizebox{\textwidth}{!}{
    \begin{tikzpicture}
      \begin{axis} [
          xmin=1, xmax=4500, ymin=10^(-0.2), ymax=10^(5.5),
          smooth, xmode=log, ymode=log,
          ylabel=frequency,
          xlabel=triangle degree (log),
          ylabel style={
            yshift=-0.2cm
          },
          legend style={at={(0,0)}, anchor=south west, at={(axis description cs:0.025,0.025)}}, legend cell align={left},
        ]
        \addplot[black, only marks, mark options={scale=0.3}] table [x index=0, y index=1, header=false, col sep=comma] {data/model_sample/ho_degree_distribution_2/value_counts.csv};
        \addplot[red, ultra thick, no marks, domain={30:1000}, forget plot] {10^(6.75)*x^(-3.1417)};  
      \end{axis}
    \end{tikzpicture}}
  \end{subfigure}
  \hspace*{\fill}
  \caption{Degree distributions of the ADRCM}
  \label{sample_degree_distributions}
\end{figure}

The vertex, edge, and triangle degree distributions of a generated network sample with $100\,000$ vertices and $\g = 0.7$ can be seen in Figure~\ref{sample_degree_distributions} which illustrates the challenges in estimating power-law exponents for degree distributions.
In the small-degree range, the power-law tail of the distribution is hidden due to the  Poisson distribution.
However, as the degrees exceed $\sim 30$, the power-law tail is apparent.

In Theorem~\ref{thm:gen_deg}, we demonstrated that both the ordinary and the higher-order degree distributions follow a power-law tail.
However, this result is rigorously established only for infinitely large networks.
To apply this theorem to real data sets of finite size, it is essential to investigate the extent to which these findings hold for finite networks.

To address this, we conducted Monte Carlo simulations for finite network sizes.
For each network size, we generated $100$ networks with a parameter $\g = 0.7$.
The power-law distribution was then fitted to their degree distributions using the described method.
This process yielded $100$ exponents for vertex, edge, and triangle degree distributions.
Given that the parameters of the underlying ADRCM remained constant, this set of exponents provided a basis for statistical analysis.
By repeating this procedure for networks of varying size, we assessed the convergence of degree distribution exponents to the theoretical limit established in Theorem~\ref{thm:gen_deg}.

Additionally, we examined the simulation of Palm distributions using the same approach.
For this case, $100\,000$ infinite networks were simulated to fit the degree distribution exponents, equivalent to $100\,000$ typical vertices.
The edges and triangles considered in the simulation of the Palm distribution were those involving a special vertex placed at the origin.

The results of the simulations are depicted in Figure~\ref{ho_degree_exponents}, presenting three sets of boxplots summarizing the distribution of the fitted exponents.
The three subfigures visually illustrate the convergence of fitted exponents towards the theoretical limit, indicated by a red horizontal line.
From the observed results, the following conclusions can be drawn.
\begin{figure} [h] \centering
  \begin{subfigure}[b]{0.3\textwidth} \resizebox{\textwidth}{!}{
      \begin{tikzpicture}
        \begin{axis}[
            xmin=0, xmax=6,
            ymin=1.5, ymax=5,
            xlabel=network size,
            axis x line=middle,
            axis y line=middle,
            xtick={1, 2, ..., 5},
            ytick={2, 2.4286, 3, 4, 5},
            xticklabels={$10^2$, $10^3$, $10^4$, $10^5$, $\infty$,},
          ]
          \addplot+[boxplot, boxplot/draw direction=y, black, solid, boxplot prepared={draw position=1, lower whisker=1.8259, lower quartile=2.2788, median=2.9269, upper quartile=4.4992, upper whisker=60.4986, sample size= 7199}] coordinates {};
          \addplot+[boxplot, boxplot/draw direction=y, black, solid, boxplot prepared={draw position=2, lower whisker=1.6487, lower quartile=2.3105, median=2.5512, upper quartile=2.8551, upper whisker=14.7628, sample size=10000}] coordinates {};
          \addplot+[boxplot, boxplot/draw direction=y, black, solid, boxplot prepared={draw position=3, lower whisker=2.2003, lower quartile=2.4364, median=2.5060, upper quartile=2.5881, upper whisker= 3.0121, sample size= 1000}] coordinates {};
          \addplot+[boxplot, boxplot/draw direction=y, black, solid, boxplot prepared={draw position=4, lower whisker=2.4336, lower quartile=2.4818, median=2.4996, upper quartile=2.5356, upper whisker= 2.6149, sample size=  100}] coordinates {};
          \addplot+[boxplot, boxplot/draw direction=y, black, solid, boxplot prepared={draw position=5, lower whisker=2.3624, lower quartile=2.4181, median=2.4470, upper quartile=2.4776, upper whisker= 2.5452, sample size=  100}] coordinates {};
      \addplot[red, thick, domain=0:7] {2.4286};
        \end{axis}
    \end{tikzpicture}} \caption{Vertex degree}
  \end{subfigure}
  \hfill
  \begin{subfigure}[b]{0.3\textwidth} \resizebox{\textwidth}{!}{
      \begin{tikzpicture}
        \begin{axis}[
            xmin=0, xmax=6,
            ymin=1.5, ymax=10,
            xlabel=network size,
            axis x line=middle,
            axis y line=middle,
            xtick={1, 2, ..., 5},
            ytick={2, 2.8571, 4, 5, 10},
            xticklabels={$10^2$, $10^3$, $10^4$, $10^5$, $\infty$,},
          ]
          \addplot+[boxplot, boxplot/draw direction=y, black, solid, boxplot prepared={draw position=1, lower whisker=1.8329, lower quartile=3.0542, median=4.3207, upper quartile=8.0437, upper whisker=60.4986, sample size= 2511}] coordinates {};
          \addplot+[boxplot, boxplot/draw direction=y, black, solid, boxplot prepared={draw position=2, lower whisker=1.5691, lower quartile=2.9418, median=3.5152, upper quartile=4.5820, upper whisker=60.4986, sample size= 9662}] coordinates {};
          \addplot+[boxplot, boxplot/draw direction=y, black, solid, boxplot prepared={draw position=3, lower whisker=2.4888, lower quartile=2.9327, median=3.1218, upper quartile=3.3702, upper whisker= 5.0759, sample size= 1000}] coordinates {};
          \addplot+[boxplot, boxplot/draw direction=y, black, solid, boxplot prepared={draw position=4, lower whisker=2.7523, lower quartile=2.9626, median=3.0462, upper quartile=3.0969, upper whisker= 3.3239, sample size=  100}] coordinates {};
          \addplot+[boxplot, boxplot/draw direction=y, black, solid, boxplot prepared={draw position=5, lower whisker=2.7251, lower quartile=2.8644, median=2.9033, upper quartile=2.9465, upper whisker= 3.0804, sample size=  100}] coordinates {};
          \addplot[red, thick, domain=0:7] {2.8571};
        \end{axis}
    \end{tikzpicture}} \caption{Edge degree}
  \end{subfigure}
  \hfill
  \begin{subfigure}[b]{0.3\textwidth} \resizebox{\textwidth}{!}{
      \begin{tikzpicture}
        \begin{axis}[
            xmin=0, xmax=6,
            ymin=1.5, ymax=10,
            xlabel=network size,
            axis x line=middle,
            axis y line=middle,
            xtick={1, 2, ..., 5},
            ytick={2, 3.2857, 4, 5, 10},
            xticklabels={$10^2$, $10^3$, $10^4$, $10^5$, $\infty$,},
          ]
          \addplot+[boxplot, boxplot/draw direction=y, black, solid, boxplot prepared={draw position=1, lower whisker=1.8401, lower quartile=3.5013, median=5.4145, upper quartile=9.9192, upper whisker=60.4986, sample size=  521}] coordinates {};
          \addplot+[boxplot, boxplot/draw direction=y, black, solid, boxplot prepared={draw position=2, lower whisker=1.5083, lower quartile=3.4999, median=4.6512, upper quartile=7.3912, upper whisker=60.4986, sample size= 5961}] coordinates {};
          \addplot+[boxplot, boxplot/draw direction=y, black, solid, boxplot prepared={draw position=3, lower whisker=2.1959, lower quartile=3.4492, median=3.8596, upper quartile=4.4700, upper whisker=15.8237, sample size=  999}] coordinates {};
          \addplot+[boxplot, boxplot/draw direction=y, black, solid, boxplot prepared={draw position=4, lower whisker=3.1716, lower quartile=3.4737, median=3.6056, upper quartile=3.7754, upper whisker= 4.1080, sample size=  100}] coordinates {};
          \addplot+[boxplot, boxplot/draw direction=y, black, solid, boxplot prepared={draw position=5, lower whisker=3.1072, lower quartile=3.2963, median=3.3612, upper quartile=3.4452, upper whisker= 3.8259, sample size=  100}] coordinates {};
          \addplot[red, thick, domain=0:7] {3.2857};
        \end{axis}
    \end{tikzpicture}} \caption{Triangle degree}
  \end{subfigure}
  \caption{Degree distribution exponents}
  \label{ho_degree_exponents}
\end{figure}
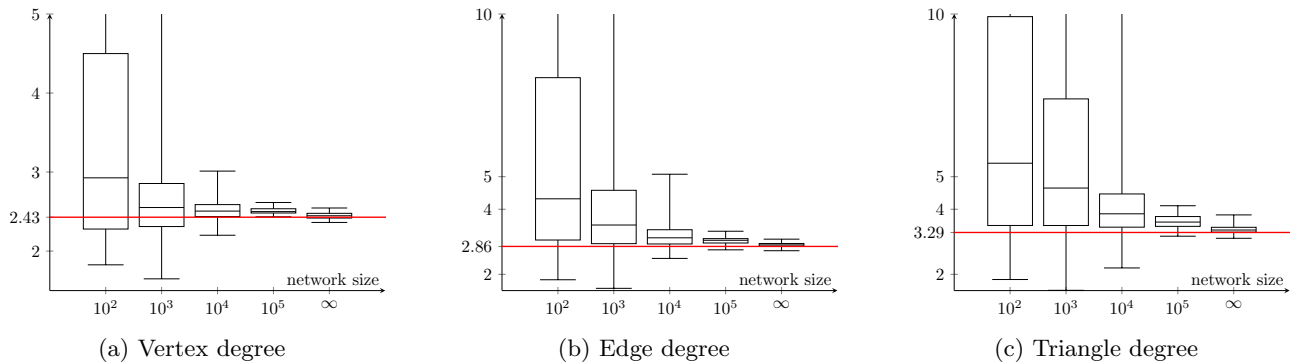
\begin{itemize}
\item{
  As the network size increases, the fluctuation of fitted exponents decreases.
  Smaller networks (with fewer than $1\,000$ vertices) exhibit significant fluctuations, while larger networks (with over $10\,000$ vertices) tend to approach the theoretical limit more closely.
  Infinite networks display the least fluctuations.
}
\item{
  For a given network size, higher dimensions lead to larger fluctuations in the fitted exponents.
  This suggests that considering higher dimensions introduces more variability in estimating the exponents of degree distributions.
}
\item{
  Fitted exponents for finite networks tend to be higher than the theoretical values, indicating a bias in the estimation process.
  This bias is attributed to the constraint on the maximum degree in each dimension due to the finite size, which results in the truncation of degree distribution tails.
  For small degrees, such truncation is absent.
  These effects lead to higher degree distribution exponents.
  The negligible bias observed in the distribution of exponents for ``infinite'' networks supports this explanation.
}
\end{itemize}

\FloatBarrier
\subsection{Edge count of the ADRCM}
\label{ss:edge}

In Theorem~\ref{thm:ass}, we demonstrated that the edge count in large networks follows a normal distribution if $\g < 0.5$.
Conversely, Theorem~\ref{thm:alpha} established that the edge count distribution can be described by a stable distribution if $\g > 0.5$.
To validate these claims in finite networks, we conducted an analysis of the edge count distribution in finite networks containing $100\,000$ vertices.

For each of the three selected values of parameter $\g$ ($0.25$, $0.50$, and $0.60$), we simulated $1000$ networks with $\b = 1$.
Then, we examined the distributions of the edge counts for each of the three cases by fitting both a normal and a stable distribution to the empirical values.

% normal, stable distributions
To fit a normal distribution, we estimated the expectation as the sample mean and the variance as the sample variance.
When fitting the stable distribution, we utilized the insights from Theorem~\ref{thm:alpha} to set the $\alpha$ and $\b$ parameters directly: $\alpha = 1 / \g$ (if $\g < 0.5$, otherwise $\alpha = 2$) and $\b = 1$.
The ``location'' and ``scale'' parameters needed to be estimated from the empirical distribution.
For this purpose, we employed maximum likelihood estimation~\citep{nolan2001maximum}.

\begin{figure} [h] \centering
  \begin{subfigure}[b]{0.3\textwidth} \resizebox{\textwidth}{!}{
    \begin{tikzpicture}
      \begin{axis} [
          grid=none,
          xmin=129600, xmax=137000,
          ymin=0, ymax=0.00045,
          smooth,
          xlabel=edge count,
          ylabel=density,
          ylabel style={
            yshift=0.2cm
          },
          xtick={131000, 133000, 135000},
          scaled y ticks=false,
        ]
        \addplot[ybar interval, mark=no, fill=gray!50!white,] table [x=25_bin_left_limit, y=25_value, col sep=comma] {data/model_test/num_of_edges_normal_mle/linear_histograms.csv};
        \addplot[ultra thick, red, samples=100, domain=129000:137000]{normalpdf(x, 133356.732, 1128.9282)};
      \end{axis}
    \end{tikzpicture}} \caption{Distribution of edge count ($\g = 0.25$, normal distribution)}
  \end{subfigure}
  \hfill
  \begin{subfigure}[b]{0.3\textwidth} \resizebox{\textwidth}{!}{
    \begin{tikzpicture}
      \begin{axis} [
          grid=none,
          xmin=192000, xmax=216000,
          ymin=0., ymax=0.0002,
          smooth,
          xlabel=\phantom a,
          xtick={195000, 205000, 215000},
          scaled y ticks=false,
        ]
        \addplot[ybar interval, mark=no, fill=gray!50!white,] table [x=50_bin_left_limit, y=50_value, col sep=comma] {data/model_test/num_of_edges_normal_mle/linear_histograms.csv};
        \addplot[ultra thick, red, samples=100, domain=190000:292000]{normalpdf(x, 199963.745, 2916.7542)};
      \end{axis}
    \end{tikzpicture}} \caption{Distribution of edge count ($\g = 0.50$, normal distribution)}
  \end{subfigure}
  \hfill
  \begin{subfigure}[b]{0.3\textwidth} \resizebox{\textwidth}{!}{
    \begin{tikzpicture}
      \begin{axis} [
          grid=none,
          xmin=235000, xmax=320000,
          ymin=0, ymax=0.0001,
          smooth,
          xlabel=\phantom a,
          xtick={240000, 270000, 300000},
          scaled y ticks=false,
        ]
        \addplot[ybar interval, mark=no, fill=gray!50!white,] table [x=60_bin_left_limit, y=60_value, col sep=comma] {data/model_test/num_of_edges_normal_mle/linear_histograms.csv};
        \addplot[ultra thick, red, domain=235000:320000] table [x=60_value, y=60_pdf, col sep=comma] {data/model_test/num_of_edges_stable/theoretical_pdf.csv};
      \end{axis}
    \end{tikzpicture}} \caption{Distribution of edge count ($\g = 0.60$, stable distribution)}
  \end{subfigure}
  \\
  \begin{subfigure}[b]{0.33\textwidth} \resizebox{\textwidth}{!}{
    \begin{tikzpicture}
      \begin{axis} [
          grid=none,
          xmin=-3.0, xmax=3.8,
          ymin=-3.0, ymax=3.8,
          clip mode=individual,
          smooth,
          xlabel=normal distribution (standardized),
          ylabel=empirical distribution (standardized),
          xtick={-2, 0, 2},
          scaled y ticks=false,
          y tick label style={/pgf/number format/fixed}
        ]
        \addplot[black, only marks, mark=*, mark size=1.5] table [x=25_theoretical, y=25_empirical, col sep=comma] {data/model_test/num_of_edges_normal_mle/qq_plot.csv};
        \addplot[ultra thick, red, samples=100, domain=-3.0:3.8]{x};
      \end{axis}
    \end{tikzpicture}} \caption{Standardized Q-Q plot with normal distribution ($\g = 0.25$)}
  \end{subfigure}
  \hfill
  \begin{subfigure}[b]{0.3\textwidth} \resizebox{\textwidth}{!}{
    \begin{tikzpicture}
      \begin{axis} [
          grid=none,
          xmin=-3.0, xmax=3.8,
          ymin=-3.0, ymax=5.8,
          clip mode=individual,
          xlabel=\phantom a,
          smooth,
          xtick={-2, 0, 2},
          scaled y ticks=false,
          y tick label style={/pgf/number format/fixed}
        ]
        \addplot[black, only marks, mark=*, mark size=1.5] table [x=50_theoretical, y=50_empirical, col sep=comma] {data/model_test/num_of_edges_normal_mle/qq_plot.csv};
        \addplot[ultra thick, red, samples=100, domain=-3.0:5.8]{x};
      \end{axis}
    \end{tikzpicture}} \caption{Standardized Q-Q plot with normal distribution ($\g = 0.50$)}
  \end{subfigure}
  \hfill
  \begin{subfigure}[b]{0.3\textwidth} \resizebox{\textwidth}{!}{
    \begin{tikzpicture}
      \begin{axis} [
          grid=none,
          xmin=-3.0, xmax=3.5,
          ymin=-3.0, ymax=7.0,
          clip mode=individual,
          xlabel=\phantom a,
          smooth,
          xtick={-2, 0, 2, 4, 6, 8},
          ytick={-2, 0, 2, 4, 6, 8},
          scaled y ticks=false,
          y tick label style={/pgf/number format/fixed}
        ]
        \addplot[black, only marks, mark=*, mark size=1.5] table [x=60_theoretical, y=60_empirical, col sep=comma] {data/model_test/num_of_edges_normal_mle/qq_plot.csv};
        \addplot[ultra thick, red, samples=100, domain=-4.0:7.0]{x};
      \end{axis}
    \end{tikzpicture}} \caption{Standardized Q-Q plot with normal distribution ($\g = 0.60$)}
  \end{subfigure}
  \\
  \begin{subfigure}[b]{0.33\textwidth} \resizebox{\textwidth}{!}{
    \begin{tikzpicture}
      \begin{axis} [
          grid=none,
          xmin=-3.2, xmax=3.8,
          ymin=-3.2, ymax=3.8,
          clip mode=individual,
          xlabel=\phantom a,
          smooth,
          xlabel=stable distribution (standardized),
          ylabel=empirical distribution (standardized),
          xtick={-2, 0, 2, 4},
          scaled y ticks=false,
          y tick label style={/pgf/number format/fixed}
        ]
        \addplot[black, only marks, mark=*, mark size=1.5] table [x=25_theoretical, y=25_empirical, col sep=comma] {data/model_test/num_of_edges_stable/qq_plot.csv};
        \addplot[ultra thick, red, samples=100, domain=-3.2:3.8]{x};
      \end{axis}
    \end{tikzpicture}} \caption{Standardized Q-Q plot with stable distribution ($\g = 0.25$)}
  \end{subfigure}
  \hfill
  \begin{subfigure}[b]{0.3\textwidth} \resizebox{\textwidth}{!}{
    \begin{tikzpicture}
      \begin{axis} [
          grid=none,
          xmin=-3.0, xmax=7.0,
          ymin=-3.0, ymax=7.0,
          clip mode=individual,
          xlabel=\phantom a,
          smooth,
          xtick={-2, 0, 2, 4, 6},
          scaled y ticks=false,
          y tick label style={/pgf/number format/fixed}
        ]
        \addplot[black, only marks, mark=*, mark size=1.5] table [x=50_theoretical, y=50_empirical, col sep=comma] {data/model_test/num_of_edges_stable/qq_plot.csv};
        \addplot[ultra thick, red, samples=100, domain=-3.0:7.0]{x};
      \end{axis}
    \end{tikzpicture}} \caption{Standardized Q-Q plot with stable distribution ($\g = 0.50$)}
  \end{subfigure}
  \hfill
  \begin{subfigure}[b]{0.3\textwidth} \resizebox{\textwidth}{!}{
    \begin{tikzpicture}
      \begin{axis} [
          grid=none,
          xmin=-2.5, xmax=10.5,
          ymin=-2.5, ymax=10.5,
          clip mode=individual,
          xlabel=\phantom a,
          smooth,
          xtick={-2, 0, 2, 4, 6, 8, 10},
          scaled y ticks=false,
          y tick label style={/pgf/number format/fixed}
        ]
        \addplot[black, only marks, mark=*, mark size=1.5] table [x=60_theoretical, y=60_empirical, col sep=comma] {data/model_test/num_of_edges_stable/qq_plot.csv};
        \addplot[ultra thick, red, samples=100, domain=-2.5:10.5]{x};
      \end{axis}
    \end{tikzpicture}} \caption{Standardized Q-Q plot with stable distribution ($\g = 0.60$)}
  \end{subfigure}
  \caption{Distribution of the edge count for different $\g$ parameters}
  \label{simplex_count_distributions}
\end{figure}
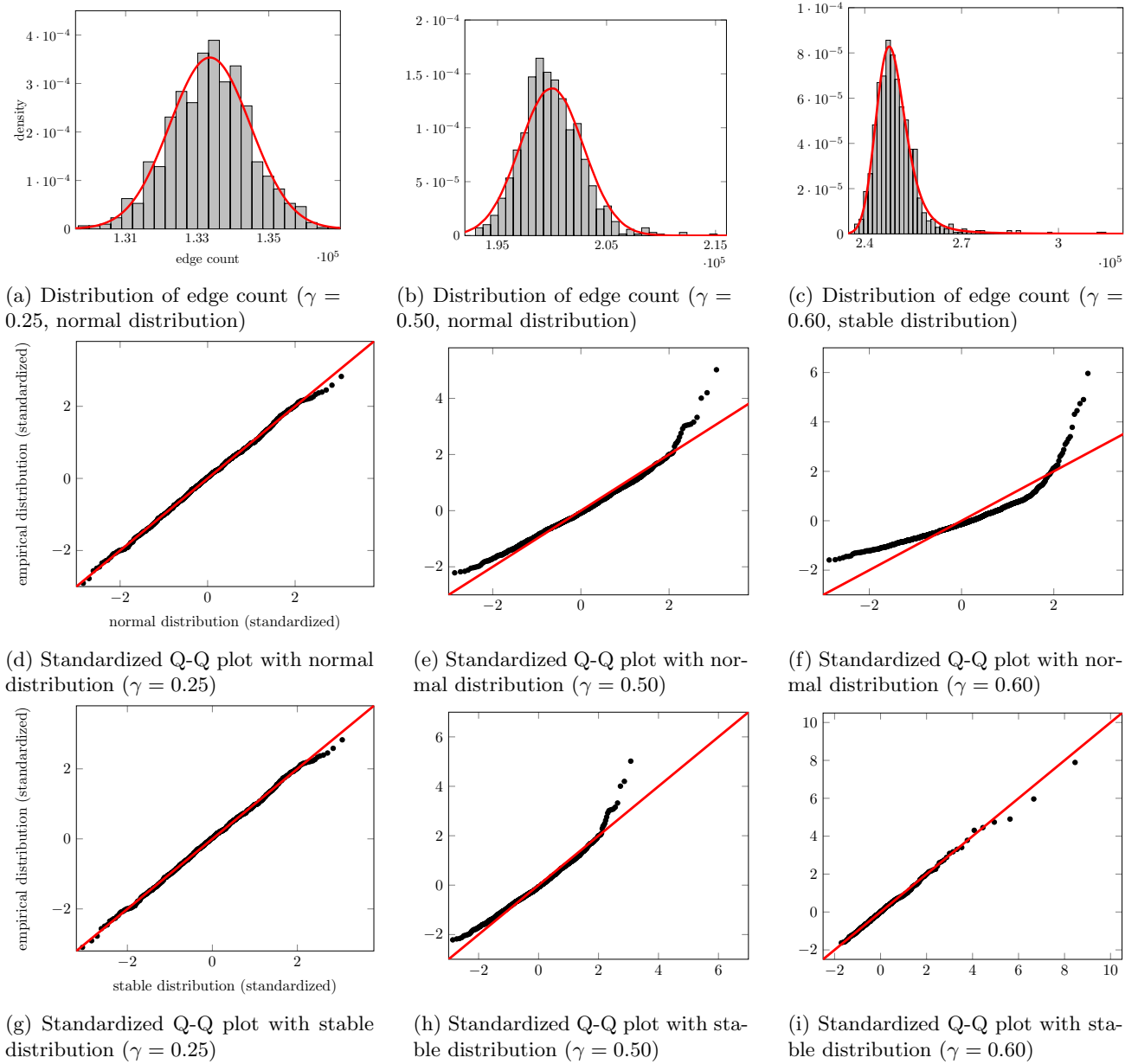

Figure~\ref{simplex_count_distributions} visually represents the results of our analysis, showing the distributions of the edge counts for each of the three cases: $\g = 0.25$, $\g = 0.50$, and $\g = 0.60$.
The subfigures in Figure~\ref{simplex_count_distributions} provide a comprehensive view of the empirical and fitted distributions of the edge counts, along with Q-Q (Quantile-Quantile) plots for comparing the empirical and fitted distributions.
The top row displays the empirical and fitted distributions, while the second and third rows present the Q-Q plots for the fitted normal distributions and the fitted stable distributions.

When $\g = 0.25$, the distribution of edge counts appears symmetric, and the fitted normal distribution closely aligns with the empirical data. However, for $\g = 0.6$, a fat right tail is clearly visible in the empirical distribution.
This heavy-tailed behavior is not adequately captured by the fitted normal distribution, as evidenced by the deviation from the diagonal line in the Q-Q plot for the normal distribution.
In contrast, the stable distribution provides a better fit, aligning well with the data points in the Q-Q plot. Interestingly, for $\g = 0.5$, the normal distribution does not describe the data as effectively as it does for $\g = 0.25$.
We offer two potential reasons for this observation:
\begin{itemize}
\item{
The finite size of the network: In this case, a few high-degree vertices may contribute significantly to the total edge count.
However, for a sufficiently large network, these contributions would be spread among many such vertices, leading to a more normal-like distribution.
}
\item{
The boundary case of $\g = 0.5$: At this value, the degree distributions have an infinite variance, which can affect the distribution characteristics and may not be accurately captured by a normal distribution.
}
\end{itemize}

Supporting the validity of Theorem~\ref{thm:alpha}, the above observations suggest that the normal distribution appears to be a reasonably good fit when $\g < 0.5$, but the stable distribution explains the data more accurately if $\g > 0.5$.

\FloatBarrier
\subsection{Betti numbers of the ADRCM}

To establish the validity of Theorem~\ref{thm:bet} for finite networks, we conducted simulations on finite networks containing $100\,000$ vertices.
In this case, due to computational costs of computing Betti numbers, we performed $100$ simulations for each of the three different values of parameter $\g$: $0.25$, $0.50$, and $0.67$.

For values of $\g = 0.25$, aligning with the findings of Theorem~\ref{thm:bet}, we approximated the empirical values of the first Betti numbers with a normal distribution.
For values of $\g = 0.6$, we observed the distribution of the Betti numbers and conjectured that they follow a stable distribution with stability parameter $\alpha = 1 / \g$.
We posit this based on the expectation that the infinite variance of the simplex count leads to a corresponding infinite variance of the Betti numbers.

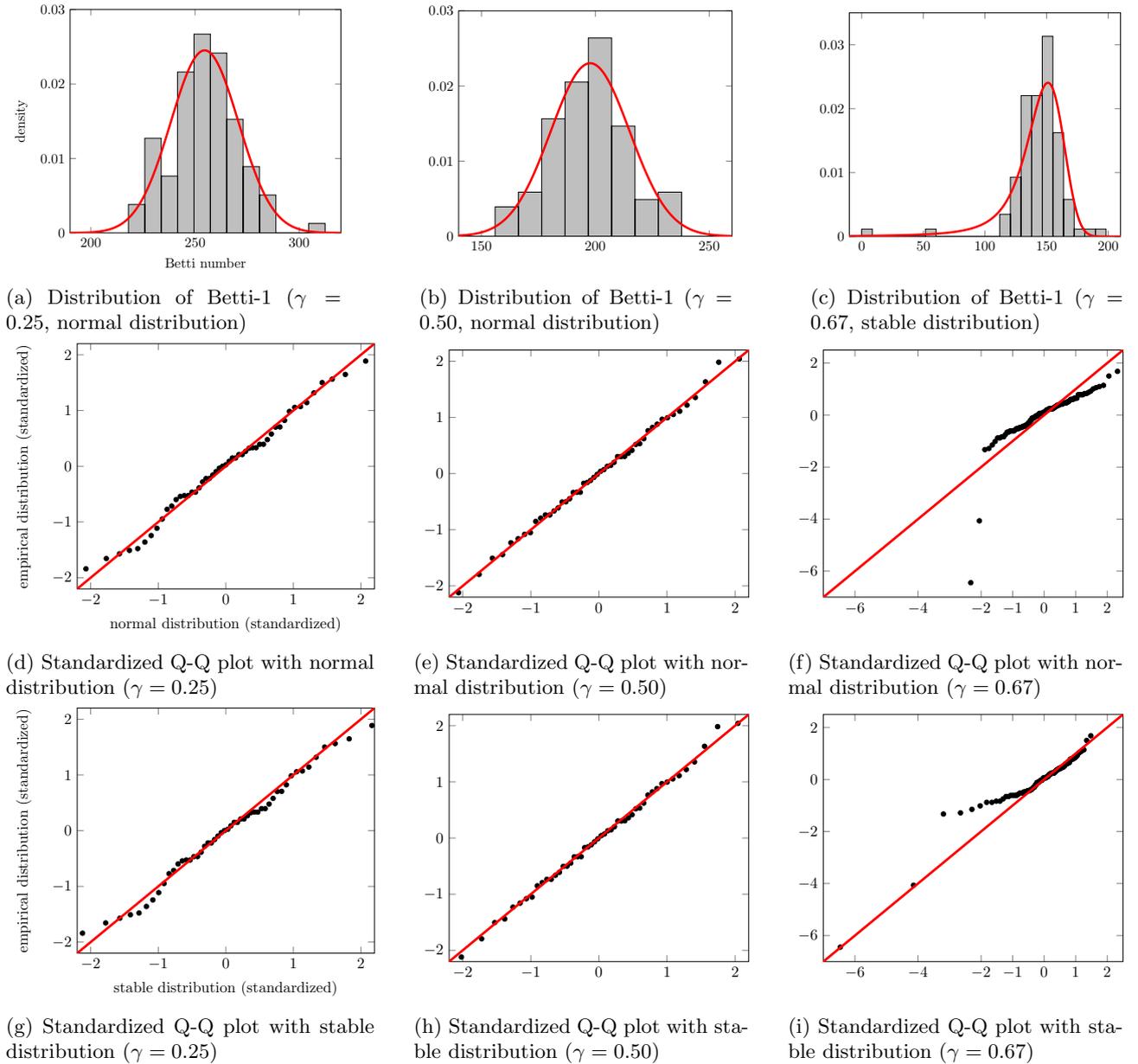
\begin{figure} [h!] \centering
  \begin{subfigure}[b]{0.3\textwidth} \resizebox{\textwidth}{!}{
    \begin{tikzpicture}
      \begin{axis} [
          grid=none,
          xmin=190, xmax=320,
          ymin=0, ymax=0.03,
          smooth,
          xlabel=Betti number,
          ylabel=density,
          xtick={200, 250, 300},
          scaled y ticks=false,
          ylabel style={
            yshift=-0.0cm
          },
          y tick label style={/pgf/number format/fixed}
        ]
        \addplot[ybar interval, mark=no, fill=gray!50!white,] table [x=25_bin_left_limit, y=25_value, col sep=comma] {data/model_test/betti_number_1_normal/linear_histograms.csv};
        \addplot[ultra thick, red, samples=100, domain=190:320]{normalpdf(x, 254.57, 16.2747)};
      \end{axis}
    \end{tikzpicture}} \caption{Distribution of Betti-1 ($\g = 0.25$, normal distribution)}
  \end{subfigure}
  \hfill
  \begin{subfigure}[b]{0.28\textwidth} \resizebox{\textwidth}{!}{
    \begin{tikzpicture}
      \begin{axis} [
          grid=none,
          xmin=140, xmax=260,
          ymin=0, ymax=0.03,
          smooth,
          xlabel=\phantom a,
          xtick={150, 200, 250},
          scaled y ticks=false,
          y tick label style={/pgf/number format/fixed}
        ]
        \addplot[ybar interval, mark=no, fill=gray!50!white,] table [x=50_bin_left_limit, y=50_value, col sep=comma] {data/model_test/betti_number_1_normal/linear_histograms.csv};
        \addplot[ultra thick, red, samples=100, domain=140:260]{normalpdf(x, 197.75, 17.3282)};
      \end{axis}
    \end{tikzpicture}} \caption{Distribution of Betti-1 ($\g = 0.50$, normal distribution)}
  \end{subfigure}
  \hfill
  \begin{subfigure}[b]{0.28\textwidth} \resizebox{\textwidth}{!}{
    \begin{tikzpicture}
      \begin{axis} [
          grid=none,
          xmin=-10, xmax=210,
          ymin=0, ymax=0.035,
          smooth,
          xlabel=\phantom a,
          xtick={0, 50, 100, 150, 200},
          scaled y ticks=false,
          y tick label style={/pgf/number format/fixed}
        ]
        \addplot[ybar interval, mark=no, fill=gray!50!white,] table [x=67_bin_left_limit, y=67_value, col sep=comma] {data/model_test/betti_number_1_normal/linear_histograms.csv};
        \addplot[ultra thick, red, domain=-10:250] table [x=67_value, y=67_pdf, col sep=comma] {data/model_test/betti_number_1_stable/theoretical_pdf.csv};
      \end{axis}
    \end{tikzpicture}} \caption{Distribution of Betti-1 ($\g = 0.67$, stable distribution)}
  \end{subfigure}
  \\
  \begin{subfigure}[b]{0.33\textwidth} \resizebox{\textwidth}{!}{
    \begin{tikzpicture}
      \begin{axis} [
          grid=none,
          xmin=-2.2, xmax=2.2,
          ymin=-2.2, ymax=2.2,
          clip mode=individual,
          smooth,
          xlabel=normal distribution (standardized),
          ylabel=empirical distribution (standardized),
          xtick={-2, -1, 0, 1, 2},
          scaled y ticks=false,
          y tick label style={/pgf/number format/fixed}
        ]
        \addplot[black, only marks, mark=*, mark size=1.5] table [x=25_theoretical, y=25_empirical, col sep=comma] {data/model_test/betti_number_1_normal/qq_plot.csv};
        \addplot[ultra thick, red, samples=100, domain=-2.2:2.2]{x};
      \end{axis}
    \end{tikzpicture}} \caption{Standardized Q-Q plot with normal distribution ($\g = 0.25$)}
  \end{subfigure}
  \hfill
  \begin{subfigure}[b]{0.3\textwidth} \resizebox{\textwidth}{!}{
    \begin{tikzpicture}
      \begin{axis} [
          grid=none,
          xmin=-2.2, xmax=2.2,
          ymin=-2.2, ymax=2.2,
          clip mode=individual,
          smooth,
          xlabel=\phantom a,
          xtick={-2, -1, 0, 1, 2},
          scaled y ticks=false,
          y tick label style={/pgf/number format/fixed}
        ]
        \addplot[black, only marks, mark=*, mark size=1.5] table [x=50_theoretical, y=50_empirical, col sep=comma] {data/model_test/betti_number_1_normal/qq_plot.csv};
        \addplot[ultra thick, red, samples=100, domain=-2.2:2.2]{x};
      \end{axis}
    \end{tikzpicture}} \caption{Standardized Q-Q plot with normal distribution ($\g = 0.50$)}
  \end{subfigure}
  \hfill
  \begin{subfigure}[b]{0.3\textwidth} \resizebox{\textwidth}{!}{
    \begin{tikzpicture}
      \begin{axis} [
          grid=none,
          xmin=-7.0, xmax=2.5,
          ymin=-7.0, ymax=2.5,
          clip mode=individual,
          smooth,
          xlabel=\phantom a,
          xtick={-8, -6, -4, -2, -1, 0, 1, 2},
          scaled y ticks=false,
          y tick label style={/pgf/number format/fixed}
        ]
        \addplot[black, only marks, mark=*, mark size=1.5] table [x=67_theoretical, y=67_empirical, col sep=comma] {data/model_test/betti_number_1_normal/qq_plot.csv};
        \addplot[ultra thick, red, samples=100, domain=-10.0:2.5]{x};
      \end{axis}
    \end{tikzpicture}} \caption{Standardized Q-Q plot with normal distribution ($\g = 0.67$)}
  \end{subfigure}
  \\
  \begin{subfigure}[b]{0.33\textwidth} \resizebox{\textwidth}{!}{
    \begin{tikzpicture}
      \begin{axis} [
          grid=none,
          xmin=-2.2, xmax=2.2,
          ymin=-2.2, ymax=2.2,
          clip mode=individual,
          smooth,
          xlabel=stable distribution (standardized),
          ylabel=empirical distribution (standardized),
          xtick={-2, -1, 0, 1, 2},
          scaled y ticks=false,
          y tick label style={/pgf/number format/fixed}
        ]
        \addplot[black, only marks, mark=*, mark size=1.5] table [x=25_theoretical, y=25_empirical, col sep=comma] {data/model_test/betti_number_1_stable/qq_plot.csv};
        \addplot[ultra thick, red, samples=100, domain=-2.2:2.2]{x};
      \end{axis}
    \end{tikzpicture}} \caption{Standardized Q-Q plot with stable distribution ($\g = 0.25$)}
  \end{subfigure}
  \hfill
  \begin{subfigure}[b]{0.3\textwidth} \resizebox{\textwidth}{!}{
    \begin{tikzpicture}
      \begin{axis} [
          grid=none,
          xmin=-2.2, xmax=2.2,
          ymin=-2.2, ymax=2.2,
          clip mode=individual,
          smooth,
          xlabel=\phantom a,
          xtick={-2, -1, 0, 1, 2},
          scaled y ticks=false,
          y tick label style={/pgf/number format/fixed}
        ]
        \addplot[black, only marks, mark=*, mark size=1.5] table [x=50_theoretical, y=50_empirical, col sep=comma] {data/model_test/betti_number_1_stable/qq_plot.csv};
        \addplot[ultra thick, red, samples=100, domain=-2.2:2.2]{x};
      \end{axis}
    \end{tikzpicture}} \caption{Standardized Q-Q plot with stable distribution ($\g = 0.50$)}
  \end{subfigure}
  \hfill
  \begin{subfigure}[b]{0.3\textwidth} \resizebox{\textwidth}{!}{
    \begin{tikzpicture}
      \begin{axis} [
          grid=none,
          xmin=-7.0, xmax=2.5,
          ymin=-7.0, ymax=2.5,
          clip mode=individual,
          smooth,
          xlabel=\phantom a,
          xtick={-8, -6, -4, -2, -1, 0, 1, 2},
          scaled y ticks=false,
          y tick label style={/pgf/number format/fixed}
        ]
        \addplot[black, only marks, mark=*, mark size=1.5] table [x=67_theoretical, y=67_empirical, col sep=comma] {data/model_test/betti_number_1_stable/qq_plot.csv};
        \addplot[ultra thick, red, samples=100, domain=-10.0:2.5]{x};
      \end{axis}
    \end{tikzpicture}} \caption{Standardized Q-Q plot with stable distribution ($\g = 0.67$)} \label{betti_qq_plot_067}
  \end{subfigure}
  \caption{Distribution of the first Betti numbers with different $\g$ parameters}
  \label{betti_number_distributions}
\end{figure}

In all cases, the parameter $\b$ of the stable distribution remained constant at $-1$. As in Section \ref{ss:edge}, we estimated the remaining parameters of the fitted distributions via maximum likelihood.
The results are visualized in Figure~\ref{betti_number_distributions}.

From the Q-Q plots it is evident that the fitted normal distribution provides a satisfactory approximation to the distribution of the Betti numbers for the simulations with $\g = 0.25$ and $\g = 0.50$.
The points on the Q-Q plots are closely aligned with the diagonal line, indicating a good fit.

However, for $\g = 0.67$, the distribution displays a heavy left tail, which is clearly visible both from the histogram and the Q-Q plot against the normal distribution: the points on the Q-Q plot significantly deviate from the diagonal line in the lower quantiles.
The shallow slope of the points in the central section suggests that the standard deviation is not accurately captured by the normal distribution, which is also an artifact of the heavy left tail.
In contrast, the stable distribution fits the histogram more accurately, as visualized in the Q-Q plot shown in Plot~\ref{betti_qq_plot_067}.
We can also see that the left tail is not entirely accurate in the stable distribution case.
This is explained by two effects:
\begin{itemize}
\item{the simulation number is low, thus there are not enough values in the left tail to precisely estimate the distribution;}
\item{the minimum value of in the distribution is $0$, suggesting the presence of finite-size effects.}
\end{itemize}
All in all, the points on the Q-Q plot against the stable distribution follow more closely the diagonal line both in the central region and in the left tail of the distribution, reinforcing our earlier conjecture about the stable distribution of Betti numbers for $\g > 0.5$.

\FloatBarrier
\section{Analysis of collaboration networks}
\label{sec:dat}

\FloatBarrier

% introduction
In this section, we analyze four datasets collected from arXiv to showcase the applications of our results and to further motivate our model extensions.
As higher-order relationships appear naturally in the case of scientific collaborations, we chose to analyze a publicly available dataset of scientific papers.
The authors of the papers are represented as vertices in a simplicial complex, whereas each paper represents a higher-order interaction of the authors.

% Patania: The shape of collaborations
\cite{patania} also investigates higher-order collaboration networks on the arxiv data and extend the concept of triadic closure to higher dimensions.
However, their analysis was purely of empirical nature and did not consider the question of using a stochastic higher-order network model.
In contrast, we are compare the arxiv dataset with the ADRCM, and also perform hypothesis tests.
We also note that although we consider a different time frame, the Betti numbers we found are largely consistent with the results published by~\cite{patania}.

\FloatBarrier
\subsection{Datasets}
\label{ss:dat}

% Data collection
We analyze all available documents uploaded to arXiv from various scientific fields.
For each document, we extracted the author names, the publication time and its primary category.
The datasets were built using the primary categories of the documents the authors specified.
\begin{itemize}
\item{\textbf{Computer Science (cs)}:
The {computer science} dataset is the largest we analyze with more than $400\,000$ authors.
}
\item{\textbf{Engineering (eess)}:
The second dataset we analyze consists of documents from the scientific field of electrical engineering, which is built from around $80\,000$ authors.
}
\item{\textbf{Mathematics (math)}:
The {mathematics} dataset encompasses around $200\,000$ authors.
}
\item{\textbf{Statistics (stat)}:
The smallest dataset we analyze contains documents from the field of {statistics}, including around $45\,000$ authors.
}
\end{itemize}

\begin{figure} [h] \centering
  \begin{subfigure}[b]{0.24\textwidth}
    \includegraphics[width=\textwidth]{data/data/computer_science_network.png} \caption{cs}
  \end{subfigure}
  \hfill
  \begin{subfigure}[b]{0.24\textwidth}
    \includegraphics[width=\textwidth]{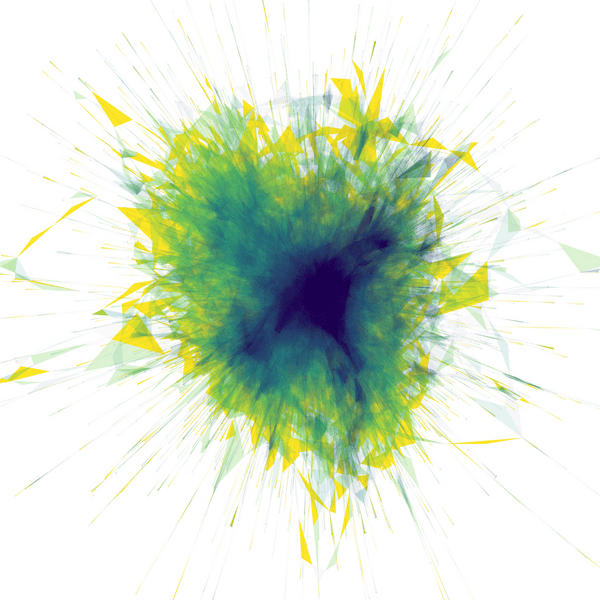}      \caption{eess}
  \end{subfigure}
  \hfill
  \begin{subfigure}[b]{0.24\textwidth}
    \includegraphics[width=\textwidth]{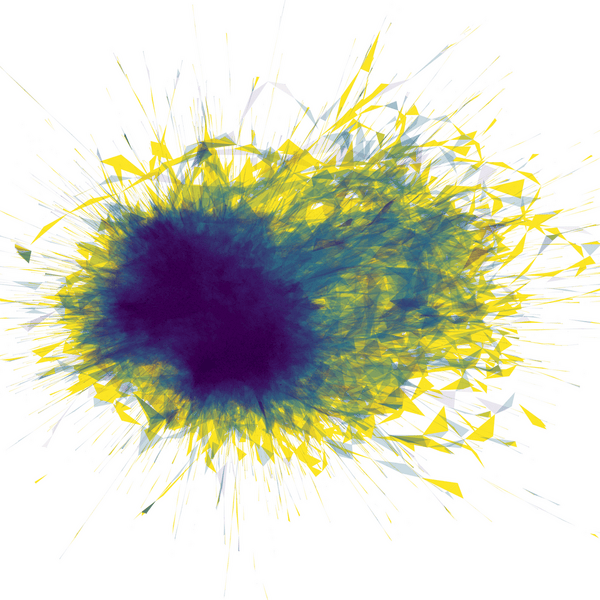}      \caption{math}
  \end{subfigure}
  \hfill
  \begin{subfigure}[b]{0.24\textwidth}
    \includegraphics[width=\textwidth]{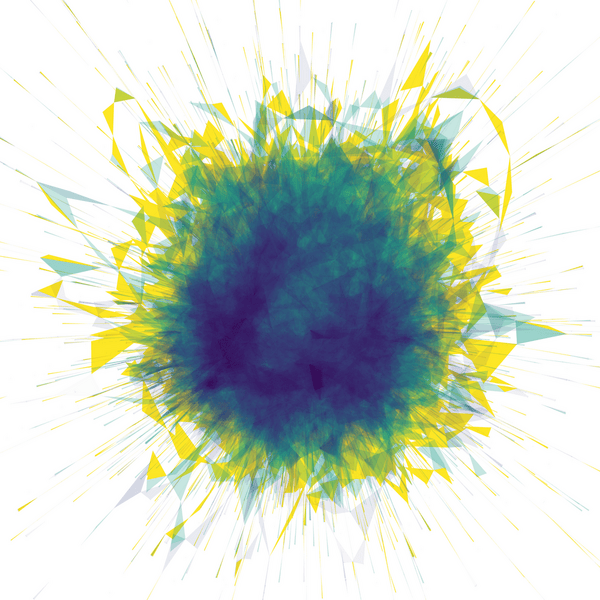}       \caption{stat}
  \end{subfigure}
  \caption{The largest component of the simplicial complexes built from the data sets}
  \label{network_plots}
\end{figure}

The largest components of the datasets are visualized in Figure~\ref{network_plots}, and their most important characteristics are summarized in Table~\ref{dataset_properties}.

% Merge authors with the same names
As arXiv does not uniquely identify authors, we chose to use their full names as identifiers.
Although in the case of a common full name, this method will result in treating distinct authors as if they were the same, the effect of these identifier collisions is greatly reduced as we consider scientific fields separately.

% Build data set
After identifying the authors as vertices, each document is considered as a higher-order interaction of the authors.
This means that every document with $n+1$ authors is represented by an $n$-simplex.
Furthermore, as our goal is to build a simplicial complex, every lower-dimensional face of this $n$-simplex is also added to the simplicial complex to ensure that it is closed under taking subsets.

% combinatorial explosion
A document with $n + 1$ authors has ${n+1 \choose m+1}$ $m$-faces, so in total $\sum_{m=0}^n {n+1 \choose m+1} = 2^{n+1} - 1$ number of simplices needs to be considered.
This poses a twofold practical implementation challenge.
\begin{enumerate}
\item{Due to computational reasons, we must limit the maximum dimension of the simplicial complex.}
\item{
The more authors a document has, the higher its influence is on the simplicial complex, as the number of simplices grows exponentially with the number of authors.
\cite{carstens} have also found the same problem when analyzing collaboration networks.
They tackled this problem by weighting the simplices: they assigned greater weights to smaller simplices, and to those, in which the represented collaboration was frequent.
Although introducing weighted simplices is possible, it is beyond the scope of our present work.
}
\end{enumerate}
Taking into account the above aspects, we consider interactions with a dimension of at most 20 (which, including all the faces, means more than 2 million simplices in total for a document with 21 authors).
To further reduce computational complexity, we analyzed the 2-skeleton of the collaboration network, with the triangles being the highest dimensional simplices.

% network plot images, main characteristics
Using this procedure, we built four separate datasets, each representing publications of a specific scientific field published up to August 4, 2023.
As we will see, the nature of collaborations significantly differ in the four cases, thus, by considering the four scientific fields separately, we can examine how the ADRCM model behaves for four distinct scientific communities.
\begin{table} [h] \centering \caption{Main properties of the datasets} \label{dataset_properties}
  \begin{tabular}{|l|r|r|r|r|} \hline
    dataset &  authors & documents & components & size of largest component \\ \hline
    cs      & 433\,244 &  452\,881 &    22\,576 &                  370\,494 \\
    eess    &  77\,686 &   69\,594 &     5\,533 &                   54\,147 \\
    math    & 198\,601 &  466\,428 &    26\,197 &                  152\,441 \\
    stat    &  44\,380 &   36\,689 &     4\,049 &                   32\,373 \\ \hline
  \end{tabular}
\end{table}

\begin{figure} [h] \centering \resizebox{0.5\textwidth}{!}{
  \begin{tikzpicture}
    \begin{axis} [
        xmin=0, xmax=21,
        smooth,
        ymode=log,
        xlabel=authors per document,
        ylabel=frequency,
        legend style={
          at={(0,0)},
          anchor=south west,
          at={(axis description cs:0.025,0.025)},
          nodes={scale=0.75, transform shape}
        },
        legend cell align={left},
      ]
      \addplot [orange,         only marks, mark=triangle, mark options={scale=1.2}] table [x=value, y=computer_science, col sep=comma] {data/data/interaction_dimension_distribution/value_counts.csv};
      \addplot [red,            only marks, mark=diamond,  mark options={scale=1.2}] table [x=value, y=engineering     , col sep=comma] {data/data/interaction_dimension_distribution/value_counts.csv};
      \addplot [green!70!black, only marks, mark=asterisk, mark options={scale=1.0}] table [x=value, y=mathematics     , col sep=comma] {data/data/interaction_dimension_distribution/value_counts.csv};
      \addplot [blue,           only marks, mark=oplus,    mark options={scale=1.0}] table [x=value, y=statistics      , col sep=comma] {data/data/interaction_dimension_distribution/value_counts.csv};

      \legend{cs,eess,math,stat}
    \end{axis}
  \end{tikzpicture}} \caption{Distribution of authors per documents} \label{interaction_dimension_distribution}
\end{figure}
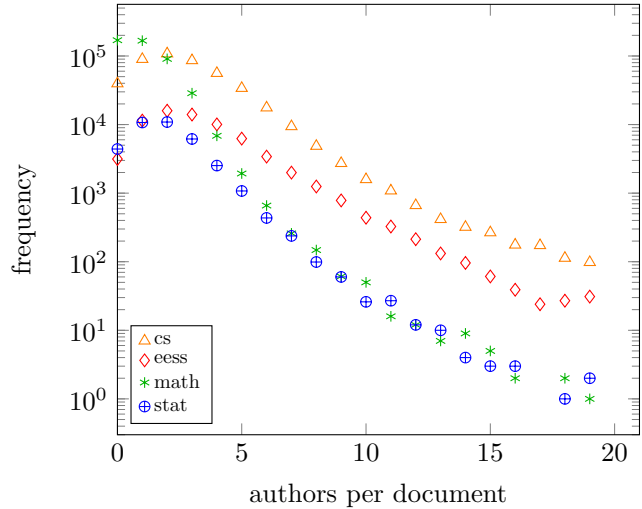

It is also interesting to analyze the distribution of the dimension of the higher-order interactions, or, equivalently, of the per-document author count.
Figure~\ref{interaction_dimension_distribution} visualizes the distribution of the per-document author count for each datasets revealing the typical size of the collaborations scientists participate in within each of the examined scientific fields.
% diverse datasets -> comprehensive examination
The distribution related to the cs and eess datasets have the fattest tails, i.e., relatively higher number of documents have more authors.
On the other hand, the opposite is true for the math and stat datasets, where most papers tend to have a lower number of authors.
The dataset diversity  of the different fields opens the opportunity to comprehensively examine the application of the theorems stated in Section~\ref{sec:mod}.

\begin{figure} [h!] \centering
  \hspace*{\fill}
  \begin{subfigure}[b]{0.24\textwidth} \resizebox{\textwidth}{!}{
    \begin{tikzpicture}
      \begin{axis} [
          xmin=1, xmax=3000, ymin=10^(-0.2), ymax=10^(5.5),
          clip mode=individual,
          smooth, xmode=log, ymode=log,
          xlabel=vertex degree, ylabel=frequency,
          ylabel style={
            yshift=-0.4cm
          },
          legend style={at={(0,0)}, anchor=south west, at={(axis description cs:0.025,0.025)}}, legend cell align={left},
        ]
        \addplot[only marks, mark options={scale=0.4}] table [x=value, y=computer_science, col sep=comma] {data/data/total_degree_distribution/value_counts.csv};
        \addplot[red, ultra thick, no marks, domain={10:3000}] {10^(6.5)*x^(-2.3908)};
      \end{axis}
    \end{tikzpicture}}
  \end{subfigure}
  \hfill
  \begin{subfigure}[b]{0.24\textwidth} \resizebox{\textwidth}{!}{
    \begin{tikzpicture}
      \begin{axis} [
          xmin=1, xmax=600, ymin=10^(-0.2), ymax=10^(4.5),
          clip mode=individual,
          smooth, xmode=log, ymode=log,
          xlabel=vertex degree,
          ylabel style={
            yshift=-0.6cm
          },
          legend style={at={(0,0)}, anchor=south west, at={(axis description cs:0.025,0.025)}}, legend cell align={left},
        ]
        \addplot[only marks, mark options={scale=0.4}] table [x=value, y=engineering, col sep=comma] {data/data/total_degree_distribution/value_counts.csv};
        \addplot[red, ultra thick, no marks, domain={10:3000}, forget plot] {10^(6.5)*x^(-2.9849)};
      \end{axis}
    \end{tikzpicture}}
  \end{subfigure}
  \hfill
  \begin{subfigure}[b]{0.24\textwidth} \resizebox{\textwidth}{!}{
    \begin{tikzpicture}
      \begin{axis} [
          xmin=1, xmax=600, ymin=10^(-0.2), ymax=10^(5.25),
          clip mode=individual,
          smooth, xmode=log, ymode=log,
          xlabel=vertex degree,
          ylabel style={
            yshift=-0.6cm
          },
          legend style={at={(0,0)}, anchor=south west, at={(axis description cs:0.025,0.025)}}, legend cell align={left},
        ]
        \addplot[only marks, mark options={scale=0.4}] table [x=value, y=mathematics, col sep=comma] {data/data/total_degree_distribution/value_counts.csv};
        \addplot[red, ultra thick, no marks, domain={10:3000}] {10^(6.3)*x^(-2.7897)};
      \end{axis}
    \end{tikzpicture}}
  \end{subfigure}
  \hfill
  \begin{subfigure}[b]{0.24\textwidth} \resizebox{\textwidth}{!}{
    \begin{tikzpicture}
      \begin{axis} [
          xmin=1, xmax=220, ymin=10^(-0.2), ymax=10^(4.5),
          clip mode=individual,
          smooth, xmode=log, ymode=log,
          xlabel=vertex degree,
          ylabel style={
            yshift=-0.6cm
          },
          legend style={at={(0,0)}, anchor=south west, at={(axis description cs:0.025,0.025)}}, legend cell align={left},
        ]
        \addplot[only marks, mark options={scale=0.4}] table [x=value, y=statistics, col sep=comma] {data/data/total_degree_distribution/value_counts.csv};
        \addplot[red, ultra thick, no marks, domain={10:3000}] {10^(5.85)*x^(-2.9635)};
      \end{axis}
    \end{tikzpicture}}
  \end{subfigure}
  \hspace*{\fill}
  \\
  \hspace*{\fill}
  \begin{subfigure}[b]{0.24\textwidth} \resizebox{\textwidth}{!}{
    \begin{tikzpicture}
      \begin{axis} [
          xmin=1, xmax=300, ymin=10^(-0.2), ymax=10^(6),
          clip mode=individual,
          smooth, xmode=log, ymode=log,
          xlabel=edge degree, ylabel=frequency,
          ylabel style={
            yshift=-0.4cm
          },
          legend style={at={(0,0)}, anchor=south west, at={(axis description cs:0.025,0.025)}}, legend cell align={left},
        ]
        \addplot[only marks, mark options={scale=0.4}] table [x=value, y=computer_science, col sep=comma] {data/data/ho_degree_distribution_1/value_counts.csv};
        \addplot[red, ultra thick, no marks, domain={10:300}] {10^(8.37)*x^(-3.7591)};
      \end{axis}
    \end{tikzpicture}} \caption{cs}
  \end{subfigure}
  \hfill
  \begin{subfigure}[b]{0.24\textwidth} \resizebox{\textwidth}{!}{
    \begin{tikzpicture}
      \begin{axis} [
          xmin=1, xmax=110, ymin=10^(-0.2), ymax=10^(5.5),
          clip mode=individual,
          smooth, xmode=log, ymode=log,
          xlabel=edge degree,
          ylabel style={
            yshift=-0.6cm
          },
          legend style={at={(0,0)}, anchor=south west, at={(axis description cs:0.025,0.025)}}, legend cell align={left},
        ]
        \addplot[only marks, mark options={scale=0.4}] table [x=value, y=engineering, col sep=comma] {data/data/ho_degree_distribution_1/value_counts.csv};
        \addplot[red, ultra thick, no marks, domain={10:300}] {10^(7.9)*x^(-4.1374)};
      \end{axis}
    \end{tikzpicture}} \caption{eess}
  \end{subfigure}
  \hfill
  \begin{subfigure}[b]{0.24\textwidth} \resizebox{\textwidth}{!}{
    \begin{tikzpicture}
      \begin{axis} [
          xmin=1, xmax=80, ymin=10^(-0.2), ymax=10^(5.5),
          clip mode=individual,
          smooth, xmode=log, ymode=log,
          xlabel=edge degree,
          ylabel style={
            yshift=-0.6cm
          },
          legend style={at={(0,0)}, anchor=south west, at={(axis description cs:0.025,0.025)}}, legend cell align={left},
        ]
        \addplot[only marks, mark options={scale=0.4}] table [x=value, y=mathematics, col sep=comma] {data/data/ho_degree_distribution_1/value_counts.csv};
        \addplot[red, ultra thick, no marks, domain={10:300}] {10^(7.5)*x^(-4.4681)};
      \end{axis}
    \end{tikzpicture}} \caption{math}
  \end{subfigure}
  \hfill
  \begin{subfigure}[b]{0.24\textwidth} \resizebox{\textwidth}{!}{
    \begin{tikzpicture}
      \begin{axis} [
          xmin=1, xmax=70, ymin=10^(-0.2), ymax=10^(5.0),
          clip mode=individual,
          smooth, xmode=log, ymode=log,
          xlabel=edge degree,
          ylabel style={
            yshift=-0.6cm
          },
          legend style={at={(0,0)}, anchor=south west, at={(axis description cs:0.025,0.025)}}, legend cell align={left},
        ]
        \addplot[only marks, mark options={scale=0.4}] table [x=value, y=statistics, col sep=comma] {data/data/ho_degree_distribution_1/value_counts.csv};
        \addplot[red, ultra thick, no marks, domain={10:300}] {10^(8)*x^(-4.8626)};
      \end{axis}
    \end{tikzpicture}} \caption{stat}
  \end{subfigure}
  \hspace*{\fill}
  \caption{Vertex-degree distributions (top) and edge-degree distributions (bottom) of the data sets}
  \label{ho_degree_distributions}
\end{figure}
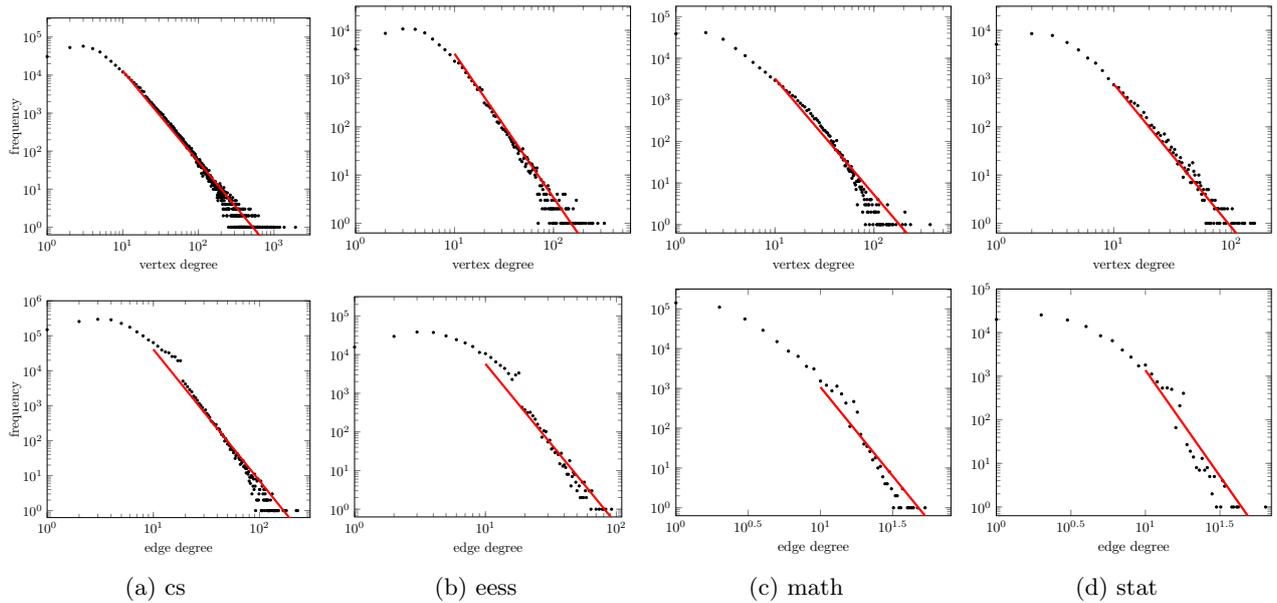

% fitting gamma, exponents
To fit the ADRCM to the datasets, we need to set two model parameters.
First, we can use Theorem~\ref{thm:gen_deg} to estimate the parameter $\g$ describing the datasets based on their vertex or higher-order degree distributions.
The vertex and edge-degree distributions are visualized in Figure~\ref{ho_degree_distributions}.

% explain drop in the distributions
All plots exhibit a drop in the empirical distributions at the value 20.
This is explained by the exclusion of documents with more than 21 authors.
Due to combinatorial reasons, this discontinuity is more pronounced in the case of the edge-degree distributions.
For heavier tailed distributions, a larger number of documents have more, than 21 authors, leading to a greater drop for the cs and eess datasets.
As explained in the beginning of Section~\ref{ss:dat}, including thesedocuments would lead to the problem of the high influence of a few high-dimensional interactions as described earlier.

\subsection{Higher-order degree distributions}

% fitting the power-law distribution
Just as in the case of the simulations, fitting the parameters of the power-law distribution poses computational challenges once again.
When determining the minimum value $x_{\textrm{min}}$ from which the power law is visible, the goal is to find a balance between two conflicting interests.
\begin{itemize}
\item{
On the one hand, choosing a low minimum degree value would ensure enough data points in the degree distributions so that the estimates of the exponent is less noisy. 
}
\item{
On the other hand, choosing a high minimum degree value would remove the noise from light-tailed components of the degree.
}
\end{itemize}
Considering both effects, we found that setting the minimum value $x_{\textrm{min}} = 10$ is a good compromise for fitting the power-law distributions. We note that this choice is more conservative than the one used in Section~\ref{sec:sim}, where we set $x_{\textrm{min}} = 30$. This is because we found the datasets to be more noisy than the simulated networks. Hence, we chose to use a smaller minimum value to enlarge the number of data points used for the fitting. After fitting the power-law distributions, we can use Theorem~\ref{thm:gen_deg} to infer the $\g$ model parameter based on the fitted power-law exponents.

The fitted exponents and the $\g$ model parameters inferred from these exponents are summarized in Table~\ref{dataset_exponents}.
\begin{table} [h] \centering \caption{Fitted exponents of the degree distributions and the inferred $\g$ model parameters} \label{dataset_exponents}
  \begin{tabular}{|l|r|r|r|r|} \hline
    % exponent = m - (m+1) / gamma - 1
    % inferred gamma = (m+1) / (m-1 - exponent)
    \multirow{2}{*}{dataset } & \multicolumn{2}{|c|}{vertex degree} & \multicolumn{2}{|c|}{edge degree} \\ %\cline{2-5}
                     &  exponent & inferred $\g$ &  exponent & inferred $\g$ \\ \hline
    cs & -2.39 &      0.72 & -3.76 &      0.53 \\
    eess      & -2.98 &      0.50 & -4.14 &      0.48 \\
    math      & -2.79 &      0.56 & -4.47 &      0.45 \\
    stat       & -2.96 &      0.51 & -4.86 &      0.41 \\ \hline
  \end{tabular}
\end{table}
In general, the edge-degree distributions have a thinner tail compared to that of the related vertex-degree distributions.
We can see that the parameter $\g$ differs substantially when inferred from the vertex- and edge-degree distributions, respectively.
This observation shows that the ADRCM, with the connection kernel we apply it, is not flexible enough to capture the binary and the higher-order features at the same time.
We henceforth infer $\g$ from the {vertex-degree distributions} due to the following reasons.
First, they are less affected by the high-dimensional interactions: a document with $n$ authors contributes with $n$ values in case of the vertex degrees, while it is represented by ${n \choose 2}$ values in the edge-degree distribution.
Additionally, the computation of the vertex-degree distribution only requires the consideration of pairwise relationships, which the original ADRCM was designed to describe.

% fitting beta, edge density
As discussed by~\cite{glm2}, the parameter $\b$ governs the asymptotic edge density (the expected number of edges containing a vertex) of the generated networks through the formula $\E[d_{0,1}]= \b / (1 - \g)$.
Thus, using the above formula, we estimate the parameter $\b$ from the mean vertex degree of the datasets.
The mean vertex degrees and the estimated parameter $\hat \b$ are shown in Table~\ref{dataset_betas}.
\begin{table} [h] \centering \caption{Mean vertex degree \& $\hat\b$} \label{dataset_betas}
  \begin{tabular}{|l|r|r|} \hline
    dataset & mean vertex degree & $\hat \b$ \\ \hline
    cs      &               9.57 &      2.69 \\
    eess    &               7.13 &      3.54 \\
    math    &               4.58 &      2.02 \\
    stat    &               5.14 &      2.52 \\ \hline
  \end{tabular}
\end{table}

\begin{figure} [h!] \centering
  \begin{subfigure}[b]{0.24\textwidth}
    \includegraphics[width=\textwidth]{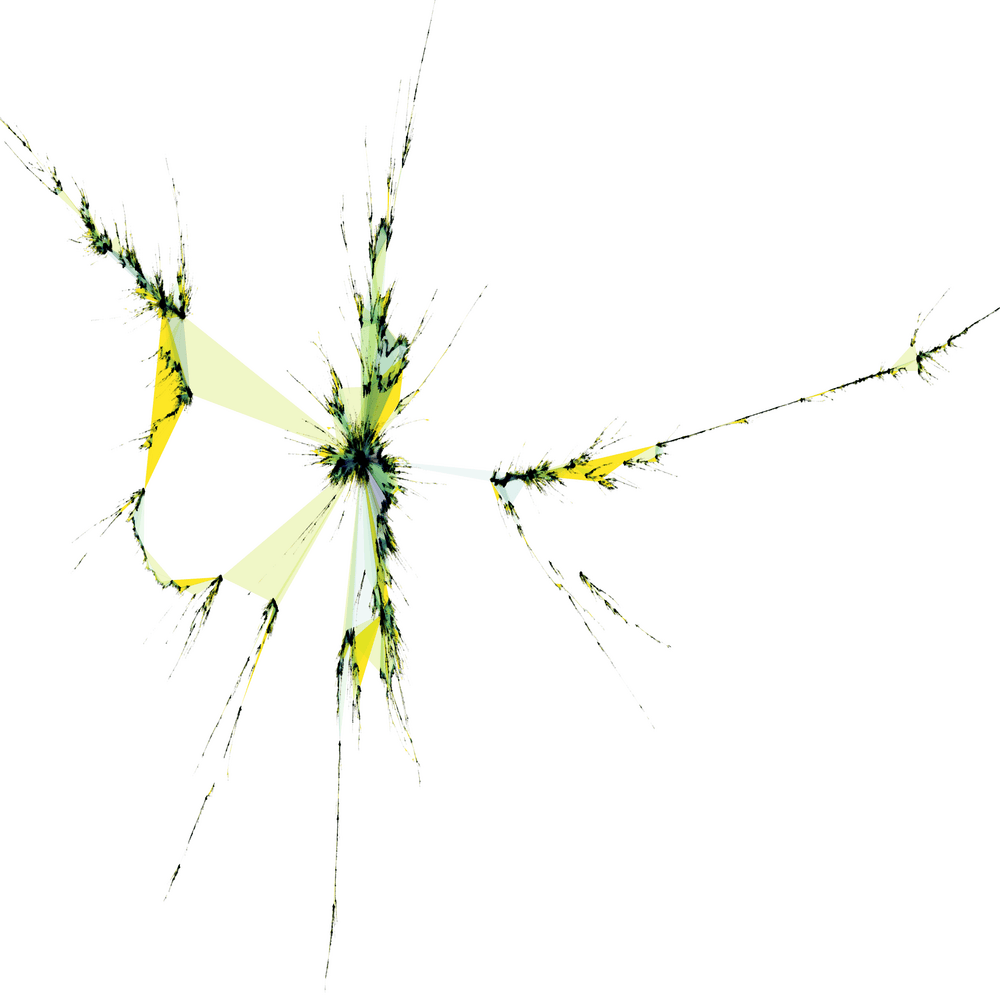} \caption{cs}
  \end{subfigure}
  \hfill
  \begin{subfigure}[b]{0.24\textwidth}
    \includegraphics[width=\textwidth]{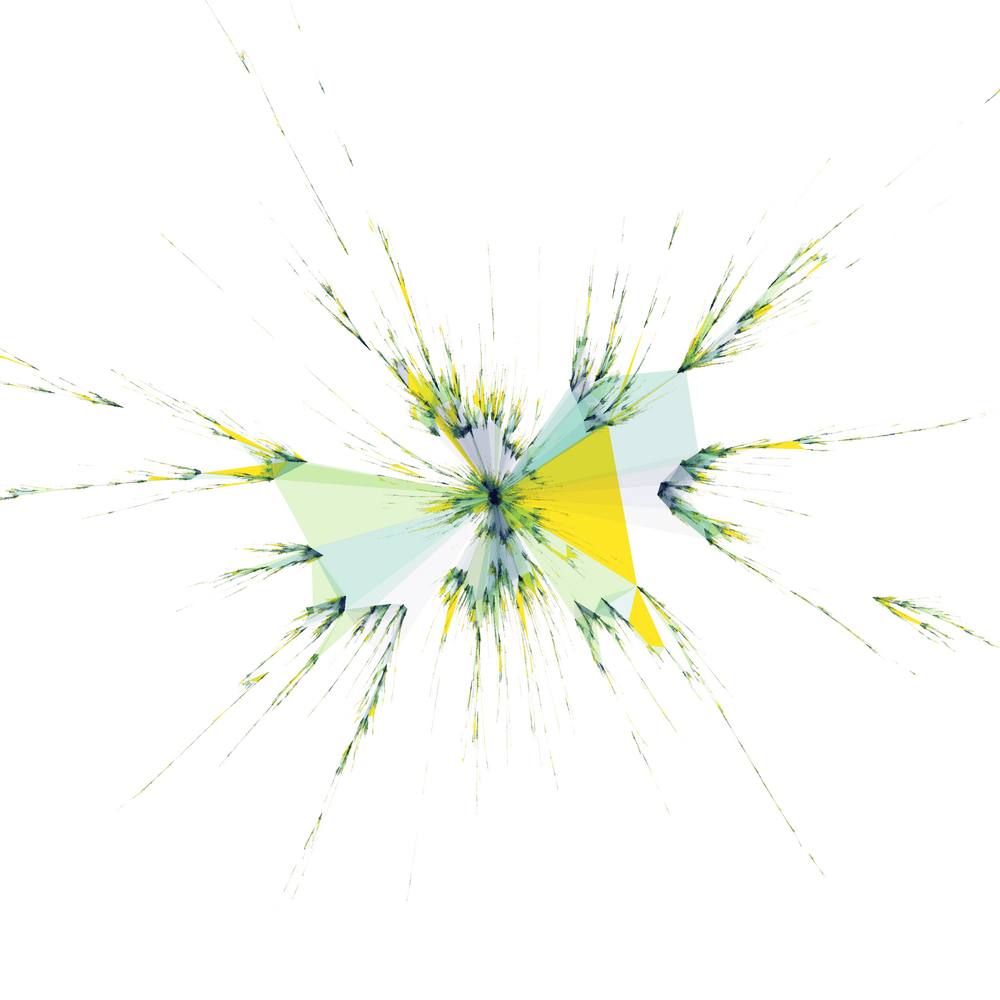}      \caption{eess}
  \end{subfigure}
  \hfill
  \begin{subfigure}[b]{0.24\textwidth}
    \includegraphics[width=\textwidth]{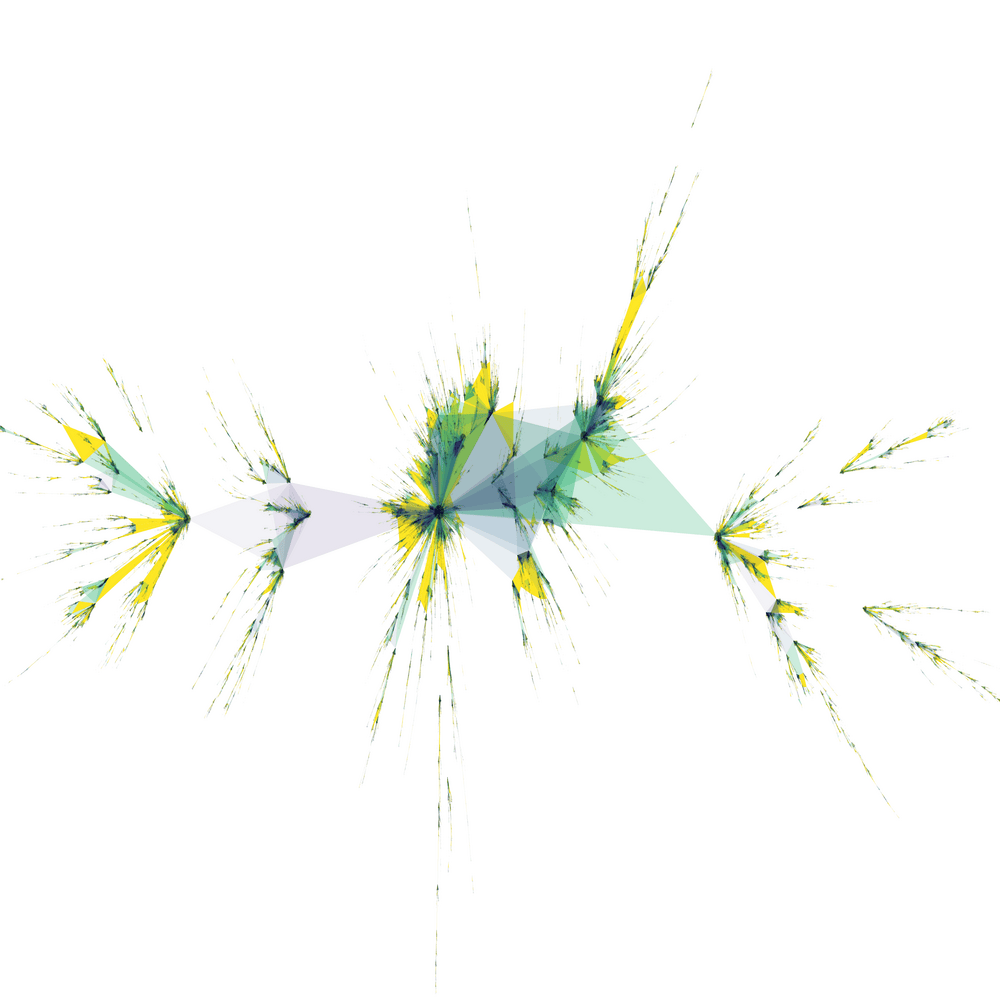}      \caption{math}
  \end{subfigure}
  \hfill
  \begin{subfigure}[b]{0.24\textwidth}
    \includegraphics[width=\textwidth]{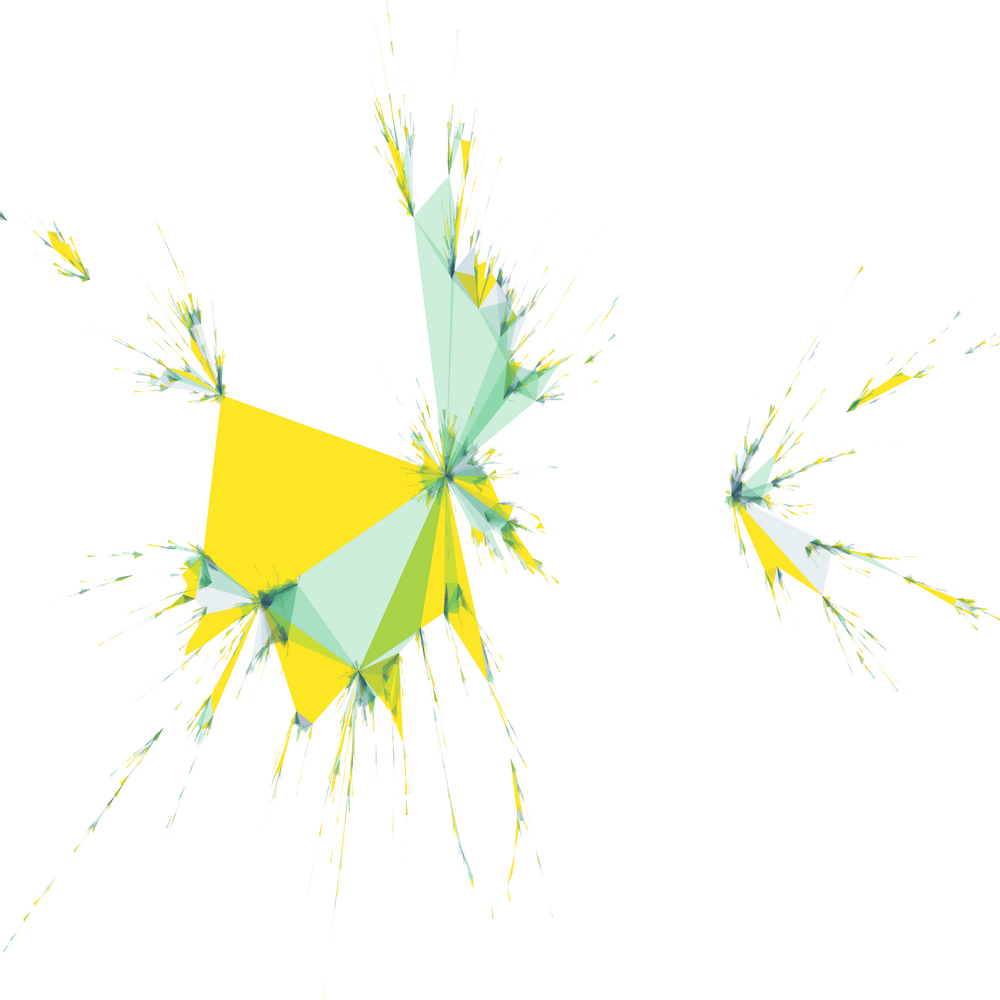}       \caption{stat}
  \end{subfigure}
  \caption{
  The largest component of simulated ADRCMs with fitted parameters.
 }
  \label{network_plot_samples_of_adrcm}
\end{figure}

% sample networks
After fitting the model parameters, we can generate synthetic networks.
We simulated a representative network for each datasets, whose largest components are visualized in Figure~\ref{network_plot_samples_of_adrcm}.
When comparing with the plots of the actual datasets, we observe that although the ADRCM is capable of generating triangles and tetrahedrons, have a tendency to produce globally tree-like structures.

\FloatBarrier

%
%SS TRI
%
\subsection{Triangle counts}
\label{ss:tri}

% triangle count, hypothesis tests
Next, we examine if the simplex counts in the ADRCM match those in the datasets.
The simplex counts of the datasets are presented in Table~\ref{simplex_counts_table}.
\begin{table}[!htb] \centering \caption{Number of simplices of different dimensions in the datasets} \label{simplex_counts_table}
  \begin{tabular}{|l|r|r|r|} \hline
    dataset & vertices &       edges &   triangles \\ \hline
    cs      & 433\,244 & 2\,073\,235 & 4\,055\,220 \\
    eess    &  77\,686 &    276\,947 &    562\,382 \\
    math    & 198\,601 &    455\,130 &    321\,406 \\
    stat    &  44\,380 &    114\,003 &    135\,800 \\ \hline
  \end{tabular}
\end{table}

\begin{table}[!htb] \centering \caption{Estimated parameters of the stable distributions for triangle counts} \label{triangle_count_stable_distribution_parameters}
    \begin{tabular}{|l|r|r|r|r|r|r|} \hline
      dataset & $\hat\alpha$ & $\hat\b$ &     location &    scale \\ \hline
      cs      &         1.39 &      1.0 & 18\,785\,263 & 504\,582 \\
      eess    &         1.98 &      1.0 &  1\,911\,396 &  38\,527 \\
      math    &         1.79 &      1.0 &  2\,027\,542 &  28\,774 \\
      stat    &         1.96 &      1.0 &     566\,665 &  15\,352 \\ \hline
    \end{tabular}
\end{table}

The number of vertices is matched by the model on expectation as we choose the size of the sampling window accordingly.
It is also irrelevant to examine the edge count, being asymptotically fixed through the parameter $\b$.
Consequently, the first nontrivial dimension to consider is the triangle count.

% conjecture: stable distribution of triangle counts
As shown in Theorem~\ref{thm:ass}, the number of edges follows a stable distribution if the model parameter $\g$ is larger than $0.5$, which is the case for our datasets.
We conjecture that the distributions of higher-dimensional simplex counts also follow a stable distribution for $\g > 0.5$.

% simulation, fitting of stable distribution parameters
To study the distribution of the triangle counts, we simulated $100$ networks with estimated parameters $\hat\b$ and $\hat\gamma$ determined according to the datasets.
For fitting stable distributions to the triangle counts, we use the method detailed in Section~\ref{sec:sim}: while the parameters $\alpha$ and $\b$ of the stable distribution are predicted based on our mathematical conjecture, the location and scale parameters are estimated via maximum likelihood.
The fitted parameters are presented in Table~\ref{triangle_count_stable_distribution_parameters}.

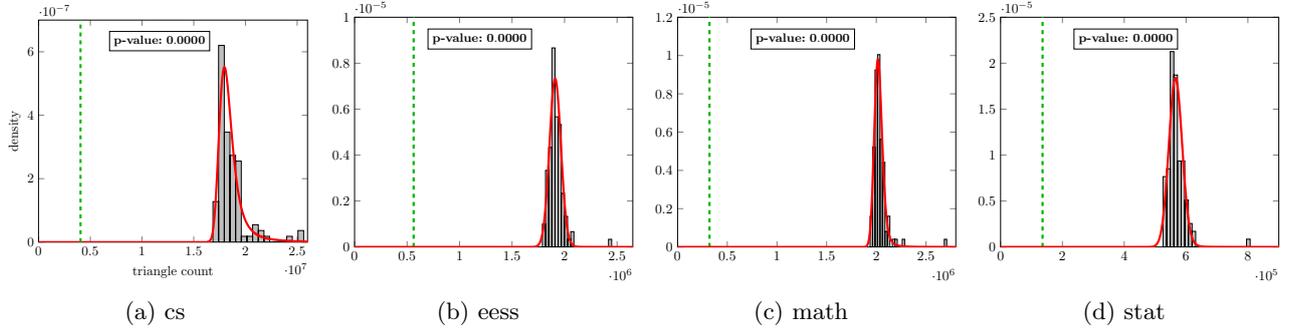
\begin{figure} [h!] \centering
  \hfill
  \begin{subfigure}[b]{0.24\textwidth} \resizebox{\textwidth}{!}{
      \begin{tikzpicture}
        \begin{groupplot}[
            group style={
              group size=1 by 1,
              xticklabels at=edge bottom,
              horizontal sep=0pt,
            },
            height=7.5cm,
            ymin=0, ymax=7e-7,
          ]

          \nextgroupplot[
            xmin=0, xmax=26000000,
            axis line style={-},
            xlabel=triangle count,
            ylabel=density,
            ylabel style={
              yshift=-0.6cm
            },
          ]
          
          \addplot[ybar interval, mark=no, fill=gray!50!white] table [x=computer_science_bin_left_limit, y=computer_science_value, col sep=comma] {data/hypothesis_test/num_of_triangles_stable/linear_histograms.csv};
          \addplot[ultra thick, red] table [x=computer_science_value, y=computer_science_pdf, col sep=comma] {data/hypothesis_test/num_of_triangles_stable/theoretical_pdf.csv};
          \addplot[ultra thick, green!70!black, mark=none, dashed] coordinates {(4055220, 0) (4055220, 1e-6)};
          \node[anchor=north, draw] at (rel axis cs: 0.45, 0.95) {\small{\textbf{p-value:} $\mathbf{0.0000}$}};
        \end{groupplot}
    \end{tikzpicture}} \caption{cs}
  \end{subfigure}
  \hfill
  \begin{subfigure}[b]{0.24\textwidth} \resizebox{\textwidth}{!}{
      \begin{tikzpicture}
        \begin{groupplot}[
            group style={
              group size=1 by 1,
              xticklabels at=edge bottom,
              horizontal sep=0pt,
            },
            height=7.5cm,
            ymin=0, ymax=1e-5,
          ]

          \nextgroupplot[
            xmin=0, xmax=2650000,
            axis line style={-},
		xlabel=\phantom a,
            ylabel style={
              yshift=-0.6cm
            },
          ]
          
          \addplot[ultra thick, green!70!black, mark=none, dashed] coordinates {(562382, 0) (562382, 1e-5)};
          \addplot[ybar interval, mark=no, fill=gray!50!white] table [x=engineering_bin_left_limit, y=engineering_value, col sep=comma] {data/hypothesis_test/num_of_triangles_stable/linear_histograms.csv};
          \addplot[ultra thick, red, domain=1700000:2500000] table [x=engineering_value, y=engineering_pdf, col sep=comma] {data/hypothesis_test/num_of_triangles_stable/theoretical_pdf.csv};
          \node[anchor=north, draw] at (rel axis cs: 0.45, 0.95) {\small{\textbf{p-value:} $\mathbf{0.0000}$}};
        \end{groupplot}
    \end{tikzpicture}} \caption{eess}
  \end{subfigure}
  \hfill
  \begin{subfigure}[b]{0.24\textwidth} \resizebox{\textwidth}{!}{
      \begin{tikzpicture}
        \begin{groupplot}[
            group style={
              group size=1 by 1,
              xticklabels at=edge bottom,
              horizontal sep=0pt,
            },
            height=7.5cm,
            ymin=0, ymax=1.2e-5,
          ]

          \nextgroupplot[
            xmin=0, xmax=2800000,
            axis line style={-},
		xlabel=\phantom a,
            ylabel style={
              yshift=-0.6cm
            },
          ]

          \addplot[ultra thick, green!70!black, mark=none, dashed] coordinates {(321406, 0) (321406, 1e-4)};
          \addplot[ybar interval, mark=no, fill=gray!50!white,] table [x=mathematics_bin_left_limit, y=mathematics_value, col sep=comma] {data/hypothesis_test/num_of_triangles_stable/linear_histograms.csv};
          \addplot[ultra thick, red] table [x=mathematics_value, y=mathematics_pdf, col sep=comma] {data/hypothesis_test/num_of_triangles_stable/theoretical_pdf.csv};
          \node[anchor=north, draw] at (rel axis cs: 0.45, 0.95) {\small{\textbf{p-value:} $\mathbf{0.0000}$}};
        \end{groupplot}
    \end{tikzpicture}} \caption{math}
  \end{subfigure}
  \hfill
  \begin{subfigure}[b]{0.24\textwidth} \resizebox{\textwidth}{!}{
      \begin{tikzpicture}
        \begin{groupplot}[
            group style={
              group size=1 by 1,
              xticklabels at=edge bottom,
              horizontal sep=0pt,
            },
            height=7.5cm,
            ymin=0, ymax=2.5e-5,
          ]

          \nextgroupplot[
            xmin=0, xmax=900000,
            axis line style={-},
		xlabel=\phantom a,
            ylabel style={
              yshift=-0.6cm
            },
          ]

          \addplot[ultra thick, green!70!black, mark=none, dashed] coordinates {(135800, 0) (135800, 1e-4)};
          \addplot[ybar interval, mark=no, fill=gray!50!white,] table [x=statistics_bin_left_limit, y=statistics_value, col sep=comma] {data/hypothesis_test/num_of_triangles_stable/linear_histograms.csv};
          \addplot[ultra thick, red]      table [x=statistics_value, y=statistics_pdf, col sep=comma] {data/hypothesis_test/num_of_triangles_stable/theoretical_pdf.csv};
          \node[anchor=north, draw] at (rel axis cs: 0.45, 0.95) {\small{\textbf{p-value:} $\mathbf{0.0000}$}};
        \end{groupplot}
    \end{tikzpicture}} \caption{stat}
  \end{subfigure}
  \hfill
  \caption{    Stable distribution and hypothesis testing of the triangle counts for the datasets.
    The model parameters were determined based on the parameters of the datasets.
  }
  \label{simplex_count_hypothesis_testing}
\end{figure}

% hypothesis tests
After empirically verifying the simplex-count distribution, we conduct a hypothesis test based on the triangle counts.
Our null model is the ADRCM model with the connection kernel from Section~\ref{sec:mod}.
The dataset values are marked by vertical green dashed lines in Figure~\ref{simplex_count_hypothesis_testing}.
We conclude that in all cases, the ADRCM contains substantially more triangles compared to the dataset.
In particular, the null hypothesis is rejected at the 5\% level.

\FloatBarrier
\subsection{Betti numbers}

The presence of loops is important features of collaboration networks as they quantify its interconnectedness.
In Section~\ref{sec:sim}, we provided numerical evidence for the conjecture that the Betti numbers follow a stable distribution if $\g > 0.5$.
On this basis, we can conduct a similar hypothesis test as above on the first Betti numbers, where we use the ADRCM as the null model.

% Betti numbers of datasets
The Betti numbers of the datasets we aim to test are presented in Table~\ref{betti_numbers_table}.
\begin{table}[!htb] \centering \caption{Betti numbers of the datasets} \label{betti_numbers_table}
      \begin{tabular}{|l|r|r|} \hline
        dataset & Betti-0 &  Betti-1 \\ \hline
        cs      & 22\,576 & 168\,770 \\
        eess    &  5\,533 &   7\,419 \\
        math    & 26\,197 &  78\,009 \\
        stat    &  4\,049 &   7\,275 \\ \hline
      \end{tabular}
\end{table}

\begin{table}[!htb] \centering \caption{Parameter estimates of the stable distributions for Betti-1} \label{betti_number_stable_distribution_parameters}
    \begin{tabular}{|l|r|r|r|r|} \hline
      dataset & $\hat\a$ & $\hat\b$ & location & scale \\ \hline
      cs      &     1.39 &     -1.0 &       37 & 12.83 \\
      eess    &     1.98 &     -1.0 &      105 &  9.66 \\
      math    &     1.79 &     -1.0 &      490 & 19.16 \\
      stat    &     1.96 &     -1.0 &      126 &  8.99 \\ \hline
    \end{tabular}
\end{table}

\begin{figure} [h] \centering
  \hfill
  \begin{subfigure}[b]{0.24\textwidth} \resizebox{\textwidth}{!}{
      \begin{tikzpicture}
        \begin{groupplot}[
            group style={
              group size=2 by 1,
              xticklabels at=edge bottom,
              horizontal sep=0pt,
            },
            height=7.5cm,
            ymin=0, ymax=0.028,
          ]

          \nextgroupplot[
            xmin=-10, xmax=110,
            xtick={0, 50, 100},
            axis y line=left,
            axis line style={-},
            xlabel=Betti-1,
            xlabel style={
              xshift=0.75cm
            },
            ylabel=density,
            ylabel style={
              yshift=-0.6cm
            },
            xtick scale label code/.code={},
            width=7.0cm,
          ]
          \addplot[ybar interval, mark=no, fill=gray!50!white,] table [x=computer_science_bin_left_limit, y=computer_science_value, col sep=comma] {data/hypothesis_test/betti_number_1_stable/linear_histograms.csv};
          \addplot[ultra thick, red]      table [x=computer_science_value, y=computer_science_pdf, col sep=comma] {data/hypothesis_test/betti_number_1_stable/theoretical_pdf.csv};

          \nextgroupplot[
            xmin=168769, xmax=168771,
            axis x discontinuity=parallel,
            xtick={168770},
		xlabel=\phantom a,
            xticklabels={$168\,770$},
            scaled x ticks=false,
            axis y line=right,
            axis line style={-},
            width=3.5cm
          ]
          \addplot[ultra thick, green!70!black, mark=none, dashed] coordinates {(168770, 0) (168770, 0.04)};
        \end{groupplot}
        \node[anchor=north east, draw] at (3.25, 5.5) {\small{\textbf{p-value:} $\mathbf{0.0000}$}};
    \end{tikzpicture}} \caption{cs}
  \end{subfigure}
  \hfill
  \begin{subfigure}[b]{0.22\textwidth} \resizebox{\textwidth}{!}{
      \begin{tikzpicture}
        \begin{groupplot}[
            group style={
              group size=2 by 1,
              xticklabels at=edge bottom,
              horizontal sep=0pt,
            },
            height=7.5cm,
            ymin=0, ymax=0.0375,
          ]

          \nextgroupplot[
            xmin=-10, xmax=150,
            xtick={0, 50, 100},
            axis y line=left,
            axis line style={-},
		xlabel=\phantom a,
            xlabel style={
              xshift=0.75cm
            },
            ylabel style={
              yshift=-0.6cm
            },
            xtick scale label code/.code={},
            width=7.0cm,
          ]
          \addplot[ybar interval, mark=no, fill=gray!50!white,] table [x=engineering_bin_left_limit, y=engineering_value, col sep=comma] {data/hypothesis_test/betti_number_1_stable/linear_histograms.csv};
          \addplot[ultra thick, red]      table [x=engineering_value, y=engineering_pdf, col sep=comma] {data/hypothesis_test/betti_number_1_stable/theoretical_pdf.csv};

          \nextgroupplot[
            xmin=7418, xmax=7420,
            axis x discontinuity=parallel,
            xtick={7419},
            axis y line=right,
            axis line style={-},
            width=3.5cm
          ]
          \addplot[ultra thick, green!70!black, mark=none, dashed] coordinates {(7419, 0) (7419, 0.05)};
        \end{groupplot}
        \node[anchor=north east, draw] at (3.25, 5.5) {\small{\textbf{p-value:} $\mathbf{0.0000}$}};
    \end{tikzpicture}} \caption{eess}
  \end{subfigure}
  \hfill
  \begin{subfigure}[b]{0.24\textwidth} \resizebox{\textwidth}{!}{
      \begin{tikzpicture}
        \begin{groupplot}[
            group style={
              group size=2 by 1,
              xticklabels at=edge bottom,
              horizontal sep=0pt,
            },
            height=7.5cm,
            ymin=0, ymax=0.0175,
          ]

          \nextgroupplot[
            xmin=150, xmax=600,
            xtick={200, 300, 400, 500},
		xlabel=\phantom a,
            axis y line=left,
            axis line style={-},
            xlabel style={
              xshift=0.75cm
            },
            ylabel style={
              yshift=-0.6cm
            },
            xtick scale label code/.code={},
            width=7.0cm,
          ]
          \addplot[ybar interval, mark=no, fill=gray!50!white,] table [x=mathematics_bin_left_limit, y=mathematics_value, col sep=comma] {data/hypothesis_test/betti_number_1_stable/linear_histograms.csv};
          \addplot[ultra thick, red]      table [x=mathematics_value, y=mathematics_pdf, col sep=comma] {data/hypothesis_test/betti_number_1_stable/theoretical_pdf.csv};

          \nextgroupplot[
            xmin=78008, xmax=78010,
            axis x discontinuity=parallel,
            scaled x ticks=false,
            xtick={78009},
		xlabel=\phantom a,
            axis y line=right,
            axis line style={-},
            width=3.5cm
          ]
          \addplot[ultra thick, green!70!black, mark=none, dashed] coordinates {(78009, 0) (78009, 0.05)};
        \end{groupplot}
        \node[anchor=north east, draw] at (3.25, 5.5) {\small{\textbf{p-value:} $\mathbf{0.0000}$}};
    \end{tikzpicture}} \caption{math}
  \end{subfigure}
  \hfill
  \begin{subfigure}[b]{0.22\textwidth} \resizebox{\textwidth}{!}{
      \begin{tikzpicture}
        \begin{groupplot}[
            group style={
              group size=2 by 1,
              xticklabels at=edge bottom,
              horizontal sep=0pt,
            },
            height=7.5cm,
            ymin=0, ymax=0.035,
          ]

          \nextgroupplot[
            xmin=-10, xmax=180,
            xtick={0, 50, 100, 150},
            axis y line=left,
		xlabel=\phantom a,
            axis line style={-},
            xlabel style={
              xshift=0.75cm
            },
            ylabel style={
              yshift=-0.6cm
            },
            xtick scale label code/.code={},
            width=7.0cm,
          ]
          \addplot[ybar interval, mark=no, fill=gray!50!white,] table [x=statistics_bin_left_limit, y=statistics_value, col sep=comma] {data/hypothesis_test/betti_number_1_stable/linear_histograms.csv};
          \addplot[ultra thick, red]      table [x=statistics_value, y=statistics_pdf, col sep=comma] {data/hypothesis_test/betti_number_1_stable/theoretical_pdf.csv};

          \nextgroupplot[
            xmin=7274, xmax=7276,
            axis x discontinuity=parallel,
            scaled x ticks=false,
            xtick={7275},
            axis y line=right,
            axis line style={-},
            width=3.5cm
          ]
          \addplot[ultra thick, green!70!black, mark=none, dashed] coordinates {(7275, 0) (7275, 0.05)};
        \end{groupplot}
        \node[anchor=north east, draw] at (3.25, 5.5) {\small{\textbf{p-value:} $\mathbf{0.0000}$}};
    \end{tikzpicture}} \caption{stat}
  \end{subfigure}
  \hfill
  \caption{Hypothesis testing of Betti-1 for the datasets.}
  \label{betti_number_hypothesis_testing}
\end{figure}
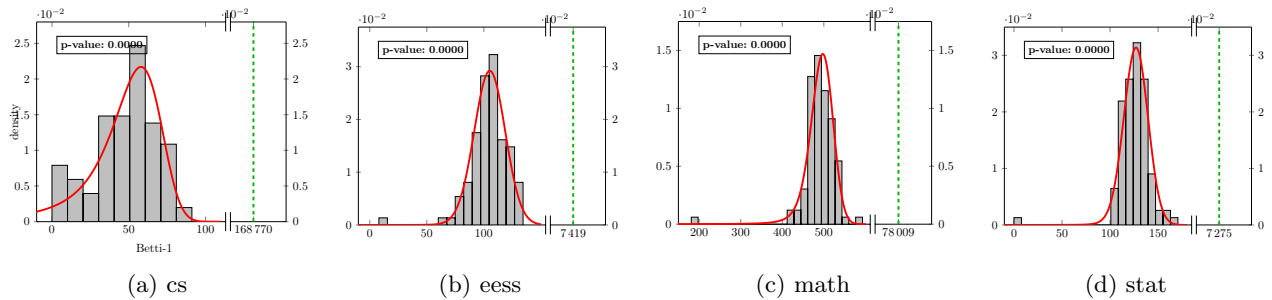

% fitting the stable distribution
We again simulate 100 networks using the ADRCM with the fitted model parameters.
As in Section~\ref{ss:tri}, the parameters $\alpha$ and $\b$ of the stable distributions are predicted by our conjecture, while the location and scale parameters are fitted via maximum likelihood.
After fitting the stable distributions, we visualize the hypothesis testing in Figure~\ref{betti_number_hypothesis_testing}.
The parameters of the considered stable distributions are given in Table~\ref{betti_number_stable_distribution_parameters}.
In particular, the real datasets contain a significantly greater amount of loops than the networks generated by the ADRCM, thus the null hypothesis is rejected.

As explained by~\cite{patania}, the loops in the network can be interpreted as bridges between communities.
Hence, they are important features of scientific collaboration networks.
Our analysis indicates that this community structure features a rich spatial correlation pattern, which cannot be captured fully by a simple two-parameter model such as the ADRCM.

We believe that the reason for this phenomena is at least partly due to the chosen connection kernel.
The vertices connect to many vertices within their neighborhood with probability 1, thereby making it difficult to form loops. 
This suggests that the ADRCM generates network that appear tree-like on a global level with relatively few large loops.

To illustrate this idea, we carried out a pilot study, where we examined the influence of the connection kernel on the first Betti numbers.
More precisely, we can employ the more general connection kernel where two vertices $(x, u), (y, v) \in \PP$ with $u \le v$ connect with probability $1 / (2a)$ ($a \geq 1/2$), whenever $|x - y| \le a \, \b u^{-\g} v^{\g - 1}$~\citep{glm2}.
Note that for $a = 0.5$, the connection kernel coincides with the one introduced in Section~\ref{sec:mod}.
Increasing the newly introduced model parameter $a$ increases the distance of the vertices in which connections can be established.
On the other hand, to keep the expected number of connections of the vertices intact, it simultaneously reduces the connection probability.

For $\b = 1$ and $\g = 0.6$, we simulated six sets of $100$ networks each, with a network size of $100\,000$.
We then gradually increased the value of the parameter $a$ from the default value of $0.5$, and kept track of the increase of the first Betti numbers.
The results are shown in Table~\ref{influence_of_a}, and we conclude that even a slight increase of $a$ results in a drastic growth of the first Betti numbers.
\begin{table} [h!] \centering \caption{Influence of the profile function on Betti-1} \label{influence_of_a}
  \begin{tabular}{|r|r|r|r|r|r|r|} \hline
    parameter $a$            &   0.5 &  0.6 & 0.7 & 0.8 & 0.9 & 1.0 \\ \hline
    mean of Betti-1 & $170$ & $4\,873$ & $10\,976$ & $17\,786$ & $24\,914$ & $31\,961$ \\ \hline
  \end{tabular}
\end{table}

\FloatBarrier
\section{Conclusion and outlook}
\label{sec:conc}

To analyze higher-order network structures, we investigated the ADRCM as a clique complex.

First, we examined how the neighborhood of simplices of different dimensions are organized and proved that the higher-order degree distributions have a power-law tail in the limit for large networks.
Next, we proved that in the limit for large networks, the  recentered and suitably rescaled edge count follows a normal
distribution if the model parameter $\g$ is less than $0.5$, and a stable distribution for $\g > 0.5$. Turning our attention to the topological features, a CLT was proved for the Betti numbers if $\g < 0.25$.
Recognizing the limitations of the ADRCM model, we devised a ``thinning'' procedure where certain types of edges are removed independently with a given thinning probability.
This provided us with the possibility to adjust the edge degree exponents, while keeping the power-law exponent of the vertex degree distribution intact.

To show that the above theoretical results can be used in real-world data sets, we examined the extent to which the theorems are valid for finite networks by simulating several networks using identical model parameters.
We found that the convergence of specific quantities to their limiting behavior is already clearly visible in networks of reasonable size.
Furthermore, we also provided numerical evidence supporting our conjectures regarding the stable distribution of the Betti numbers when $\g > 0.5$.

Finally, after showing that the theoretical results are applicable to networks of finite size, we analyzed real-world scientific collaboration networks from arXiv.
Following an exploratory analysis of these higher-order collaboration networks, we fitted the model parameters  to the data.
Developing hypothesis tests, we showed that -- although several properties are well described by the higher-order ADRCM --, topologically important quantities, such as Betti numbers or the higher-dimensional simplex counts, are not well explained.
Looking ahead, we present several directions for future research.

One promising avenue is to introduce Dowker complexes or weighted simplices in the network representation as proposed by~\cite{bianconi}.
Similarly to binary networks, incorporating weighted connections can describe a richer set of phenomena with simplicial complex models.
Furthermore, by carefully tuning these weights, we can control and bound the influence of large simplices to avoid the large effects that high-dimensional interactions introduce due to the combinatorial explosion.

Incorporating time-dependent information into the analysis of higher-order networks would enrich our understanding of their evolution and temporal behavior.
Exploring the dynamic aspect in the arXiv data sets opens up possibilities for detecting changes in the topology of scientific fields over time.

To gain a comprehensive understanding of network structures, we can investigate different embedding spaces to examine how the embedding space influences the topological and geometric features of the generated networks.
Related to alternative embedding spaces, the investigation of alternative connection kernels could also lead to novel network models that better describe the topological properties of higher-order networks that might be missed by traditional network representations.

\section*{Acknowledgments}
This work was supported by the Danish Data Science Academy, which is funded by the Novo Nordisk Foundation (NNF21SA0069429) and Villum Fonden (40516).
We would also like express our gratitude to T.~Owada for the careful reading of an earlier version and for helpful comments.
His suggestions helped to improve both the content and the presentation of the material.
The authors thank M.~Brun for the interesting discussions and the remark on Dowker complexes.

\bibliographystyle{plainnat}
\bibliography{literature.bib}

\end{document}